\documentclass[11pt,numbers,sort&compress]{elsarticle}
\journal{}

\makeatletter
\def\ps@pprintTitle{%
\let\@oddhead\@empty
\let\@evenhead\@empty
\def\@oddfoot{\hfill\thepage}%
\let\@evenfoot\@oddfoot}
\makeatother

\usepackage{amsmath,amssymb,amsthm,graphicx}
\usepackage{xcolor}
\usepackage{hyperref}
\AtBeginDocument{
    \hypersetup{
        colorlinks=true,
        linkcolor=black,
        citecolor=black,
        filecolor=black,
        urlcolor=black,
        pdfborder={0 0 0}
    }
}
\usepackage[labelfont=bf]{caption}
\providecommand{\doi}[1]{\href{https://doi.org/#1}{DOI:#1}}
\usepackage{xurl} 
\renewcommand{\doi}[1]{%
\href{https://doi.org/#1}{\nolinkurl{DOI:#1}}%
}

\usepackage{natbib} 
\usepackage{appendix} 
\usepackage[ruled,vlined]{algorithm2e} 
\usepackage{algpseudocode}
\usepackage{booktabs} 
\usepackage{multirow} 
\usepackage{float} 
\usepackage{enumerate} 
\usepackage{dsfont} 
\usepackage{xcolor} 
\usepackage{dsfont} 
\usepackage{geometry} 
\usepackage{tabularx} 
\usepackage{subcaption} 
\usepackage{commath, mathtools} 

\numberwithin{equation}{section} 

\theoremstyle{plain}
\newtheorem{theorem}{Theorem}[section]
\newtheorem{proposition}[theorem]{Proposition}

\newtheorem{corollary}[theorem]{Corollary}
\newtheorem{assumption}{Assumption}

\newcommand{\PP}{\mathsf{P}} 
\newcommand{\EE}{\mathsf{E}} 
\newcommand{\Bias}{\mathsf{Bias}} 
\newcommand{\Var}{\mathsf{Var}} 
\newcommand{\Cov}{\mathsf{Cov}} 
\newcommand{\ISE}{\mathsf{ISE}} 
\newcommand{\MISE}{\mathsf{MISE}} 
\newcommand{\IBias}{\mathsf{IBias}} 
\newcommand{\IVar}{\mathsf{IVar}} 
\newcommand{\OO}{\mathcal{O}}
\newcommand{\oo}{\mathrm{o}}
\newcommand{\rd}{\mathrm{d}}
\newcommand{\ind}{\mathds{1}}

\allowdisplaybreaks

\makeatletter
\newcommand{\shortname}[1]{\@eadauthor={#1}}
\makeatother

\begin{document}

\begin{frontmatter}
\title{Bernstein-smoothed estimation and bootstrap inference for the lower-tail Spearman's rho curve}

\author[a1]{Guanjie Lyu}\shortname{G.~Lyu}\ead{glyu@dal.ca}
\author[a2]{Fr\'ed\'eric Ouimet}\shortname{F.~Ouimet}\ead{frederic.ouimet.23@gmail.com}
\author[a3]{Selim O.\ Susam}\shortname{S.~O.~Susam}\ead{orhunsusam@munzur.edu.tr}

\address[a1]{Department of Community Health and Epidemiology, Faculty of Medicine, Dalhousie University, Canada}
\address[a2]{D\'epartement de math\'ematiques et d'informatique, Universit\'e du Qu\'ebec \`a Trois-Rivi\`eres, Canada}
\address[a3]{Department of Industrial Engineering, Munzur University, Tunceli, Turkey}

\begin{abstract}
This paper studies a Bernstein-smoothed plug-in estimator for the lower-tail Spearman's rho curve, a rank-based measure of local concordance defined through a normalized copula integral over the lower-left square $[0,p]^2$. The estimator applies the lower-tail Spearman functional to the empirical Bernstein copula and introduces a degree parameter that controls finite-sample regularization. The contribution is target-specific estimation and inference for the lower-tail Spearman's rho curve rather than a new general-purpose copula estimator. We establish uniform strong consistency on compact intervals away from zero whenever $m \to \infty$. Under a target-specific integrated Bernstein-bias condition and $\sqrt n/m \to 0$, we derive functional weak convergence with the same first-order Gaussian limit as the empirical copula-based estimator; this degree regime includes $m = \lfloor n^{2/3}\rfloor$. For inference, we justify a smoothed beta bootstrap based on sampling from the empirical beta copula and construct pointwise confidence intervals from the absolute centered bootstrap root. Monte Carlo experiments for the Farlie--Gumbel--Morgenstern, Gaussian, Clayton, and Frank copulas assess pointwise coverage over the threshold grid and show that Bernstein smoothing generally reduces integrated variance and often lowers mean integrated squared error under weak to moderate dependence. Additional common-sample experiments compare the proposed estimator with the empirical beta and empirical checkerboard Bernstein copula estimators in terms of pointwise error, integrated error, computation time, selected degrees, numerical stability, and pointwise coverage. A sensitivity analysis shows that $m = \lfloor n^{2/3}\rfloor$ is a simple theoretically admissible default that avoids the most severe oversmoothing. A descriptive application to the Loss--ALAE insurance claims data illustrates how the lower-tail Spearman's rho curve can reveal local dependence features that are not captured by the global Spearman correlation.
\end{abstract}


\end{frontmatter}

\thispagestyle{empty}

\section{Introduction}

Classical rank-based concordance measures, such as Spearman's rho \citep{Spearman1904} and Kendall's tau \citep{Kendall1938}, provide robust summaries of bivariate dependence. Their invariance under strictly increasing marginal transformations makes them especially useful in nonparametric dependence modeling. However, a single global coefficient may be insufficient when dependence varies across the joint distribution. In risk, insurance, environmental, and financial applications, the strength of association in a tail region can differ substantially from the overall rank association. This motivates local concordance functionals that describe dependence as a function of a tail threshold.

Copulas provide a natural framework for constructing local dependence measures because they separate marginal behavior from joint dependence. By Sklar's theorem \citep{MR125600}, any continuous bivariate distribution function $H$ with marginal distribution functions $F$ and $G$ can be written as $H(x,y) = C\{F(x),G(y)\}$, where $C : [0,1]^2 \to [0,1]$ is the unique copula associated with $H$. Consequently, rank-based concordance measures can be expressed as copula functionals and studied independently of the marginal distributions. In this framework, the lower-tail Spearman's rho introduced by \citet{MR2326243} replaces the full copula integral in Spearman's rho by a normalized integral over the lower-left square $[0,p]^2$, $p \in (0,1]$. The resulting curve $p \mapsto \rho(p)$ (see \eqref{eq:rho}) describes how concordance evolves as the lower-left region expands from the lower tail to the full unit square, coinciding with the classical Spearman correlation when $p = 1$. Thus, the lower-tail Spearman's rho curve provides a scale-dependent description of local rank dependence rather than a single global summary.

We focus on the lower-tail Spearman functional because it is linear in the copula and, after Bernstein smoothing, its defining integral can be evaluated exactly through incomplete beta functions. A natural lower-tail Kendall-type functional would instead involve a nonlinear expression analogous to integration with respect to the copula measure, and its definition, smoothing bias, and weak convergence would require a separate analysis. This distinction is one of tractability and scope, and does not imply that a local Spearman measure is uniformly preferable to a local Kendall measure.

Nonparametric estimation of $\rho(p)$ is naturally based on the empirical copula. Consider independent and identically distributed observations $(X_1,Y_1), \ldots,(X_n,Y_n)$ from a bivariate distribution function $H$ with continuous marginal distribution functions $F$ and $G$ and copula $C$, and let $F_n$ and $G_n$ denote the empirical marginal distribution functions. The empirical copula, introduced by \citet{MR573609}, is defined as
\[
C_n(u,v) = \frac{1}{n} \sum_{i=1}^n \ind\{F_n(X_i) \leq u, \mskip\thickmuskip{}G_n(Y_i) \leq v\}, \qquad (u,v) \in [0,1]^2.
\]
Substituting $C_n$ into the lower-tail Spearman functional yields the empirical copula-based estimator $\widehat{\rho}_n(p)$; see \eqref{eq:LT.Spearman}. This estimator is fully rank-based, invariant under strictly increasing marginal transformations, and does not require specification of a parametric copula family. Despite these advantages, the empirical copula-based estimator inherits the finite-sample roughness of the empirical copula. This is particularly relevant for lower-tail estimation because the effective information in the region $[0,p]^2$ decreases rapidly as $p$ becomes small. As a result, the empirical lower-tail Spearman's rho curve may exhibit substantial sampling variability at small and moderate thresholds. Moreover, because $C_n$ is a step function, the plug-in curve can be irregular as a function of $p$, which is undesirable when the goal is to estimate a smooth local dependence profile and to construct stable pointwise confidence intervals.

To improve finite-sample stability while preserving the nonparametric and rank-based nature of the empirical copula approach, we replace $C_n$ in the lower-tail Spearman functional by the empirical Bernstein copula $C_{m,n}$. This yields the Bernstein-smoothed estimator $\widehat{\rho}_{m,n}(p)$; see \eqref{eq:Bernstein.estimator}. The construction retains the main advantages of the empirical copula estimator: it is invariant under strictly increasing marginal transformations and does not require specification of a parametric copula family. At the same time, the Bernstein smoothing step regularizes the empirical copula surface before integration, which can reduce finite-sample variability in tail-focused estimation. The degree $m$ therefore acts as a smoothing parameter, balancing variance reduction against smoothing bias.

Bernstein smoothing is one of several possible regularization strategies. The empirical beta copula of \citet{Segers2017empirical} is the special empirical Bernstein copula with degree $m = n$ and is always a genuine copula; it avoids degree selection, but it does not offer the same finite-sample range of smoothing levels. B-spline copulas provide local support and flexible approximation, but typically require choices of knots, penalties, and constraints that enforce the copula margins and nonnegativity \citep{MR2427366,MR2427382,Dou2025Bspline}. By contrast, the Bernstein construction uses nonnegative basis functions, has a single scalar degree in the bivariate setting considered here, yields an exact incomplete-beta representation of the target functional, and admits a direct beta-mixture interpretation. These are practical and structural trade-offs rather than a claim of uniform superiority over other smooth copula estimators.

A particularly relevant alternative is the empirical checkerboard Bernstein copula of \citet{LuGhosh2023}. Their method applies in arbitrary dimension, produces a genuine copula for dimension-specific polynomial degrees, and selects those degrees through a hierarchical empirical Bayes procedure. It is therefore a more general method for estimating a full multivariate copula. The present paper has a narrower target: a bivariate curve-valued functional evaluated in closed form with a single degree parameter, together with frequentist process and bootstrap inference tailored to that functional. Which approach is preferable depends on whether the inferential goal is a flexible estimate of the full multivariate copula or a direct estimate of the lower-tail Spearman's rho curve, and we do not claim that the proposed method is uniformly preferable.

Bernstein smoothing has been used successfully in a range of copula-based inference procedures. For instance, \citet{MR3635017} proposed nonparametric tests of independence based on the Bernstein empirical copula and Bernstein copula density, while \citet{Lyu2023Symmetry,Lyu2023Equality} developed Bernstein-polynomial tests for copula symmetry and for equality of dependence structures across two samples, respectively. Relatedly, \citet{MR4536056} combined Bernstein smoothing with weighted Cram\'er--von Mises statistics for independence testing, showing that both the polynomial degree and the choice of weight can affect power, especially under tail-dependent alternatives. Bernstein smoothing has also been used in more targeted estimation and transform-based representation problems: \citet{doi:10.1080/00949655.2025.2494140} replaced the empirical copula by the Bernstein empirical copula in a Cram\'er--von Mises minimum-distance estimator for the dependence parameter of the Farlie--Gumbel--Morgenstern copula; \citet{MR4778424} introduced a Bernstein copula characteristic-function representation; and \citet{MR4712570} constructed a Bernstein plug-in estimator of the cross-ratio function for bivariate time-to-event data. Taken together, these contributions suggest that Bernstein smoothing can provide useful finite-sample regularization, while also highlighting that the polynomial degree acts as a smoothing parameter and therefore controls a bias--variance trade-off.

The contribution of this paper is not the introduction of a new general-purpose Bernstein copula estimator. Rather, we study the estimation and inference problem for the entire lower-tail Spearman's rho curve. We give an exact representation of the Bernstein plug-in estimator, establish consistency under the sole degree condition $m \to \infty$, derive a curve-level Gaussian limit under a target-specific integrated Bernstein-bias condition that includes the practical rule $m = \lfloor n^{2/3}\rfloor$, and justify a smoothed beta bootstrap based on the empirical beta copula. We also provide common-sample comparisons with the empirical beta copula and the empirical checkerboard Bernstein copula, including adaptive degree selection, computation time, numerical diagnostics, and a separate exploratory common-bootstrap coverage experiment using a local-grid approximation to degree reselection. The passage from the copula process to the curve is obtained by applying a bounded linear map to bounded Borel measurable sample paths, and is therefore a direct special case of the functional delta method \citep{BeutnerZahle2010}. The substantive asymptotic issues are the degree-dependent Bernstein bias and the conditional validity of the resampling scheme, rather than differentiability of the lower-tail integral itself.

The remainder of the paper is organized as follows. Section~\ref{sec:estimation} defines lower-tail Spearman's rho, reviews the empirical copula-based estimator, introduces the Bernstein-smoothed estimator, and presents the main asymptotic results. Section~\ref{sec:CI_bootstrap} develops the smoothed beta bootstrap confidence intervals. Section~\ref{sec:MC.simulations} reports the simulation study, including bootstrap performance, comparisons with the empirical copula-based estimator and alternative smooth copula estimators, and sensitivity to the Bernstein degree $m$. Section~\ref{sec:application} presents the Loss--ALAE insurance claims application. Section~\ref{sec:conclusion} concludes with a summary and directions for future research. The \textsf{R} code for the simulation study is available via the link provided at the end of the manuscript.

\section{Estimation of the lower-tail Spearman's rho function}\label{sec:estimation}

Spearman's rho is a widely used concordance measure for two continuous random variables. Expressed in terms of the underlying copula $C$, it is given by
\[
\rho_S
= 12 \int_{[0,1]^2} \{C(u,v) - \Pi(u,v)\} \, \rd u \, \rd v
= 12 \int_{[0,1]^2} C(u,v) \, \rd u \, \rd v - 3,
\]
where $\Pi(u,v) = uv$ is the independence copula; see, e.g., \citet[Section~5.1.2]{MR2197664}. It offers a rank-based alternative to Pearson's correlation coefficient that is invariant under strictly increasing transformations of the margins. However, given that $\rho_S$ integrates over the entire unit square $[0,1]^2$, it may obscure joint lower-tail behavior. To isolate that region, \citet{MR2326243} restricted integration to $[0,p]^2$, $p \in (0,1]$, and defined the following (normalized) lower-tail Spearman's rho:
\begin{equation}\label{eq:rho}
\rho(p)
= \frac{\int_{[0,p]^2} C(u,v) \, \rd u \, \rd v - \int_{[0,p]^2} \Pi(u,v) \, \rd u \, \rd v}{\int_{[0,p]^2} M(u,v) \, \rd u \, \rd v - \int_{[0,p]^2} \Pi(u,v) \, \rd u \, \rd v}
= \frac{\int_{[0,p]^2} C(u,v) \, \rd u \, \rd v - p^4/4}{p^3/3 - p^4/4},
\end{equation}
where $M(u,v) = \min(u,v)$ is the Fr\'echet--Hoeffding upper bound copula. This lower-tail version of $\rho_S$ retains the interpretation of concordance while focusing exclusively on joint behavior within $[0,p]^2$.

Set
\[
D(p) = \int_{[0,p]^2} \{M(u,v) - \Pi(u,v)\} \, \rd u \, \rd v = \frac{p^3}{3} - \frac{p^4}{4}, \qquad p \in (0,1].
\]
Since $D(p) > 0$, this quantity provides the natural normalization for the restricted concordance functional. It calibrates the measure so that $\rho(p) = 0$ under independence ($C = \Pi$) and $\rho(p) = 1$ under perfect positive dependence ($C = M$). Hence $\rho(p)$ retains the usual calibration between independence and comonotonicity, while replacing the integration over the full unit square in Spearman's rho by integration over the lower-left region $[0,p]^2$.

\subsection{The Bernstein-smoothed empirical-copula estimator}\label{sec:Bernstein.EC}

Following \citet{MR2326243}, the empirical-copula estimator of $\rho(p)$ is obtained by substituting $C_n$ for $C$ in the lower-tail functional. Thus, for $p \in (0,1]$, define
\begin{equation}\label{eq:LT.Spearman}
\begin{aligned}
\widehat{\rho}_n(p)
= \frac{ \int_{[0,p]^2} C_n(u,v) \, \rd u \, \rd v - \int_{[0,p]^2} \Pi(u,v) \, \rd u \, \rd v }{ D(p) }
= \frac{ \int_{[0,p]^2} C_n(u,v) \, \rd u \, \rd v - p^4/4 }{ D(p) }.
\end{aligned}
\end{equation}

Instead of integrating the empirical copula $C_n$ directly, we first apply the bivariate Bernstein operator. For $k = 0, \ldots,m$, let
\[
P_{k,m}(u) = \binom{m}{k}u^k(1 - u)^{m-k}, \qquad u \in [0,1],
\]
and, for any bounded function $h$ on $[0,1]^2$, define
\[
B_m(h)(u,v) = \sum_{k=0}^m \sum_{\ell=0}^m h(k/m,\ell/m)P_{k,m}(u)P_{\ell,m}(v).
\]
The empirical Bernstein copula is $C_{m,n} = B_m(C_n)$. Although standard terminology refers to $C_{m,n}$ as a copula, it need not be a copula in the finite-sample sense for arbitrary $m$: in the absence of ties and for equal Bernstein degrees, \citet[Proposition~2.5]{Segers2017empirical} show that it is a copula if and only if $m$ divides $n$. This distinction is immaterial for the integral functional below. The special case $m = n$ is the empirical beta copula, denoted by $C_n^\beta = C_{n,n}$, which is always a genuine copula \citep{Segers2017empirical}.

For any $p \in (0,1]$, the Bernstein-smoothed lower-tail Spearman's rho estimator is defined as
\begin{equation}\label{eq:Bernstein.estimator}
\widehat{\rho}_{m,n}(p) = \frac{\int_{[0,p]^2} C_{m,n}(u,v) \, \rd u \, \rd v - p^4/4}{D(p)}.
\end{equation}
For numerical evaluation, note that
\[
\int_{[0,p]^2} C_{m,n}(u,v) \, \rd u \, \rd v = \sum_{k=0}^m \sum_{\ell=0}^m C_n(k/m,\ell/m) \binom{m}{k}\binom{m}{\ell} B_p(k + 1,m - k + 1) \, B_p(\ell + 1,m - \ell + 1),
\]
where, for $a,b > 0$ and $p \in [0,1]$,
\[
B_p(a,b) = \int_0^p t^{a-1}(1 - t)^{b-1} \, \rd t
\]
denotes the lower incomplete beta function.

\subsection{Asymptotic properties}\label{sec:asymptotic}

This subsection establishes the large-sample properties of the proposed Bernstein-smoothed estimator. Let $\mathcal{P} = [p_0,1]$, where $p_0 \in (0,1)$ is fixed. The restriction $p \geq p_0$ avoids the degeneracy of the normalizing constant $D(p)$ as $p \downarrow 0$.

\begin{theorem}[Strong consistency]\label{thm:strong.consistency}
Assume $m = m(n) \to \infty$. Then, as $n \to \infty$,
\[
\sup_{p \in \mathcal{P}}|\widehat{\rho}_{m,n}(p) - \rho(p)| \to 0, \qquad \text{a.s.}
\]
In particular, $\widehat{\rho}_{m,n}(p) \to \rho(p)$ almost surely for every fixed $p \in (0,1]$.
\end{theorem}

Theorem~\ref{thm:strong.consistency} requires only that the Bernstein degree diverges and therefore covers the practical choice $m = \lfloor n^{2/3}\rfloor$. We next derive the first-order distribution of the estimated curve. The following smoothness condition is standard for rank-based empirical copula processes with unknown margins.

\begin{assumption}
\label{ass:segers}
The first-order partial derivatives
\[
\dot{C}_1(u,v) = \frac{\partial C(u,v)}{\partial u}, \qquad \dot{C}_2(u,v) = \frac{\partial C(u,v)}{\partial v}
\]
exist and are continuous, respectively, on $(0,1)\times[0,1]$ and $[0,1]\times(0,1)$.
\end{assumption}

For a compact set $T$, let $\ell^\infty(T)$ denote the Banach space of bounded real-valued functions on $T$, equipped with the supremum norm
\[
\|h\|_{\infty,T} = \sup_{t \in T}|h(t)|.
\]
Weak convergence in this space is denoted by $ \rightsquigarrow $. Under Assumption~\ref{ass:segers}, the empirical copula process
\[
\mathbb{C}_n = \sqrt n(C_n - C)
\]
converges weakly in $\ell^\infty([0,1]^2)$ to the centered Gaussian process
\[
\mathbb{G}_C(u,v) = \mathbb{B}_C(u,v) - \dot{C}_1(u,v)\mathbb{B}_C(u,1) - \dot{C}_2(u,v)\mathbb{B}_C(1,v), \qquad (u,v) \in [0,1]^2,
\]
where $\mathbb{B}_C$ is the Brownian bridge with covariance function
\[
K_C((u,v),(s,t)) = C(\min(u,s),\min(v,t)) - C(u,v)C(s,t).
\]
At boundary points where the partial derivatives are not defined by Assumption~\ref{ass:segers}, they may be assigned arbitrary values because the corresponding Brownian bridge terms vanish. Write
\[
\Gamma((u,v),(s,t)) = \Cov\{\mathbb{G}_C(u,v),\mathbb{G}_C(s,t)\}.
\]

The relevant smoothing bias for the present curve is an integrated bias rather than the uniform approximation error of the full copula surface.

\begin{assumption}
\label{ass:integrated.bias}
As $r \to \infty$,
\[
\sup_{p \in \mathcal{P}}\left|\int_{[0,p]^2}\{B_r(C)(u,v) - C(u,v)\} \, \rd u \, \rd v\right| = \OO(r^{-1}).
\]
\end{assumption}

Assumption~\ref{ass:integrated.bias} is target-specific and is weaker than requiring $\|B_r(C) - C\|_{\infty,[0,1]^2} = \OO(r^{-1})$. The following proposition gives a simple sufficient condition. The integrated condition may also hold under less restrictive local or weighted smoothness conditions.

\begin{proposition}
\label{prop:bias.sufficient}
If $C$ is continuously differentiable on $[0,1]^2$ and its gradient is globally Lipschitz, then Assumption~\ref{ass:integrated.bias} holds.
\end{proposition}

\begin{theorem}[Weak convergence of the curve]
\label{thm:functional_clt_rho}
Suppose that Assumptions~\ref{ass:segers} and~\ref{ass:integrated.bias} hold, $m = m(n) \to \infty$, and $\sqrt n/m \to 0$. Then
\[
\mathbb{Z}_n(p) = \sqrt n\{\widehat{\rho}_{m,n}(p) - \rho(p)\}, \qquad p \in \mathcal{P},
\]
converges weakly in $\ell^\infty(\mathcal{P})$ to the centered Gaussian process
\[
\mathbb{Z}(p) = \frac{1}{D(p)} \int_{[0,p]^2}\mathbb{G}_C(u,v) \, \rd u \, \rd v, \qquad p \in \mathcal{P}.
\]
Moreover, for any $p,q \in \mathcal{P}$,
\[
\Cov\{\mathbb{Z}(p),\mathbb{Z}(q)\} = \frac{1}{D(p)D(q)} \int_{[0,p]^2} \int_{[0,q]^2} \Gamma((u,v),(s,t)) \, \rd u \, \rd v \, \rd s \, \rd t.
\]
\end{theorem}

The degree rule $m = \lfloor n^{2/3}\rfloor$ satisfies both $m \to \infty$ and $\sqrt n/m \to 0$. Hence, under Assumption~\ref{ass:integrated.bias}, the degree used in the numerical study is covered by Theorems~\ref{thm:strong.consistency} and~\ref{thm:functional_clt_rho}. The proof of Theorem~\ref{thm:functional_clt_rho} separates the stochastic smoothing error from the deterministic integrated Bernstein bias. Once those terms are controlled, the passage from the empirical copula process to the lower-tail Spearman's rho curve is immediate because the defining map is bounded and linear on the bounded Borel measurable functions equipped with the supremum norm. Thus, this mapping step is a direct continuous-mapping or functional delta-method argument rather than a separate differentiability difficulty.

\section{Confidence intervals via the smoothed beta bootstrap}\label{sec:CI_bootstrap}

The functional weak convergence result in Section~\ref{sec:asymptotic} provides the first-order distributional basis for $\widehat{\rho}_{m,n}(p)$. For each fixed $p$, the limiting variance depends on the unknown copula through the covariance kernel $\Gamma$. Rather than estimating this kernel directly, we use a smoothed beta bootstrap based on the empirical beta copula; see \citet{Kiriliouk2021resampling,Kojadinovic2024resampling}. The empirical beta copula is a genuine copula and is straightforward to simulate from, which avoids the ties created by resampling the observed rank pairs directly. Throughout this section, $\PP^\#$ denotes probability conditional on the original sample.

The procedure is implemented on a fixed grid
\[
\mathcal{P}_J = \{p_1, \ldots,p_J\} \subseteq (0,1], \qquad p_1 < \cdots < p_J.
\]
Let $R_i^X$ and $R_i^Y$ denote the marginal ranks of $X_i$ and $Y_i$, respectively. For each $i = 1, \ldots,n$, let $Q_i^\beta$ be the product probability law on $[0,1]^2$ under which
\[
U^\# \sim \mathrm{Beta}\{R_i^X,n + 1 - R_i^X\}, \qquad V^\# \sim \mathrm{Beta}\{R_i^Y,n + 1 - R_i^Y\},
\]
with $U^\#$ and $V^\#$ independent. The equally weighted mixture $n^{-1}\sum_{i=1}^n Q_i^\beta$ has distribution function $C_n^\beta = C_{n,n}$. For the $b$th bootstrap replication, draw $I_1^{\#,b}, \ldots,I_n^{\#,b}$ independently from the uniform distribution on $\{1, \ldots,n\}$ and, conditionally on these indices, generate independently
\[
\mathbf W_i^{\#,b} = (U_i^{\#,b},V_i^{\#,b}) \sim Q_{I_i^{\#,b}}^\beta, \qquad i = 1, \ldots,n.
\]
The resulting bootstrap sample is denoted by
\[
\mathcal{W}_b^\# = (\mathbf W_1^{\#,b}, \ldots,\mathbf W_n^{\#,b}), \qquad b = 1, \ldots,B.
\]

A more general product-beta generator can be formed by letting $A_i^X = \lceil m_\#R_i^X/n\rceil$ and $A_i^Y = \lceil m_\#R_i^Y/n\rceil$ and, conditionally on mixture component $i$, drawing independently
\[
U^\# \sim \mathrm{Beta}\{A_i^X,m_\# + 1 - A_i^X\}, \qquad V^\# \sim \mathrm{Beta}\{A_i^Y,m_\# + 1 - A_i^Y\}.
\]
Its mixture distribution is $C_{m_\#,n}$, which need not have uniform margins unless $m_\#$ divides $n$. For this reason, the validity result and the main bootstrap algorithm below use the empirical beta choice $m_\# = n$.

\begin{algorithm}[!ht]
\caption{Smoothed beta bootstrap pointwise confidence intervals}
\label{alg:bootstrap.ci}
\SetKwInOut{Input}{Input}
\SetKwInOut{Output}{Output}

\Input{Data $(X_i,Y_i)_{i=1}^n$; estimator degree $m$; grid $\mathcal{P}_J = \{p_1, \ldots,p_J\}$; number of bootstrap replications $B$; significance level $\alpha$.}
\Output{Pointwise intervals $\widehat I_{n,\alpha}^\#(p_j)$, $j = 1, \ldots,J$.}

Compute the original estimates
\[
\widehat{\rho}_j = \widehat{\rho}_{m,n}(p_j), \qquad j = 1, \ldots,J.
\]

\For{$b = 1, \ldots,B$}{
Generate the bootstrap sample
\[
\mathcal{W}_b^\# = (\mathbf W_1^{\#,b}, \ldots,\mathbf W_n^{\#,b})
\]
from the empirical beta copula using the product-beta scheme described above.

From $\mathcal{W}_b^\#$, compute the rank-based bootstrap empirical copula $C_n^{\#,b}$ and its Bernstein smoothing
\[
C_{m,n}^{\#,b}(s,t) = \sum_{k=0}^m \sum_{\ell=0}^m C_n^{\#,b}(k/m,\ell/m)P_{k,m}(s)P_{\ell,m}(t).
\]

For $j = 1, \ldots,J$, compute
\[
\widehat{\rho}_{j}^{\#,b} = \frac{\int_{[0,p_j]^2} C_{m,n}^{\#,b}(s,t) \, \rd s \, \rd t - p_j^4/4}{D(p_j)}, \qquad T_b^\#(p_j) = \sqrt n\{\widehat{\rho}_j^{\#,b} - \widehat{\rho}_j\}.
\]
}

\For{$j = 1, \ldots,J$}{
Let $c_{j,1-\alpha}^\#$ be the empirical $(1 - \alpha)$-quantile of $|T_1^\#(p_j)|, \ldots,|T_B^\#(p_j)|$. Define
\[
\widehat I_{n,\alpha}^\#(p_j) = \left[\widehat{\rho}_j - \frac{c_{j,1-\alpha}^\#}{\sqrt n}, \, \widehat{\rho}_j + \frac{c_{j,1-\alpha}^\#}{\sqrt n}\right].
\]
}
\end{algorithm}

For a single bootstrap replication, let $\widehat{\rho}_{m,n}^\#(p)$ denote the estimator computed from the bootstrap ranks and define
\[
\mathbb{Z}_n^\#(p) = \sqrt n\{\widehat{\rho}_{m,n}^\#(p) - \widehat{\rho}_{m,n}(p)\}, \qquad p \in \mathcal{P}.
\]
The following result gives the conditional justification for the bootstrap procedure.

\begin{theorem}[Bootstrap validity]
\label{thm:bootstrap.validity}
Suppose that Assumptions~\ref{ass:segers} and~\ref{ass:integrated.bias} hold, $m = m(n) \to \infty$, and $\sqrt n/m \to 0$. Conditionally on the data, in probability,
\[
\mathbb{Z}_n^\# \rightsquigarrow \mathbb{Z} \qquad \text{in }\ell^\infty(\mathcal{P}),
\]
where $\mathbb{Z}$ is the Gaussian process in Theorem~\ref{thm:functional_clt_rho}.
\end{theorem}

\begin{corollary}[Pointwise coverage]
\label{cor:pointwise.coverage}
Under the conditions of Theorem~\ref{thm:bootstrap.validity}, fix $p_j \in \mathcal{P}$ and $\alpha \in (0,1)$, and suppose that $\Var\{\mathbb{Z}(p_j)\} > 0$. Let $c_{j,1-\alpha,n}^\#$ denote the conditional $(1 - \alpha)$-quantile of $|\mathbb{Z}_n^\#(p_j)|$. Then
\[
\PP\left\{\rho(p_j) \in \left[\widehat{\rho}_{m,n}(p_j) - \frac{c_{j,1-\alpha,n}^\#}{\sqrt n}, \, \widehat{\rho}_{m,n}(p_j) + \frac{c_{j,1-\alpha,n}^\#}{\sqrt n}\right]\right\} \to 1 - \alpha.
\]
If $B = B_n \to \infty$ as $n \to \infty$, the empirical bootstrap quantile in Algorithm~\ref{alg:bootstrap.ci} may replace $c_{j,1-\alpha,n}^\#$ without changing the coverage limit.
\end{corollary}

We next explain the form of the interval in Algorithm~\ref{alg:bootstrap.ci}. At a fixed grid point $p_j$, write
\[
T_n(p_j) = \sqrt n\{\widehat{\rho}_{m,n}(p_j) - \rho(p_j)\}.
\]
If $c_{j,1-\alpha}^\#$ consistently estimates the $(1 - \alpha)$-quantile of $|T_n(p_j)|$, then
\[
\PP\{|T_n(p_j)| \leq c_{j,1-\alpha}^\#\} \approx 1 - \alpha.
\]
Rearranging this event gives
\[
\rho(p_j) \in \left[\widehat{\rho}_{m,n}(p_j) - \frac{c_{j,1-\alpha}^\#}{\sqrt n}, \, \widehat{\rho}_{m,n}(p_j) + \frac{c_{j,1-\alpha}^\#}{\sqrt n}\right].
\]
Thus, the symmetry is induced by the absolute centered bootstrap root and by the centered Gaussian first-order approximation; the interval is not a bootstrap percentile interval. Symmetric absolute-root intervals are standard bootstrap constructions; see, for example, \citet{Hall1988Symmetric} and the copula applications in \citet{Kiriliouk2021resampling}. If the signed bootstrap root is markedly asymmetric, an equal-tailed basic interval can instead be formed from its empirical quantiles $q_{j,\alpha/2}^\#$ and $q_{j,1-\alpha/2}^\#$ as
\[
\left[\widehat{\rho}_j - \frac{q_{j,1-\alpha/2}^\#}{\sqrt n}, \, \widehat{\rho}_j - \frac{q_{j,\alpha/2}^\#}{\sqrt n}\right].
\]
For comparison, let $r_{j,\gamma}^\#$ denote the empirical $\gamma$-quantile of $\widehat{\rho}_{j}^{\#,1}, \ldots,\widehat{\rho}_{j}^{\#,B}$. The corresponding percentile interval would be
\[
\left[r_{j,\alpha/2}^\#, \, r_{j,1-\alpha/2}^\#\right].
\]
The numerical study compares the symmetric absolute-root and equal-tailed basic intervals using the same bootstrap samples. The intervals considered here are pointwise in $p$ and should not be interpreted as simultaneous confidence bands for the entire curve.

\section{Simulation study}\label{sec:MC.simulations}

This section provides a computational assessment of the proposed Bernstein-smoothed lower-tail Spearman's rho estimator and the associated smooth bootstrap confidence intervals. The degree condition in Theorem~\ref{thm:functional_clt_rho} allows $m = \lfloor n^{2/3}\rfloor$ under Assumption~\ref{ass:integrated.bias}, but it does not identify a finite-sample optimal degree. In numerical implementation, $m$ is therefore treated as a smoothing parameter: it should be large enough to avoid excessive smoothing bias, but not so large that the variance-reduction benefit of Bernstein smoothing is lost. Motivated by the pointwise bias--variance expansion for the Bernstein copula estimator in \cite{MR2879763}, we use $m = \lfloor n^{2/3}\rfloor$ as a transparent practical default and examine its stability through the sensitivity analysis in Section~\ref{sec:degree}.

The numerical study evaluates the finite-sample operating characteristics of the proposed procedure. We first examine the accuracy of the empirical-beta bootstrap by reporting pointwise coverage over a grid of lower-tail thresholds, together with Monte Carlo error and a comparison of the symmetric absolute-root and equal-tailed basic intervals. We then study the bias--variance trade-off introduced by Bernstein smoothing by decomposing the integrated error into integrated squared bias, integrated variance, and MISE, and by comparing these quantities with those of the empirical copula-based estimator under the default choice $m = \lfloor n^{2/3}\rfloor$. A common-target comparison further evaluates the empirical beta copula and the empirical checkerboard Bernstein copula of \citet{LuGhosh2023} on the same samples, including adaptive degrees, computation time, convergence diagnostics, numerical failures, and a separate exploratory common-bootstrap coverage experiment. Finally, we investigate the sensitivity of the MISE to the smoothing degree through local perturbations of the default rule, in order to assess whether the proposed degree choice provides stable finite-sample performance across copula families, dependence levels, and sample sizes.

\subsection{Finite-sample assessment of bootstrap confidence intervals}

All simulations in this subsection are based on the Farlie--Gumbel--Morgenstern (FGM) copula~\citep[see, e.g.,][p.~77]{MR2197664}, defined by
\[
C_\theta(u,v) = uv\{1 + \theta(1 - u)(1 - v)\}, \qquad (u,v) \in [0,1]^2, \qquad \theta \in [-1,1].
\]
This model provides a convenient benchmark because the target lower-tail Spearman's rho curve is available in closed form. Indeed,
\[
\int_{[0,p]^2} C_\theta(u,v) \, \rd u \, \rd v = \frac{p^4}{4} + \theta \left( \frac{p^2}{2} - \frac{p^3}{3} \right)^2,
\]
and hence
\[
\rho_\theta(p) = \frac{ \theta \left( p^2/2 - p^3/3 \right)^2 }{ p^3/3 - p^4/4 }, \qquad p \in (0,1].
\]

We examine the smooth bootstrap confidence intervals. The simulation uses
\[
\theta \in \{-1, -0.5, 0, 0.5, 1\}, \qquad n \in \{50,200\},
\]
and the threshold grid
\[
\mathcal{P}_J = \{0.05, 0.10, \ldots, 1.00\}, \qquad J = 20.
\]
For each pair $(\theta,n)$, we generate $K = 500$ independent Monte Carlo samples from the FGM copula. For each sample, we compute $\widehat{\rho}_{m,n}(p_j)$ over the full grid using $m = \lfloor n^{2/3}\rfloor$. The bootstrap generator is the empirical beta copula, so $m_\# = n$, and every generated sample is reranked before the estimator is recomputed. Pointwise confidence intervals are based on $B = 399$ bootstrap replications and nominal level $1 - \alpha = 0.95$. The symmetric absolute-root and equal-tailed basic intervals are computed from the same bootstrap samples. These finite-$B$ intervals are Monte Carlo approximations to the ideal-bootstrap intervals in Corollary~\ref{cor:pointwise.coverage}.

For each grid point $p_j$, empirical pointwise coverage is estimated by
\[
\widehat{\mathrm{CP}}(p_j) = \frac{1}{K} \sum_{r=1}^K \ind\mskip-\thinmuskip{}\big\{\rho_\theta(p_j) \in \widehat I_{n,\alpha}^{\#,(r)}(p_j)\big\},
\]
where $\widehat I_{n,\alpha}^{\#,(r)}(p_j)$ denotes the bootstrap confidence interval obtained in the $r$th Monte Carlo replication. Pointwise coverage as a function of $p_j$ is the primary diagnostic. Its plug-in Monte Carlo standard error is estimated by $\{\widehat{\mathrm{CP}}(p_j)[1 - \widehat{\mathrm{CP}}(p_j)]/K\}^{1/2}$, while the figures display 95\% Wilson score intervals for the underlying coverage probability. Wilson intervals are used because, unlike untruncated Wald intervals, they remain informative when the estimated coverage is near zero or one. The values at $p = 0.05$ and $p = 0.10$ are reported separately because the effective sample size is smallest near the lower boundary. Coverage averaged over the full threshold grid is retained only as a secondary summary.

\begin{figure}[!ht]
\centering
\includegraphics[trim = 2mm 0mm 4mm 3mm, clip, width = 0.95\textwidth]{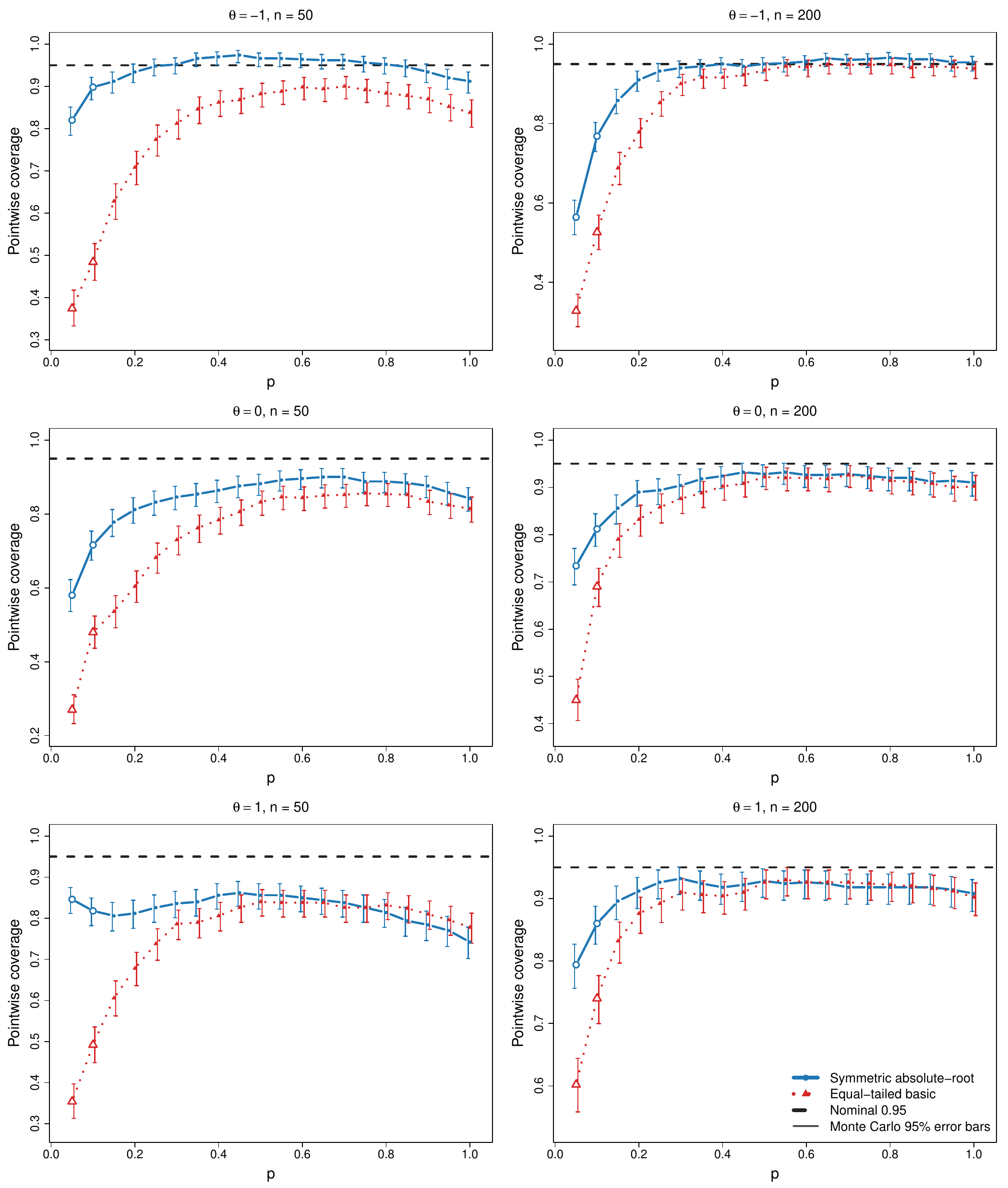}
\caption{Empirical pointwise coverage of nominal 95\% confidence intervals under the FGM copula. Rows correspond to $\theta = -1,0,1$, and columns correspond to $n = 50$ and $n = 200$. The blue solid curves represent the symmetric absolute-root intervals, the red dotted curves represent the equal-tailed basic intervals, and the horizontal dashed line gives the nominal coverage level. Vertical bars are 95\% Wilson Monte Carlo intervals. The values at $p = 0.05$ and $p = 0.10$ are highlighted. The estimator degree is $m = \lfloor n^{2/3}\rfloor$, and bootstrap samples are generated from the empirical beta copula with $m_\# = n$.}
\label{fig:fgm_pointwise_coverage}
\end{figure}

Figure~\ref{fig:fgm_pointwise_coverage} displays the pointwise coverage curves for $\theta \in \{-1,0,1\}$. The two smallest thresholds are highlighted, and the symmetric absolute-root and equal-tailed basic intervals are shown together so that the effect of bootstrap-root asymmetry can be assessed directly. Table~\ref{tab:fgm_tail_coverage} reports the corresponding coverage values and plug-in Monte Carlo standard errors at $p = 0.05$ and $p = 0.10$, together with the grid averages. The complete pointwise results for $\theta = -0.5$ and $\theta = 0.5$ are given in Figure~\ref{fig:fgm_pointwise_coverage_app} in~\ref{app:simulation}.

\newpage
\begin{table}[!ht]
\centering
\small
\setlength{\tabcolsep}{5pt}
\renewcommand{\arraystretch}{0.9}
\caption{Empirical pointwise coverage at the two smallest thresholds under the FGM copula. Plug-in binomial Monte Carlo standard errors are shown in parentheses. The final two columns report coverage averaged over the full threshold grid.}
\label{tab:fgm_tail_coverage}
\begin{tabular}{crrrrrr}
\toprule[1.5pt]
$\theta$ & $n$ & $p$ & Symmetric & Basic & Mean symmetric & Mean basic \\
\midrule[1pt]
$-1$ & $50$ & $0.05$ & 0.820 (0.017) & 0.374 (0.022) & 0.941 & 0.802 \\
$-1$ & $50$ & $0.10$ & 0.898 (0.014) & 0.484 (0.022) & 0.941 & 0.802 \\
$-1$ & $200$ & $0.05$ & 0.564 (0.022) & 0.328 (0.021) & 0.918 & 0.860 \\
$-1$ & $200$ & $0.10$ & 0.768 (0.019) & 0.526 (0.022) & 0.918 & 0.860 \\
\midrule
$-0.5$ & $50$ & $0.05$ & 0.540 (0.022) & 0.316 (0.021) & 0.881 & 0.768 \\
$-0.5$ & $50$ & $0.10$ & 0.792 (0.018) & 0.390 (0.022) & 0.881 & 0.768 \\
$-0.5$ & $200$ & $0.05$ & 0.606 (0.022) & 0.440 (0.022) & 0.913 & 0.871 \\
$-0.5$ & $200$ & $0.10$ & 0.764 (0.019) & 0.600 (0.022) & 0.913 & 0.871 \\
\midrule
$0$ & $50$ & $0.05$ & 0.580 (0.022) & 0.270 (0.020) & 0.843 & 0.746 \\
$0$ & $50$ & $0.10$ & 0.716 (0.020) & 0.480 (0.022) & 0.843 & 0.746 \\
$0$ & $200$ & $0.05$ & 0.734 (0.020) & 0.450 (0.022) & 0.900 & 0.863 \\
$0$ & $200$ & $0.10$ & 0.812 (0.017) & 0.690 (0.021) & 0.900 & 0.863 \\
\midrule
$0.5$ & $50$ & $0.05$ & 0.768 (0.019) & 0.298 (0.020) & 0.864 & 0.776 \\
$0.5$ & $50$ & $0.10$ & 0.820 (0.017) & 0.486 (0.022) & 0.864 & 0.776 \\
$0.5$ & $200$ & $0.05$ & 0.754 (0.019) & 0.544 (0.022) & 0.925 & 0.892 \\
$0.5$ & $200$ & $0.10$ & 0.840 (0.016) & 0.744 (0.020) & 0.925 & 0.892 \\
\midrule
$1$ & $50$ & $0.05$ & 0.846 (0.016) & 0.354 (0.021) & 0.824 & 0.756 \\
$1$ & $50$ & $0.10$ & 0.818 (0.017) & 0.492 (0.022) & 0.824 & 0.756 \\
$1$ & $200$ & $0.05$ & 0.794 (0.018) & 0.602 (0.022) & 0.910 & 0.885 \\
$1$ & $200$ & $0.10$ & 0.860 (0.016) & 0.740 (0.020) & 0.910 & 0.885 \\
\bottomrule[1.5pt]
\end{tabular}
\end{table}

Table~\ref{tab:fgm_tail_coverage} also shows that neither interval is uniformly well calibrated near the lower boundary. At $p = 0.05$, symmetric absolute-root coverage is as low as $0.540$ in the reported design, while equal-tailed basic coverage is as low as $0.270$; these discrepancies are far larger than the displayed Monte Carlo errors. The symmetric interval generally performs better than the basic interval at the two smallest thresholds, but it too can have severe undercoverage. The deterioration is consistent with finite-sample bias and the limited effective information available when only a small fraction of observations contributes to the lower-left region. Consequently, a grid-average coverage close to the nominal level must not be read as adequate calibration of the full curve, since averaging can conceal poor coverage at the thresholds where tail inference is most difficult. The pointwise values and Monte Carlo intervals, rather than the grid average, are therefore the primary basis for assessing the intervals.

\begin{figure}[!ht]
\centering
\includegraphics[trim = 3mm 2mm 7mm 6mm, clip, width = 0.48\textwidth]{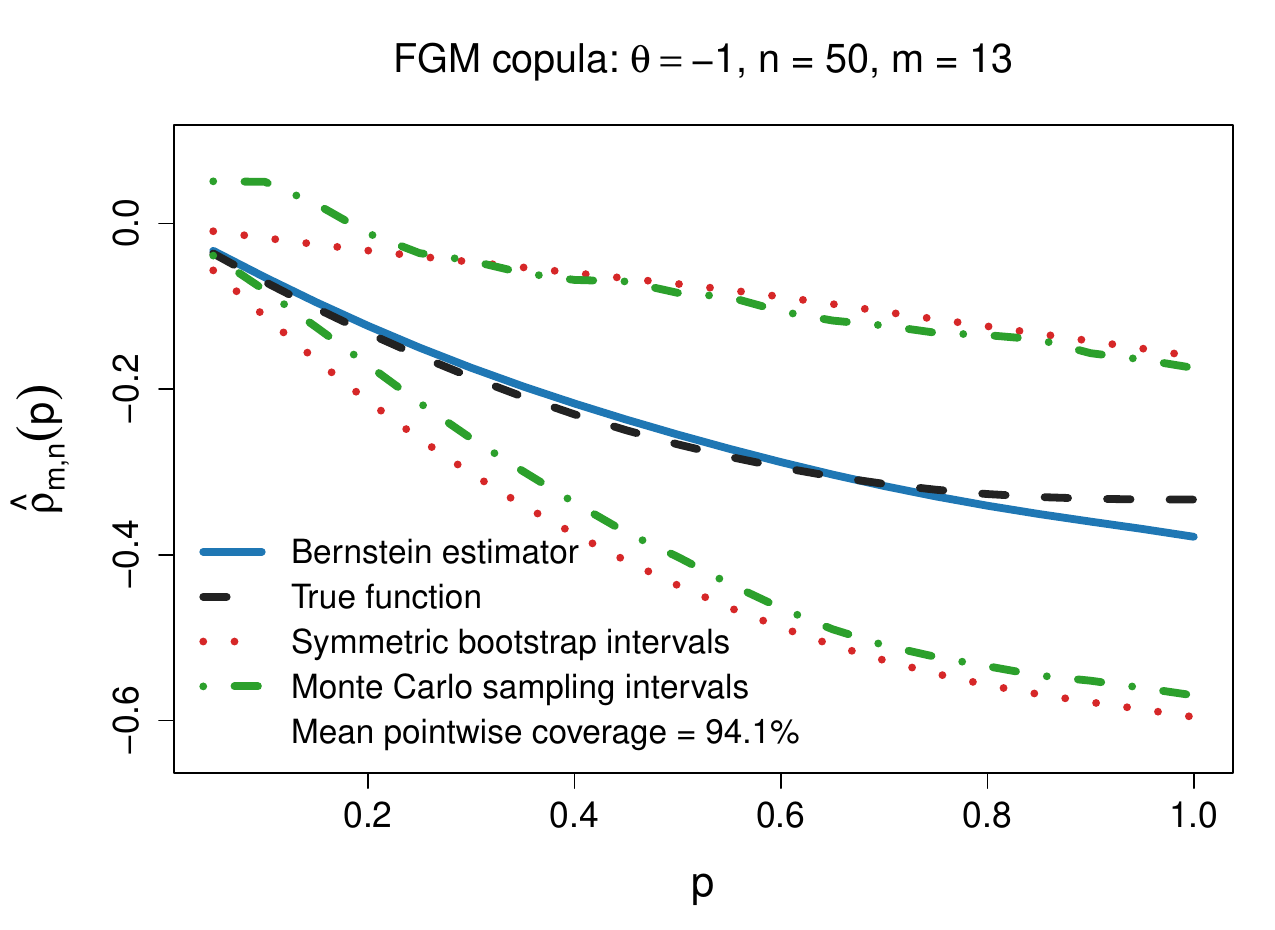}
\quad
\includegraphics[trim = 3mm 2mm 7mm 6mm, clip, width = 0.48\textwidth]{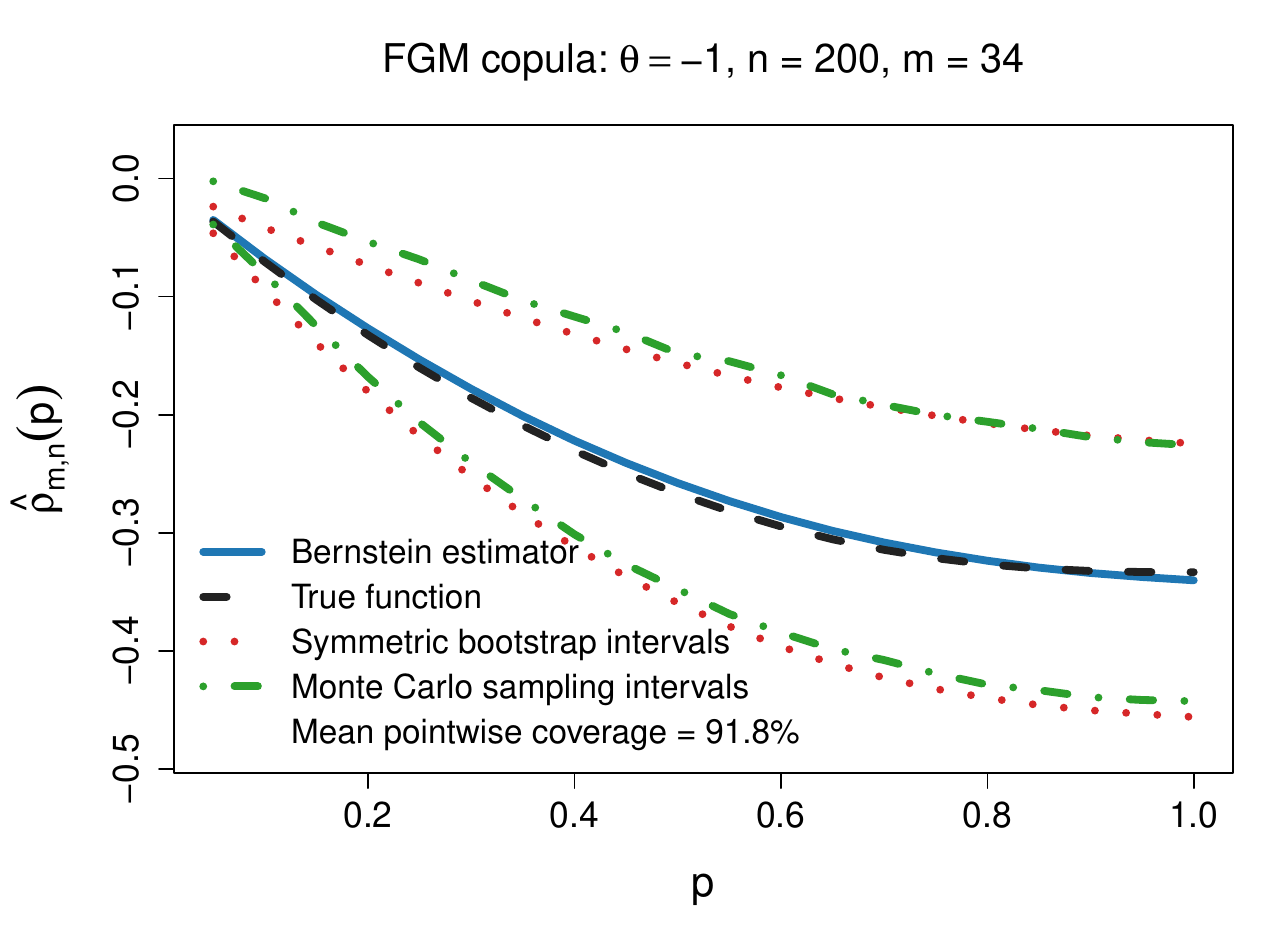} \\
\includegraphics[trim = 3mm 2mm 7mm 6mm, clip, width = 0.48\textwidth]{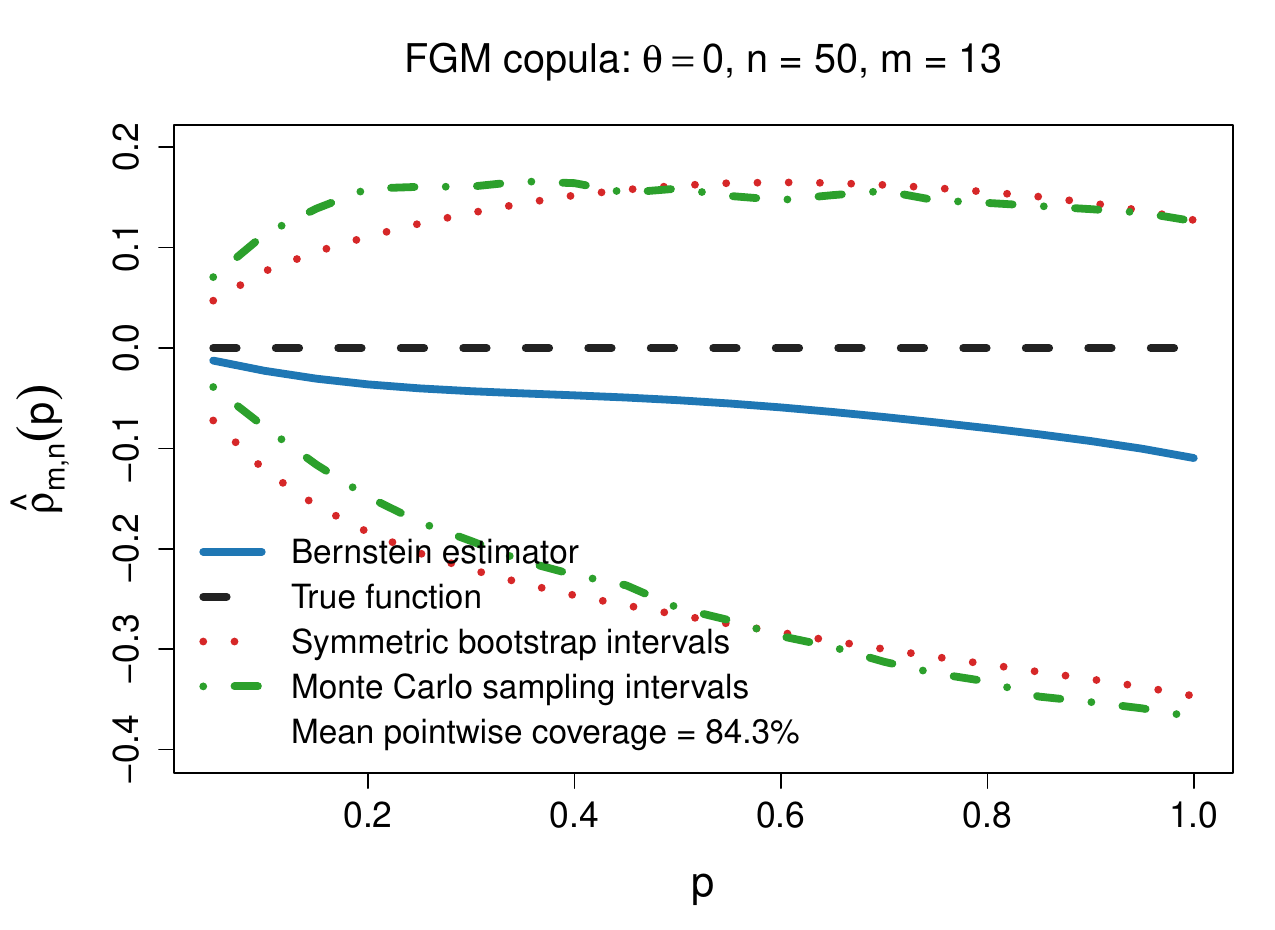}
\quad
\includegraphics[trim = 3mm 2mm 7mm 6mm, clip, width = 0.48\textwidth]{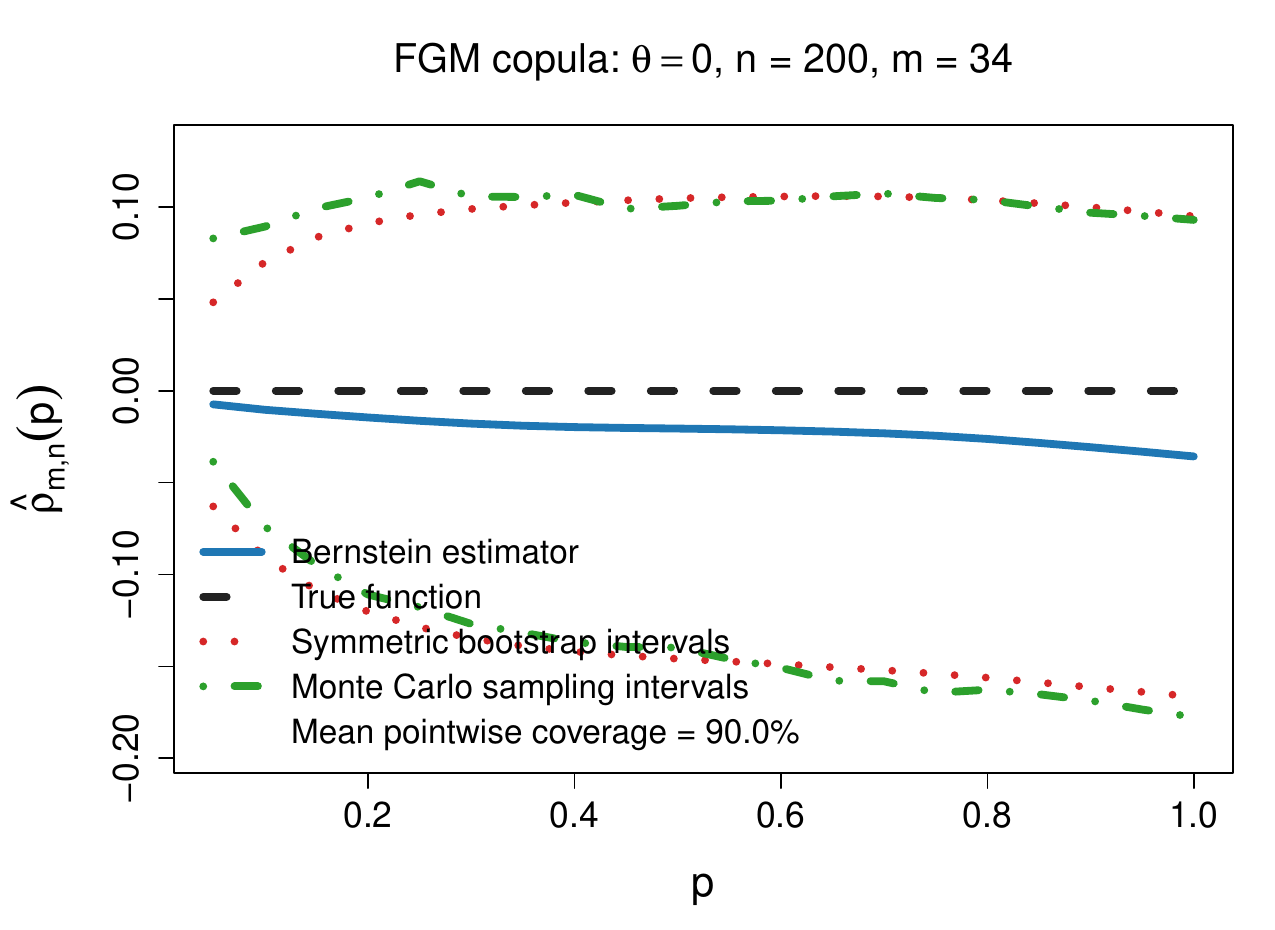} \\
\includegraphics[trim = 3mm 2mm 7mm 6mm, clip, width = 0.48\textwidth]{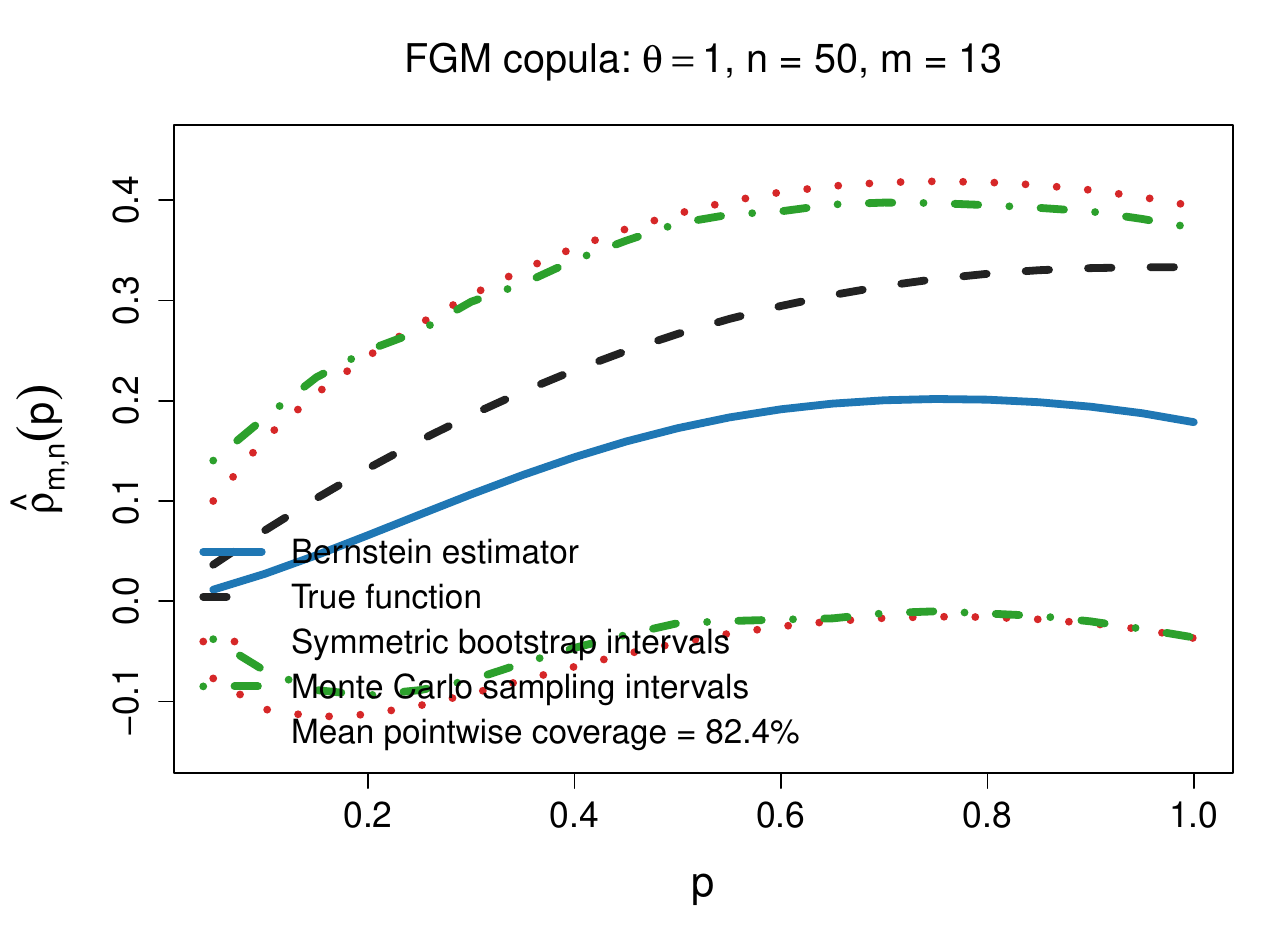}
\quad
\includegraphics[trim = 3mm 2mm 7mm 6mm, clip, width = 0.48\textwidth]{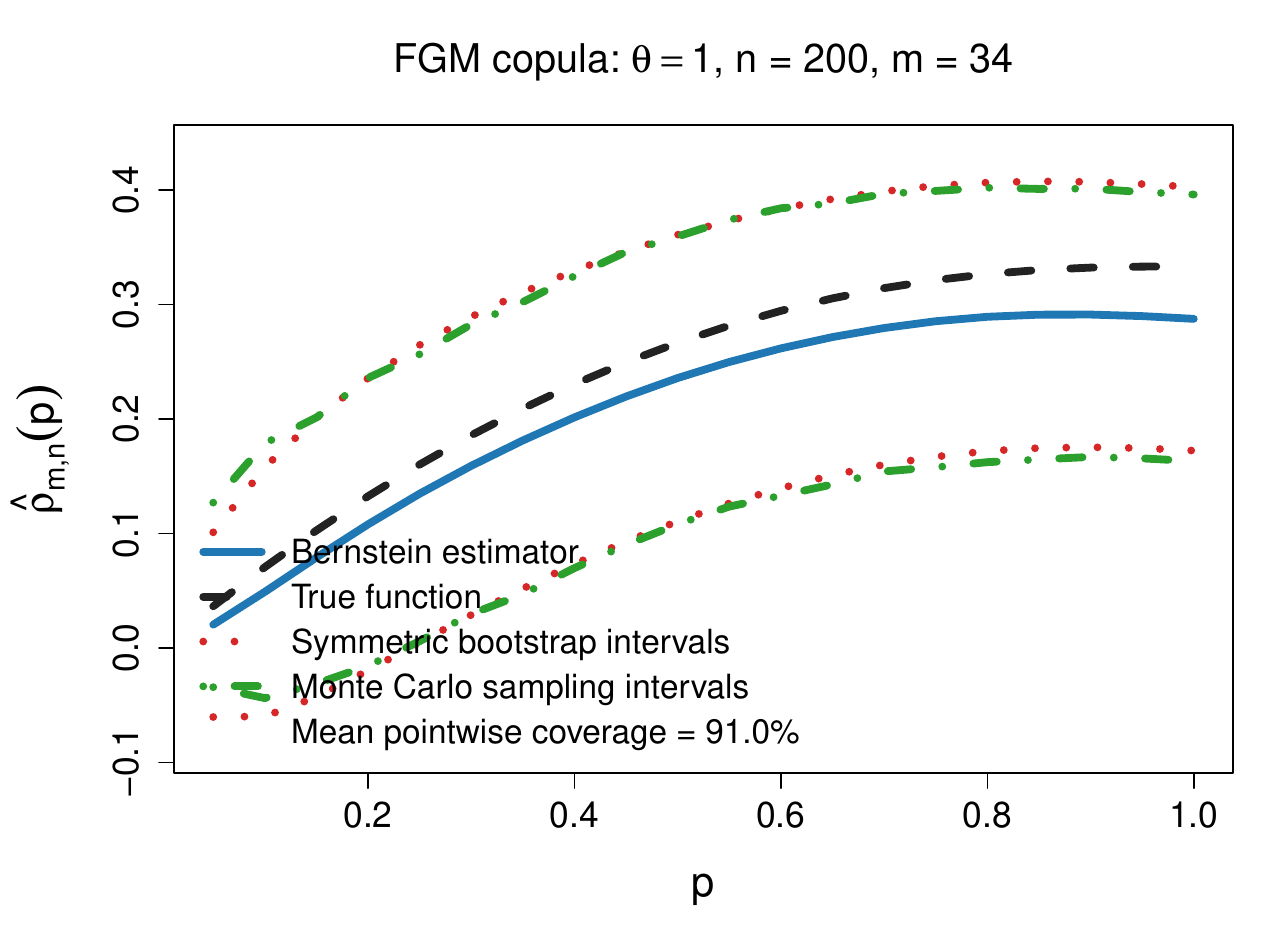} \\[-3mm]
\caption{Finite-sample assessment of the empirical-beta bootstrap pointwise confidence intervals under the FGM copula. Each panel displays the Monte Carlo average of the Bernstein-smoothed lower-tail Spearman's rho estimator $p \mapsto \widehat{\rho}_{m,n}(p)$ (solid blue), the true curve $p \mapsto \rho_\theta(p)$ (dashed black), the average lower and upper endpoints of the symmetric absolute-root intervals (red dotted), and Monte Carlo sampling intervals obtained from the empirical pointwise quantiles of $\widehat{\rho}_{m,n}(p)$ across Monte Carlo replications (green dash-dotted). Rows correspond to $\theta = -1,0,1$, while columns correspond to $n = 50$ and $n = 200$, with $m = \lfloor n^{2/3}\rfloor$ and $m_\# = n$. The mean pointwise coverage reported in each panel is the empirical coverage averaged over the threshold grid and is included as a secondary summary.}
\label{fig:simulation_CI}
\end{figure}

Figure~\ref{fig:simulation_CI} complements the coverage analysis by comparing the Monte Carlo average of the estimator and the average bootstrap endpoints with the true curve and the Monte Carlo sampling intervals. This separates interval centering and width from the binary coverage indicator. The mean pointwise coverage shown in each panel is the secondary grid average reported above. The corresponding panels for $\theta = -0.5$ and $\theta = 0.5$ are reported in Figure~\ref{fig:simulation_CI_app} in~\ref{app:simulation}.

\subsection{Comparison with the empirical copula-based estimator}\label{sec:MISE_comparison}

To further examine the finite-sample behavior of the Bernstein estimator, we compare it with the empirical copula-based estimator in a controlled Monte Carlo experiment. The purpose of this comparison is not to show uniform dominance of one estimator over the other, but rather to identify the dependence regimes in which Bernstein smoothing improves estimation accuracy and those in which the additional smoothing bias becomes non-negligible. We focus on estimation of the lower-tail Spearman's rho curve, which summarizes the strength of concordance in the lower-left quadrant of the copula as the truncation level varies. This comparison is particularly relevant because the empirical copula estimator is a natural rank-based benchmark, whereas the Bernstein estimator introduces smoothing and therefore reflects an explicit bias--variance trade-off. In the simulations below, estimation accuracy is evaluated over the interval $p \in [0.05,1]$.

\paragraph{Simulation setting}

We considered four copula families representing distinct dependence constructions: the Farlie--Gumbel--Morgenstern (FGM) copula, which has a simple polynomial form and allows only bounded dependence; the Gaussian copula, representing the elliptical class; and the Clayton and Frank copulas, representing the Archimedean class. For the FGM copula, the dependence parameter was varied over $\theta \in \{-1, -0.5, 0, 0.5, 1\}$. This choice was made because the FGM copula has a restricted Kendall's tau range. For the Gaussian, Clayton and Frank copulas, dependence was indexed by Kendall's tau, $\tau \in \{-0.5, -0.25, 0, 0.25, 0.5\}$. Two sample sizes were considered, $n = 50$ and $n = 200$, representing small and moderate samples. For each combination of copula family, dependence level, and sample size, we generated $K = 500$ independent Monte Carlo samples from the specified copula.

\paragraph{Performance measures}

Let $E$ denote the empirical copula-based estimator and $B$ denote the Bernstein estimator. For $A \in \{E, B\}$, write $\widehat{\rho}_{A}^{(k)}(p)$ for the estimate obtained from estimator $A$ in the $k$th Monte Carlo replication, $k = 1, \ldots,K$. For each replication, estimation accuracy was first summarized by the integrated squared error (ISE)
\[
\ISE^{(k)}_A = \int_{0.05}^{1} \big\{ \widehat{\rho}^{(k)}_{A}(p) - \rho(p) \big\}^2 \, \rd p, \qquad A \in \{E, B\}.
\]
In computation, this integral was approximated by the trapezoidal rule over the grid of $p$ values.

The mean integrated squared error was then estimated by averaging the ISE values over the $K$ Monte Carlo replications:
\[
\MISE_A = \frac{1}{K} \sum_{k=1}^{K} \ISE^{(k)}_A, \qquad A \in \{E, B\}.
\]
For estimator $A \in \{E, B\}$, define the pointwise Monte Carlo mean, bias, and variance, by
\[
\overline{\rho}_A(p) = \frac{1}{K} \sum_{k=1}^K \widehat{\rho}_A^{(k)}(p), \qquad
\Bias_A(p) = \overline{\rho}_A(p) - \rho(p), \qquad
\Var_A(p) = \frac{1}{K} \sum_{k=1}^K \{\widehat{\rho}_A^{(k)}(p) - \overline{\rho}_A(p)\}^2.
\]
To summarize the bias--variance decomposition over the full range of truncation levels, we integrate the pointwise squared bias and variance with respect to $p$ over $[0.05,1]$ and define
\[
\IBias^2_A = \int_{0.05}^{1} \Bias^2_A(p) \, \rd p, \qquad \IVar_A = \int_{0.05}^{1} \Var_A(p) \, \rd p.
\]
Accordingly,
\[
\MISE_A = \IBias^2_A + \IVar_A.
\]
Finally, to compare the two estimators directly, we report the MISE ratio as $\MISE_{B}/\MISE_{E} \times 100\%$. Values of the ratio below 100\% indicate that the Bernstein estimator has smaller MISE than the empirical copula-based estimator.

\begin{figure}[!ht]
\centering
\includegraphics[trim = 3mm 2mm 7mm 6mm, clip, width = 0.41\textwidth]{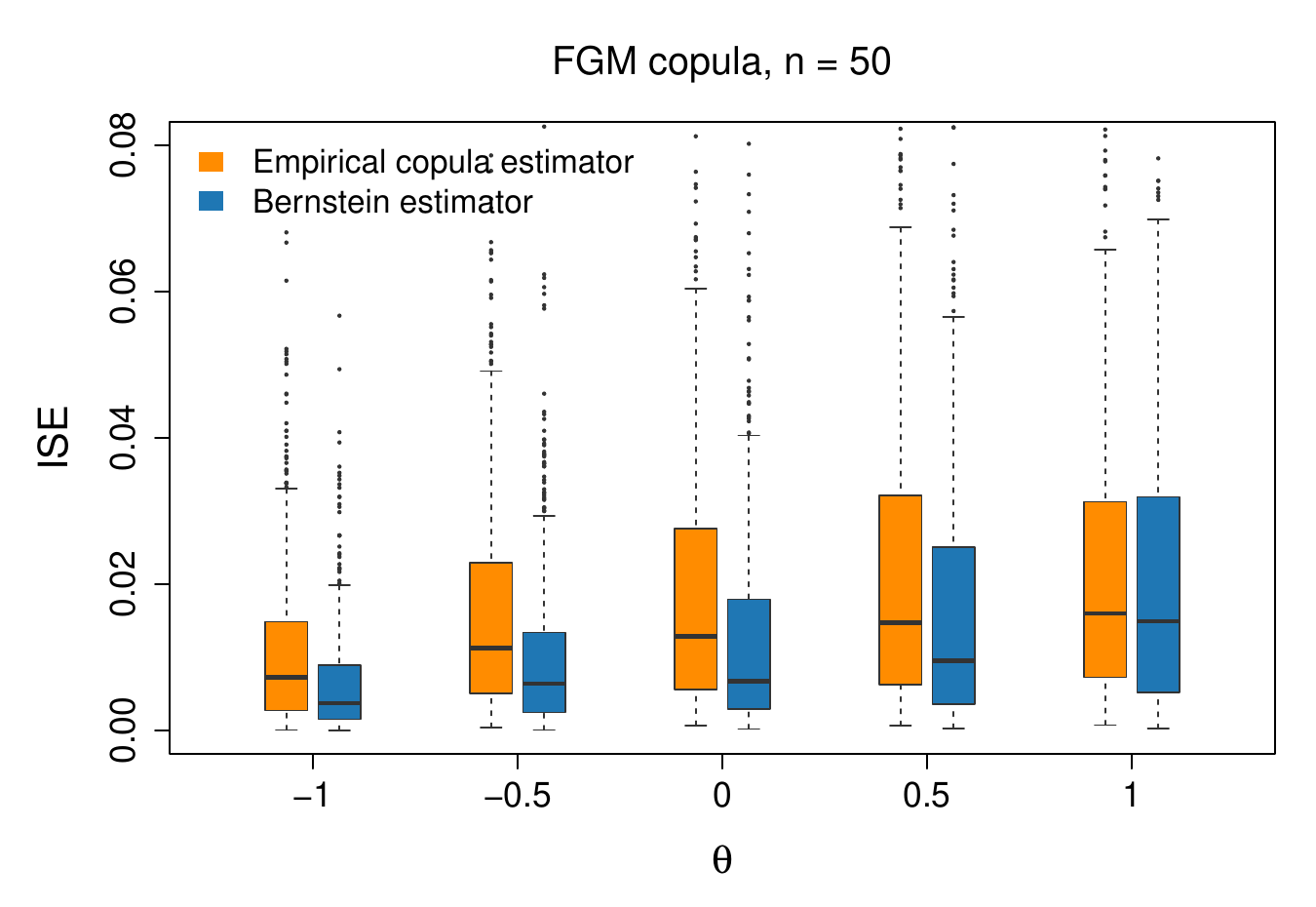}
\quad
\includegraphics[trim = 3mm 2mm 7mm 6mm, clip, width = 0.41\textwidth]{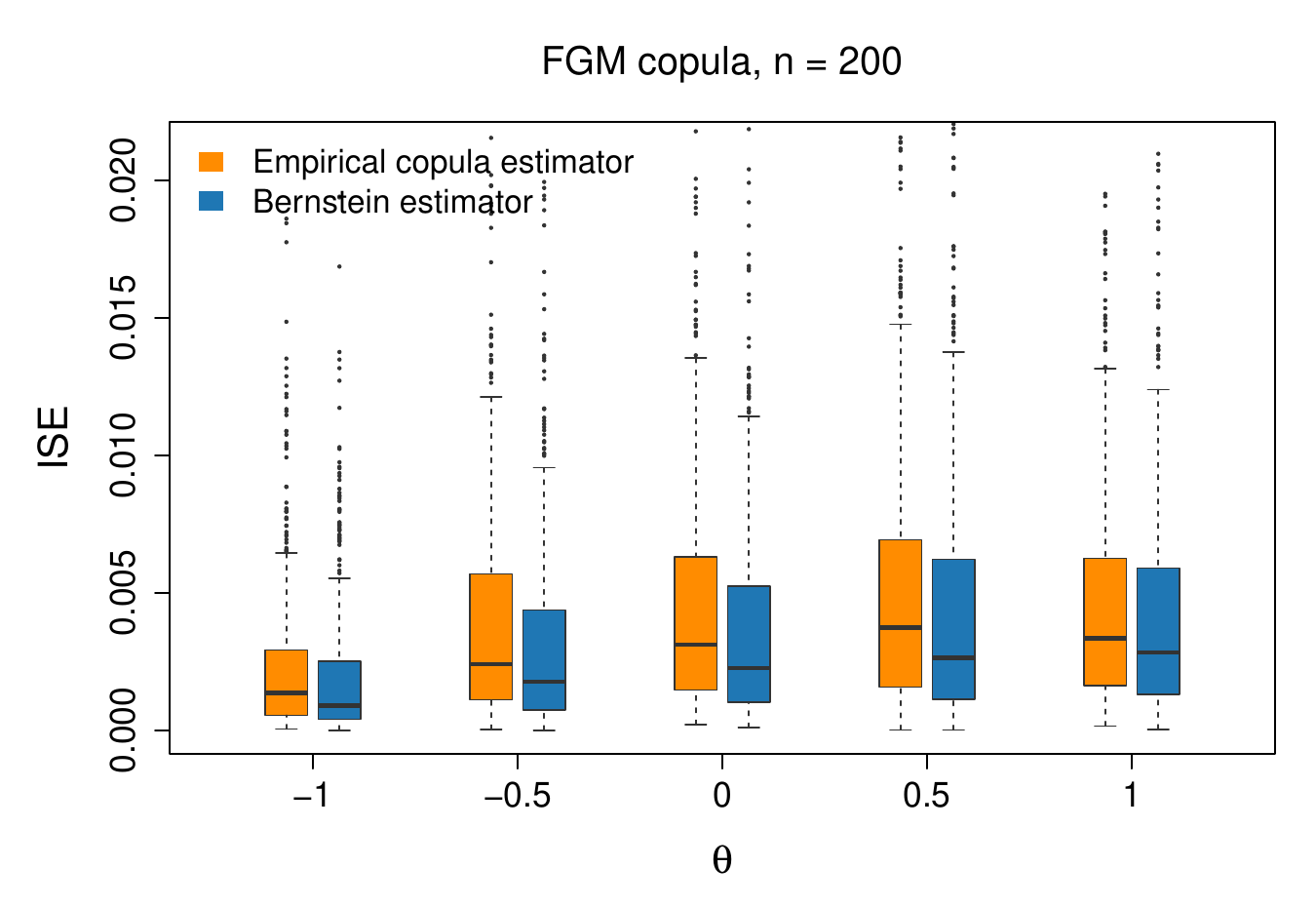} \\
\includegraphics[trim = 3mm 2mm 7mm 6mm, clip, width = 0.41\textwidth]{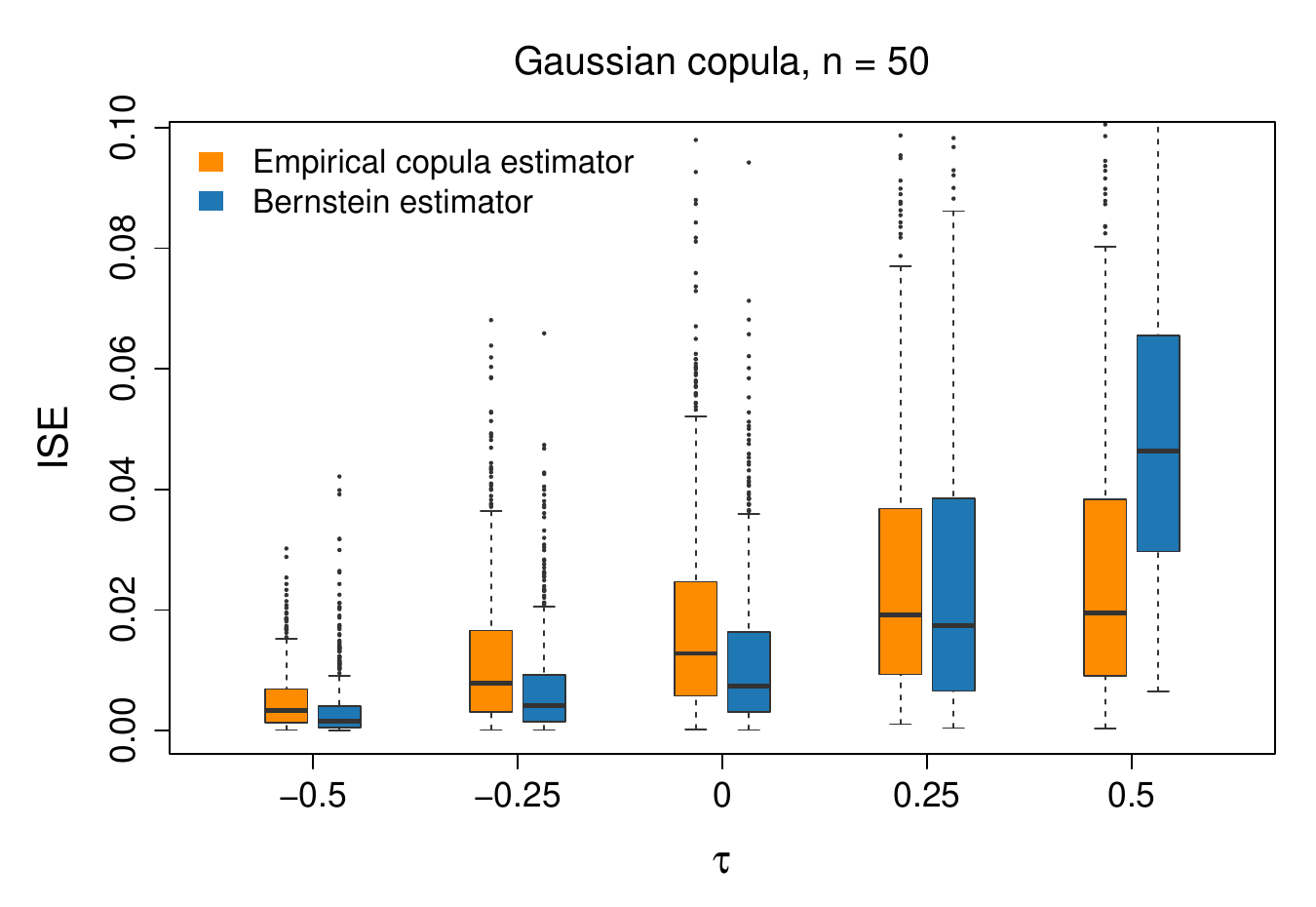}
\quad
\includegraphics[trim = 3mm 2mm 7mm 6mm, clip, width = 0.41\textwidth]{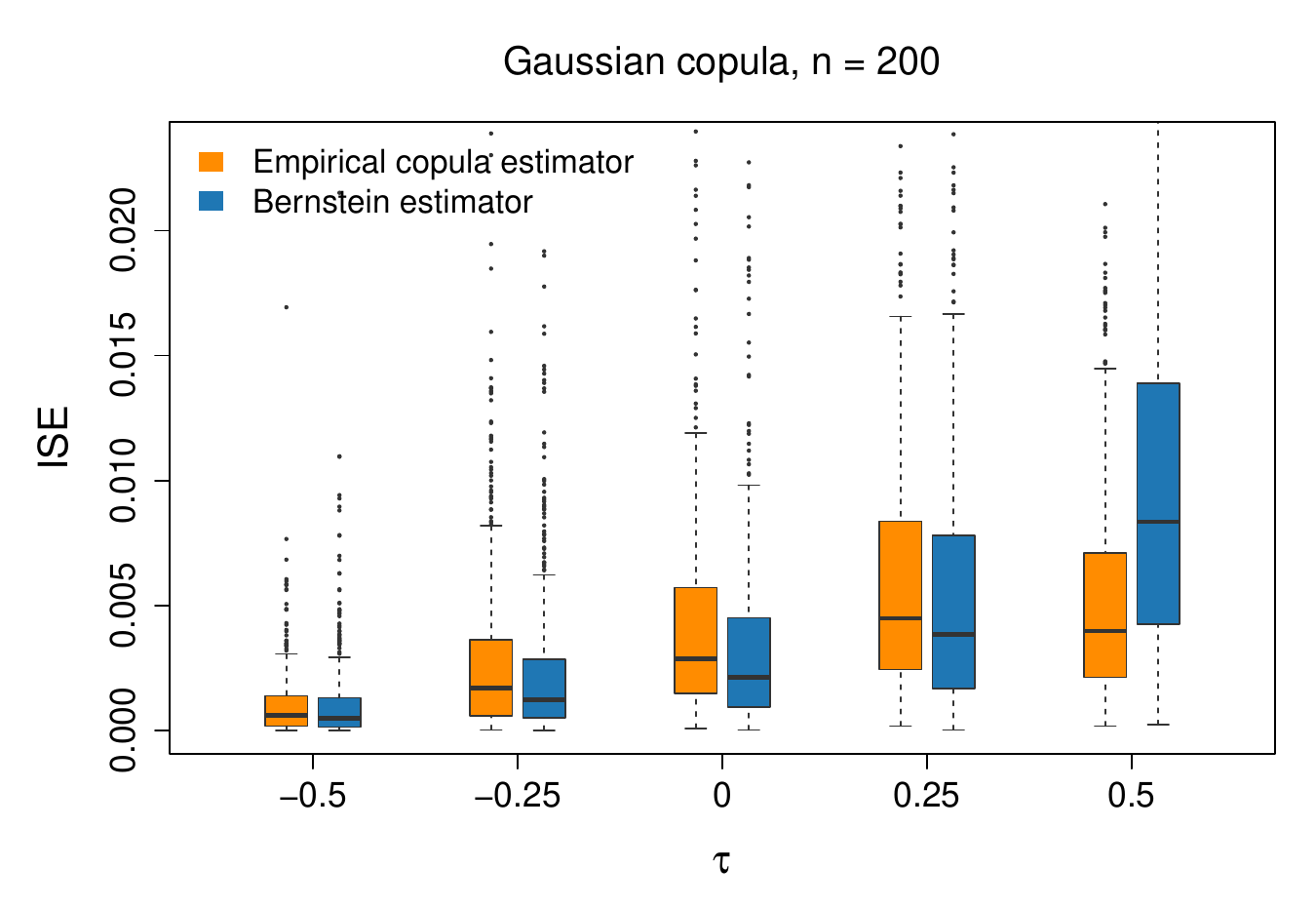} \\[-3mm]
\caption{Finite-sample comparison of the empirical copula-based estimator (orange) and the Bernstein estimator (blue) for the lower-tail Spearman's rho. Boxplots show the Monte Carlo distribution of the integrated squared error for the Farlie--Gumbel--Morgenstern (FGM) and Gaussian copulas. The horizontal axis gives the copula parameter $\theta$ for the FGM copula and Kendall's $\tau$ for the Gaussian copula. For readability, each panel's vertical range is determined by the 98th percentile of its pooled finite ISE values; the boxplots and numerical summaries are calculated from all finite evaluations.}
\label{fig:sim_boxplot}
\end{figure}

\begin{figure}[H]
\centering
\includegraphics[trim = 3mm 2mm 7mm 6mm, clip, width = 0.41\textwidth]{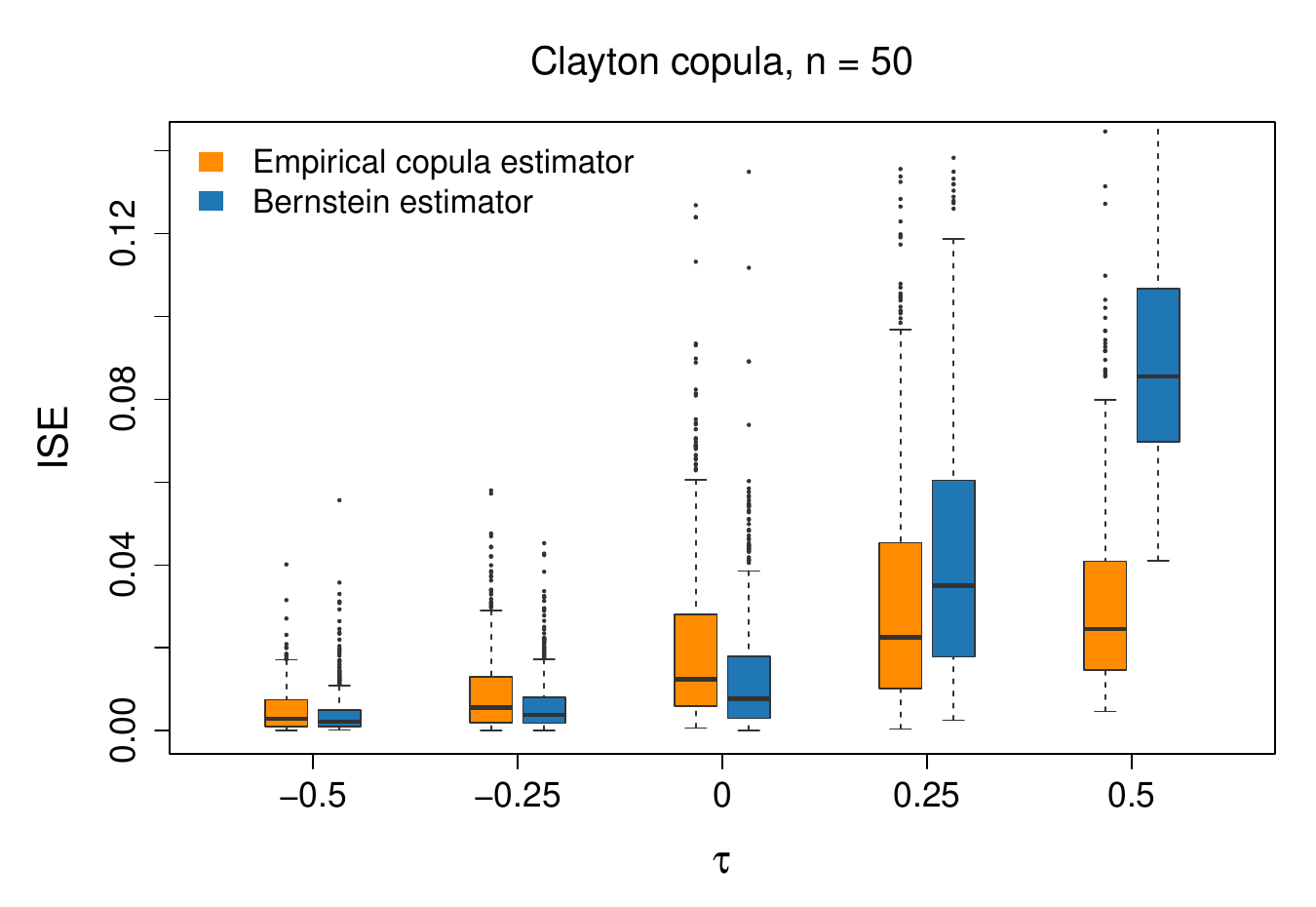}
\quad
\includegraphics[trim = 3mm 2mm 7mm 6mm, clip, width = 0.41\textwidth]{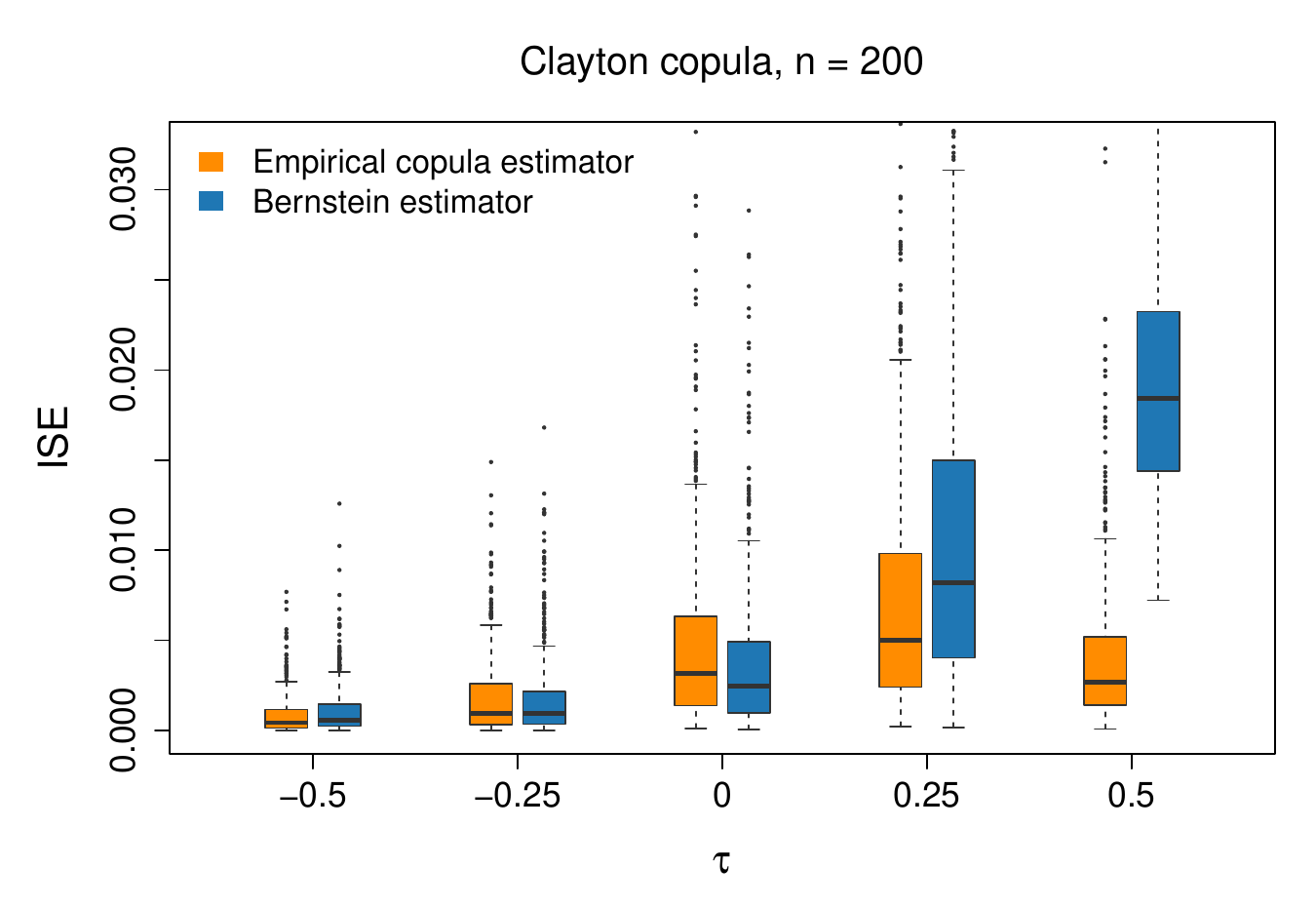} \\
\includegraphics[trim = 3mm 2mm 7mm 6mm, clip, width = 0.41\textwidth]{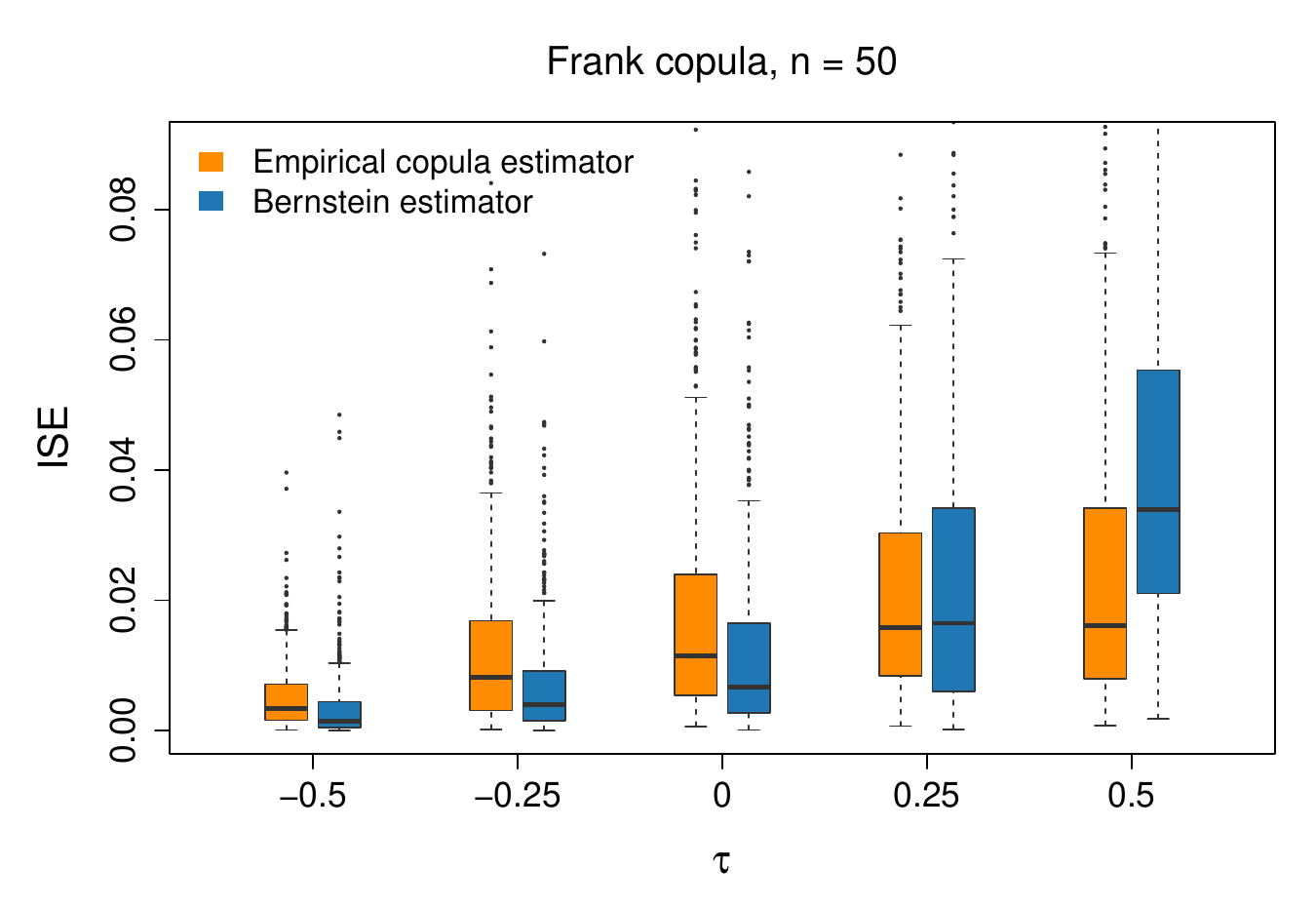}
\quad
\includegraphics[trim = 3mm 2mm 7mm 6mm, clip, width = 0.41\textwidth]{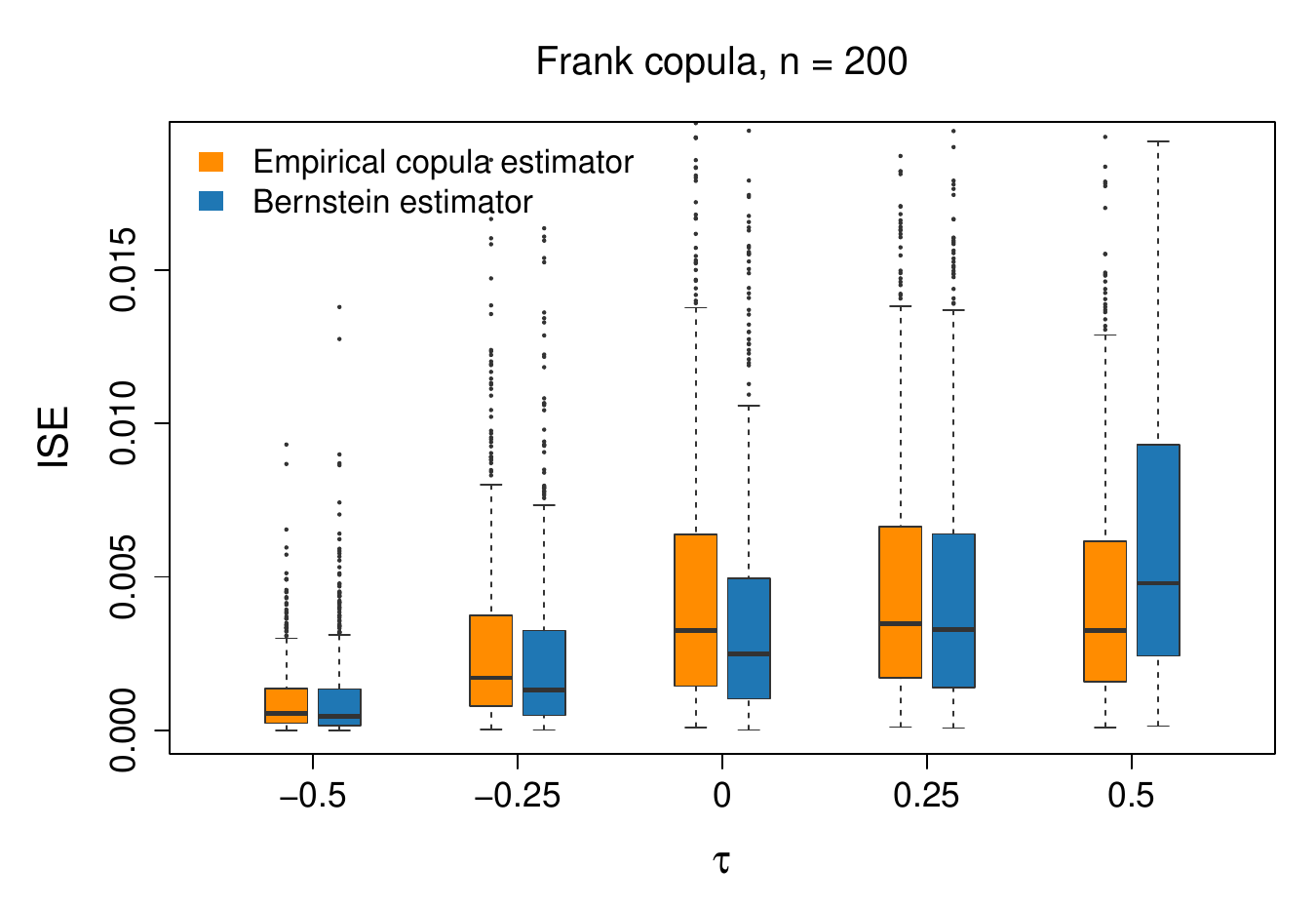} \\[-3mm]
\caption{Finite-sample comparison of the empirical copula-based estimator (orange) and the Bernstein estimator (blue) for the lower-tail Spearman's rho. Boxplots show the Monte Carlo distribution of the integrated squared error for the Clayton and Frank copulas. The horizontal axis gives Kendall's $\tau$ for the Clayton and Frank copulas. For readability, each panel's vertical range is determined by the 98th percentile of its pooled finite ISE values; the boxplots and numerical summaries are calculated from all finite evaluations.}
\label{fig:sim_boxplot_2}
\end{figure}

\paragraph{Finite-sample performance} Figures~\ref{fig:sim_boxplot}~and~\ref{fig:sim_boxplot_2} show that both estimators improve substantially as the sample size increases from $n = 50$ to $n = 200$, with noticeably smaller ISE values and tighter Monte Carlo distributions in the right column. For weak to moderate dependence, the Bernstein estimator generally yields lower and more stable ISE distributions than the empirical copula-based estimator. This pattern is most apparent for the FGM copula and for the Gaussian and Frank copulas at $\tau = -0.25,0,$ and $0.25$. However, the improvement is not uniform. For stronger positive dependence, particularly at $\tau = 0.5$ under the Gaussian, Clayton and Frank copulas, the Bernstein estimator exhibits larger ISE values than the empirical copula-based estimator. This indicates that the smoothing introduced by the Bernstein estimator can become less favorable when the copula surface has stronger curvature associated with higher dependence.

The MISE decompositions in Tables~\ref{tab:mise_decomp_n50}--\ref{tab:mise_decomp_n200} of~\ref{app:simulation} clarify this behavior. Across all reported settings, the Bernstein estimator reduces the integrated variance relative to the empirical copula-based estimator. This variance reduction is the expected benefit of smoothing. At the same time, Bernstein smoothing can increase the integrated squared bias, especially when dependence is strong. For $n = 200$, additional MISE reversals occur at $\tau = -0.5$ for the Gaussian, Clayton, and Frank copulas. Thus, the relative performance of the two estimators is governed by a bias--variance trade-off: Bernstein smoothing is advantageous when the variance reduction dominates the added bias, but it can be inferior when the smoothing bias becomes too large.

\subsection{Comparison with alternative smooth copula estimators}\label{sec:alternative.comparison}

To place the target-specific estimator relative to other smooth copula procedures, we compare the empirical copula-based estimator, the fixed-degree Bernstein estimator, the empirical beta copula estimator of \citet{Segers2017empirical}, and the empirical checkerboard Bernstein copula (ECBC) estimator of \citet{LuGhosh2023}. Here, fixed-degree refers to the deterministic rule $m = \lfloor n^{2/3}\rfloor$. The empirical beta estimator is the special case
\[
\widehat{\rho}_{\beta,n}(p) = \frac{\int_{[0,p]^2} C_n^\beta(u,v) \, \rd u \, \rd v - p^4/4}{D(p)}, \qquad p \in (0,1],
\]
where $C_n^\beta = C_{n,n}$. For ECBC, let $\widehat C_n^{\mathrm{ECBC}}$ denote the fitted bivariate copula with dimension-specific degrees $(\widehat m_1,\widehat m_2)$ selected by the empirical Bayes procedure of \citet{LuGhosh2023}, and define
\[
\widehat{\rho}_{\mathrm{ECBC},n}(p) = \frac{\int_{[0,p]^2} \widehat C_n^{\mathrm{ECBC}}(u,v) \, \rd u \, \rd v - p^4/4}{D(p)}, \qquad p \in (0,1].
\]
Thus, all four estimators are compared after application of the same lower-tail Spearman functional, rather than through a full-surface loss that may favor a method designed for a different target.

The point-estimation comparison uses the same copula families, dependence levels, sample sizes, and threshold grid as in Section~\ref{sec:MISE_comparison}. Because ECBC degree selection is substantially more computationally intensive, $K = 100$ common Monte Carlo samples are used for each family, dependence level, and sample size. These samples are the first 100 seeded replications from the larger MISE design. The fixed-degree Bernstein estimator uses $m = \lfloor n^{2/3}\rfloor$, the empirical beta estimator uses $m = n$, and ECBC selects the two marginal degrees separately. In the ECBC fit, two empirical-Bayes Markov chain Monte Carlo chains are run with 1000 burn-in iterations and 2000 retained iterations per chain. The selected degrees are the marginal posterior modes. Acceptance rates, scalar potential scale reduction factors, posterior evaluations, warnings, and numerical failures are retained for every fit. A fit for which either potential scale reduction factor is nonfinite or exceeds $1.20$ is classified as diagnostically nonconverged and is excluded from the accuracy and computation-time summaries.

\begin{figure}[!ht]
\centering
\includegraphics[trim = 1mm 2mm 3mm 4mm, clip, width = 1\textwidth]{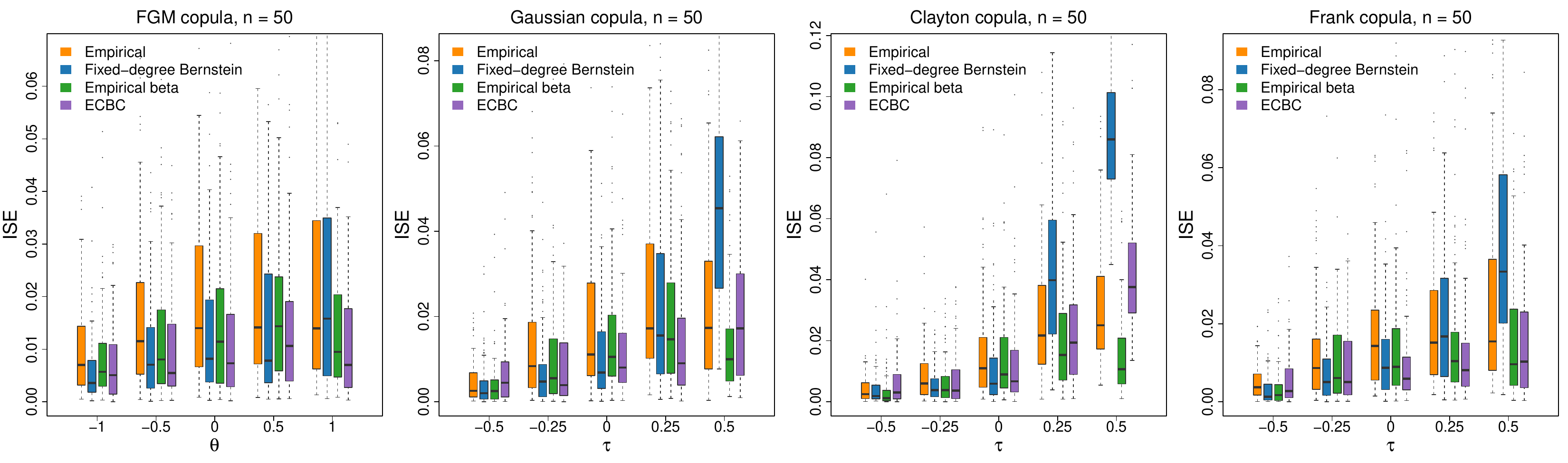} \\
\includegraphics[trim = 1mm 2mm 3mm 4mm, clip, width = 1\textwidth]{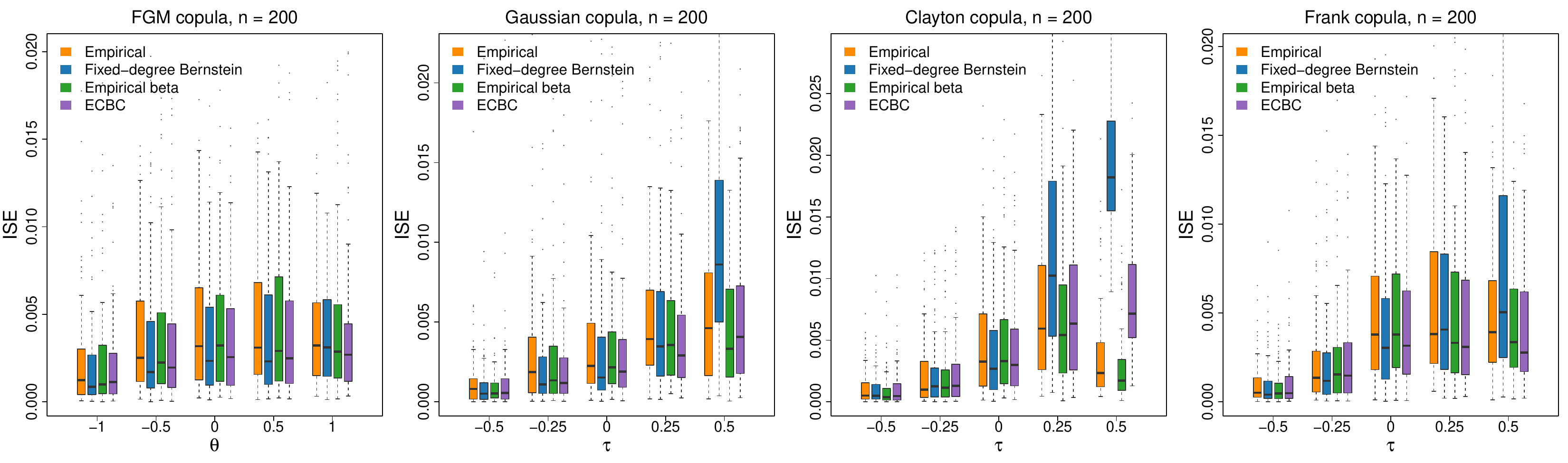} \\[-3mm]
\caption{Common-sample comparison of the empirical copula-based estimator, the fixed-degree Bernstein estimator with $m = \lfloor n^{2/3}\rfloor$, the empirical beta estimator, and the empirical checkerboard Bernstein copula estimator. Boxplots show the Monte Carlo distribution of ISE for the FGM, Gaussian, Clayton, and Frank copulas. The upper and lower rows correspond to $n = 50$ and $n = 200$, respectively. Dependence is indexed by $\theta$ for the FGM copula and by Kendall's $\tau$ for the other copula families. For readability, each panel's vertical range is determined by the 98th percentile of its pooled finite ISE values; the boxplots and numerical summaries are calculated from all finite method-specific evaluations, with diagnostically nonconverged ECBC fits excluded.}
\label{fig:alternative_ise}
\end{figure}

For each method, the comparison reports pointwise bias, variance, and MSE; replication-level ISE; integrated squared bias, integrated variance, and MISE; and computation time. Tables~\ref{tab:alternative_comparison_n50}--\ref{tab:alternative_comparison_n200} in~\ref{app:simulation} report the integrated MSE and the median ECBC degrees. Figure~\ref{fig:alternative_ise} displays the common-sample ISE distributions for all four copula families and both sample sizes. Additional pointwise bias and variance curves, ECBC degree diagnostics, computation times, convergence rates, and numerical failures are reported in Figures~\ref{fig:alternative_pointwise}--\ref{fig:alternative_degrees} and Tables~\ref{tab:alternative_runtime}--\ref{tab:alternative_ecbc_diagnostics} in~\ref{app:simulation}.

Figure~\ref{fig:alternative_ise} and Tables~\ref{tab:alternative_comparison_n50}--\ref{tab:alternative_comparison_n200} show that the fixed-degree Bernstein estimator is competitive, and often has the smallest MISE, in many weak and moderate dependence settings. Its performance deteriorates under strong positive Gaussian, Clayton, and Frank dependence because the deterministic degree introduces appreciable smoothing bias. In those cases the empirical beta estimator, and frequently ECBC, has substantially smaller MISE. ECBC does not dominate uniformly, however, and can be less accurate than simpler estimators in other settings. Thus, the numerical evidence supports the fixed-degree rule as a simple regularized choice under weak or moderate dependence, but not as an automatically optimal choice under strong dependence.

The comparison separates target accuracy from the scope and computational cost of the underlying copula estimator. The fixed-degree Bernstein estimator has one deterministic degree and a closed-form target evaluation. The empirical beta estimator is tuning-free and is a genuine copula, but fixes the degree at $n$. ECBC provides adaptive dimension-specific degrees for a full multivariate copula, at the cost of an empirical-Bayes selection step and associated convergence diagnostics. The computation times in Table~\ref{tab:alternative_runtime} are implementation-specific end-to-end timings for the supplied \textsf{R} code rather than intrinsic complexity bounds. In particular, the ECBC time includes Markov chain Monte Carlo degree selection and posterior diagnostics, whereas the other timings concern direct target evaluations. The reported timings show that ECBC adaptation can yield substantial accuracy gains in some strong-dependence settings but has a markedly higher computational cost than direct evaluation of the fixed-degree estimator. The purpose of the comparison is therefore to display these accuracy, scope, and computational trade-offs rather than to claim uniform dominance of one procedure.

A separate common-bootstrap experiment compares pointwise interval coverage. It uses $K = 50$ Monte Carlo samples, $B = 99$ bootstrap replications, $n \in \{50,200\}$, the FGM cases $\theta \in \{-1,0,1\}$, and the Gaussian, Clayton, and Frank cases with $\tau = 0.5$. Every bootstrap sample is generated from the empirical beta copula and reranked. The empirical copula, fixed-degree Bernstein, and empirical beta estimators are recomputed in every bootstrap sample; ECBC is also recomputed whenever the observed-sample ECBC fit is available and diagnostically converged. The original ECBC degrees are selected as marginal posterior modes from the empirical-Bayes Markov chain Monte Carlo output. To keep the nested experiment computationally feasible, the bootstrap ECBC degrees are approximated by the joint posterior maximizer over an adaptively expanded local integer grid centered at the original selected pair. The initial radius is four; if the maximizer lies on a non-global grid boundary, the radius is doubled, with at most two expansions. This approximation does not rerun the original Markov chain Monte Carlo selector in each bootstrap sample. Table~\ref{tab:alternative_tail_coverage} reports symmetric absolute-root coverage at $p = 0.05$ and $p = 0.10$, and Figure~\ref{fig:alternative_coverage} in~\ref{app:simulation} reports pointwise coverage over the full threshold grid. Since Theorem~\ref{thm:bootstrap.validity} concerns the deterministic-degree estimator, the coverage results for the empirical copula, empirical beta, and ECBC estimators are interpreted only as finite-sample comparisons.

\begin{table}[H]
\centering
\small
\setlength{\tabcolsep}{3pt}
\renewcommand{\arraystretch}{0.9}
\caption{Pointwise coverage of the symmetric absolute-root intervals at the two smallest thresholds in the common-bootstrap comparison. Plug-in binomial Monte Carlo standard errors, calculated from the successful interval replications, are shown in parentheses.}
\label{tab:alternative_tail_coverage}
\begin{tabular}{lcrrrrrr}
\toprule[1.5pt]
Copula & Dependence & $n$ & $p$ & Empirical & Fixed degree & Beta & ECBC \\
\midrule[1pt]
FGM & $-1$ & $50$ & $0.05$ & 0.080 (0.038) & 0.880 (0.046) & 0.540 (0.070) & 0.860 (0.049) \\
FGM & $-1$ & $50$ & $0.10$ & 0.380 (0.069) & 0.940 (0.034) & 0.720 (0.063) & 0.900 (0.042) \\
FGM & $-1$ & $200$ & $0.05$ & 0.120 (0.046) & 0.580 (0.070) & 0.320 (0.066) & 0.500 (0.071) \\
FGM & $-1$ & $200$ & $0.10$ & 0.480 (0.071) & 0.800 (0.057) & 0.600 (0.069) & 0.680 (0.066) \\
FGM & $0$ & $50$ & $0.05$ & 0.340 (0.067) & 0.580 (0.070) & 0.600 (0.069) & 0.720 (0.063) \\
FGM & $0$ & $50$ & $0.10$ & 0.640 (0.068) & 0.700 (0.065) & 0.740 (0.062) & 0.740 (0.062) \\
FGM & $0$ & $200$ & $0.05$ & 0.460 (0.070) & 0.680 (0.066) & 0.540 (0.070) & 0.660 (0.067) \\
FGM & $0$ & $200$ & $0.10$ & 0.740 (0.062) & 0.760 (0.060) & 0.760 (0.060) & 0.780 (0.059) \\
FGM & $1$ & $50$ & $0.05$ & 0.480 (0.071) & 0.920 (0.038) & 0.860 (0.049) & 0.920 (0.038) \\
FGM & $1$ & $50$ & $0.10$ & 0.880 (0.046) & 0.840 (0.052) & 0.920 (0.038) & 0.940 (0.034) \\
FGM & $1$ & $200$ & $0.05$ & 0.740 (0.062) & 0.820 (0.054) & 0.820 (0.054) & 0.820 (0.054) \\
FGM & $1$ & $200$ & $0.10$ & 0.840 (0.052) & 0.860 (0.049) & 0.920 (0.038) & 0.900 (0.042) \\
\midrule
Gaussian & $0.5$ & $50$ & $0.05$ & 0.320 (0.066) & 0.020 (0.020) & 0.740 (0.062) & 0.260 (0.062) \\
Gaussian & $0.5$ & $50$ & $0.10$ & 0.660 (0.067) & 0.080 (0.038) & 0.760 (0.060) & 0.540 (0.070) \\
Gaussian & $0.5$ & $200$ & $0.05$ & 0.800 (0.057) & 0.280 (0.063) & 0.880 (0.046) & 0.694 (0.066) \\
Gaussian & $0.5$ & $200$ & $0.10$ & 0.920 (0.038) & 0.560 (0.070) & 0.920 (0.038) & 0.816 (0.055) \\
\midrule
Clayton & $0.5$ & $50$ & $0.05$ & 0.760 (0.060) & 0.000 (0.000) & 0.840 (0.052) & 0.000 (0.000) \\
Clayton & $0.5$ & $50$ & $0.10$ & 0.860 (0.049) & 0.000 (0.000) & 0.940 (0.034) & 0.000 (0.000) \\
Clayton & $0.5$ & $200$ & $0.05$ & 0.940 (0.034) & 0.000 (0.000) & 0.960 (0.028) & 0.000 (0.000) \\
Clayton & $0.5$ & $200$ & $0.10$ & 0.920 (0.038) & 0.000 (0.000) & 0.960 (0.028) & 0.360 (0.068) \\
\midrule
Frank & $0.5$ & $50$ & $0.05$ & 0.500 (0.071) & 0.520 (0.071) & 0.920 (0.038) & 0.900 (0.042) \\
Frank & $0.5$ & $50$ & $0.10$ & 0.820 (0.054) & 0.420 (0.070) & 0.940 (0.034) & 0.900 (0.042) \\
Frank & $0.5$ & $200$ & $0.05$ & 0.820 (0.054) & 0.780 (0.059) & 0.860 (0.049) & 0.860 (0.049) \\
Frank & $0.5$ & $200$ & $0.10$ & 0.920 (0.038) & 0.820 (0.054) & 0.940 (0.034) & 0.920 (0.038) \\
\bottomrule[1.5pt]
\end{tabular}
\end{table}

The common-bootstrap experiment is intentionally exploratory. With $K = 50$, the binomial Monte Carlo standard error at the nominal coverage probability $0.95$ is approximately $0.031$, and the corresponding approximate 95\% margin of Monte Carlo error is $0.060$. Moreover, $B = 99$ yields a relatively coarse estimate of an upper bootstrap quantile. The entries in Table~\ref{tab:alternative_tail_coverage} and the curves in Figure~\ref{fig:alternative_coverage} should therefore be read as qualitative diagnostics rather than precise rankings of the four methods.

The diagnostic nevertheless reveals substantial and systematic undercoverage in several settings. The empirical-copula intervals have severe lower-boundary undercoverage in several FGM cases, most notably for $\theta = -1$, where the fixed-degree Bernstein and ECBC intervals perform substantially better. Under strong positive dependence, the fixed-degree Bernstein intervals have very low coverage near the lower boundary for the Gaussian and Clayton copulas with $\tau = 0.5$, and their coverage is also below nominal for the Frank copula. The ECBC intervals show pronounced undercoverage for the strong Gaussian and Clayton cases, although they perform considerably better for the strong Frank case. The empirical beta intervals are generally closer to the nominal level in these cases, but no estimator is uniformly well calibrated over all families and thresholds. Increasing the sample size from $n = 50$ to $n = 200$ does not uniformly eliminate the problem. These patterns are consistent with smoothing bias that is not negligible relative to the sampling error at the reported sample sizes. Assumption~\ref{ass:integrated.bias}, together with $\sqrt n/m \to 0$, makes the deterministic integrated Bernstein bias negligible on the $n^{-1/2}$ scale only asymptotically. Consequently, the finite-sample intervals based on the fixed degree should be used cautiously under strong positive dependence, especially for small $p$; reducing smoothing or applying an explicit bias correction is preferable when the coverage diagnostic indicates substantial bias.

\subsection{Sensitivity analysis of the Bernstein degree \texorpdfstring{$m$}{m}}\label{sec:degree}

The Bernstein degree $m$ controls the amount of smoothing applied to the empirical copula. Small values of $m$ produce a heavily smoothed estimator and may therefore introduce substantial smoothing bias, whereas large values of $m$ make the Bernstein estimator closer to the empirical copula and may reduce the stabilizing effect of smoothing. The asymptotic bias--variance calculations motivate the degree choice
\[
m = \lfloor c \, n^{2/3}\rfloor, \qquad c > 0.
\]
To assess whether this rule is also reasonable in finite samples, the MISE was computed over a grid of Bernstein degrees. The choice $c = 1$ is the default used in the preceding simulations, while $c \in \{0.75, 1.25\}$ gives two reference perturbations around that default.

\begin{figure}[!ht]
\centering
\includegraphics[trim = 3mm 2mm 7mm 6mm, clip, width = 0.48\textwidth]{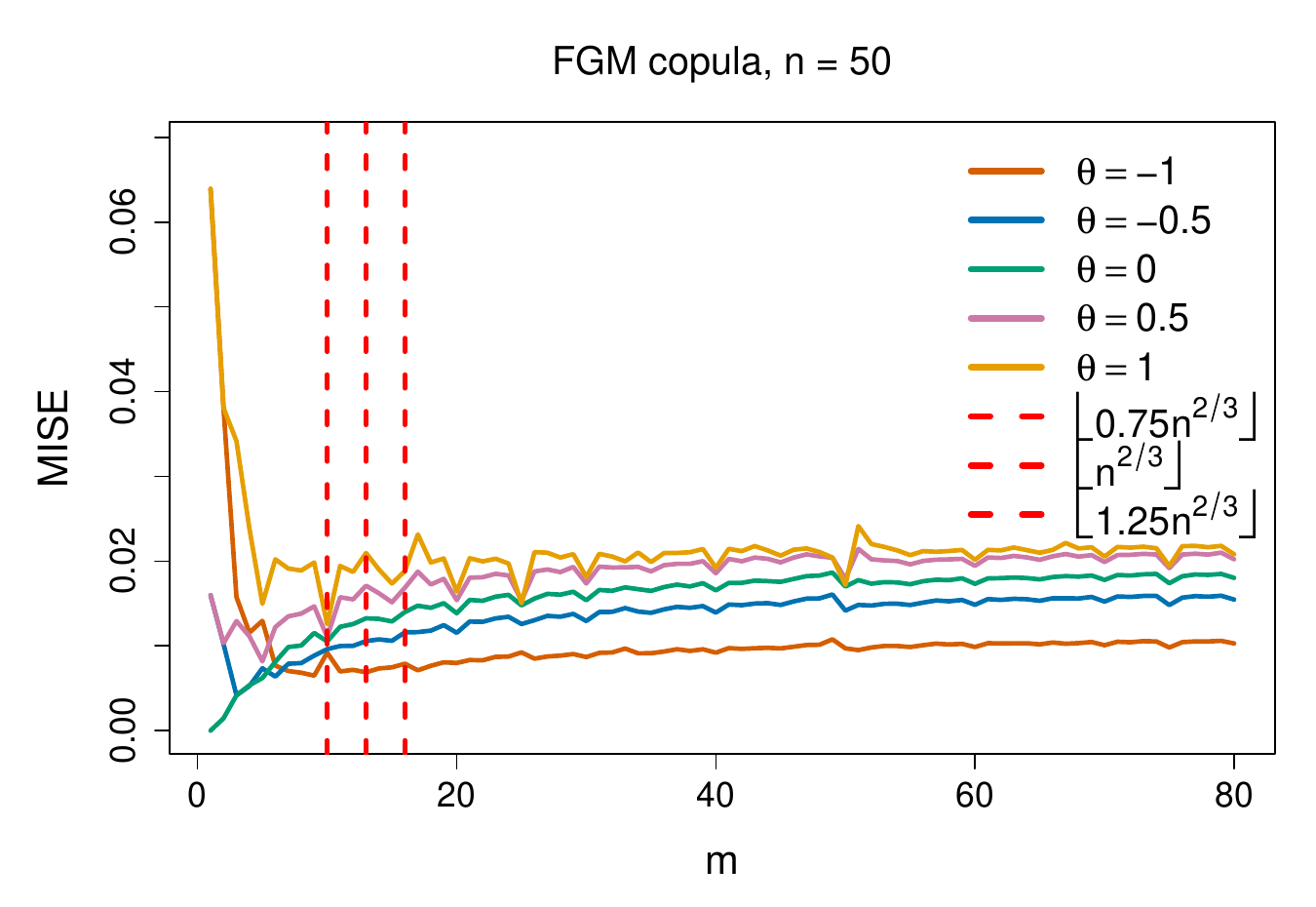}
\quad
\includegraphics[trim = 3mm 2mm 7mm 6mm, clip, width = 0.48\textwidth]{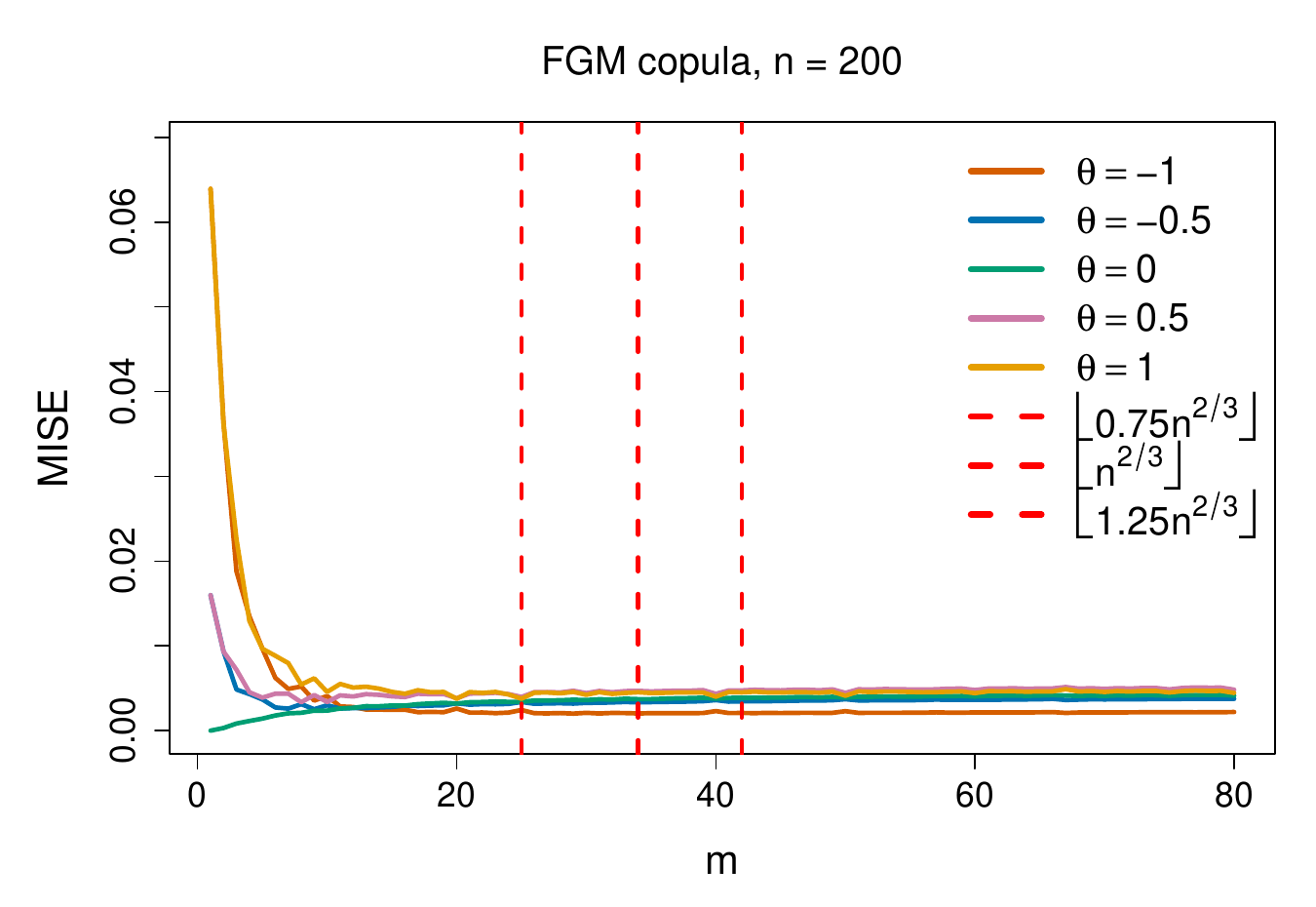} \\
\includegraphics[trim = 3mm 2mm 7mm 6mm, clip, width = 0.48\textwidth]{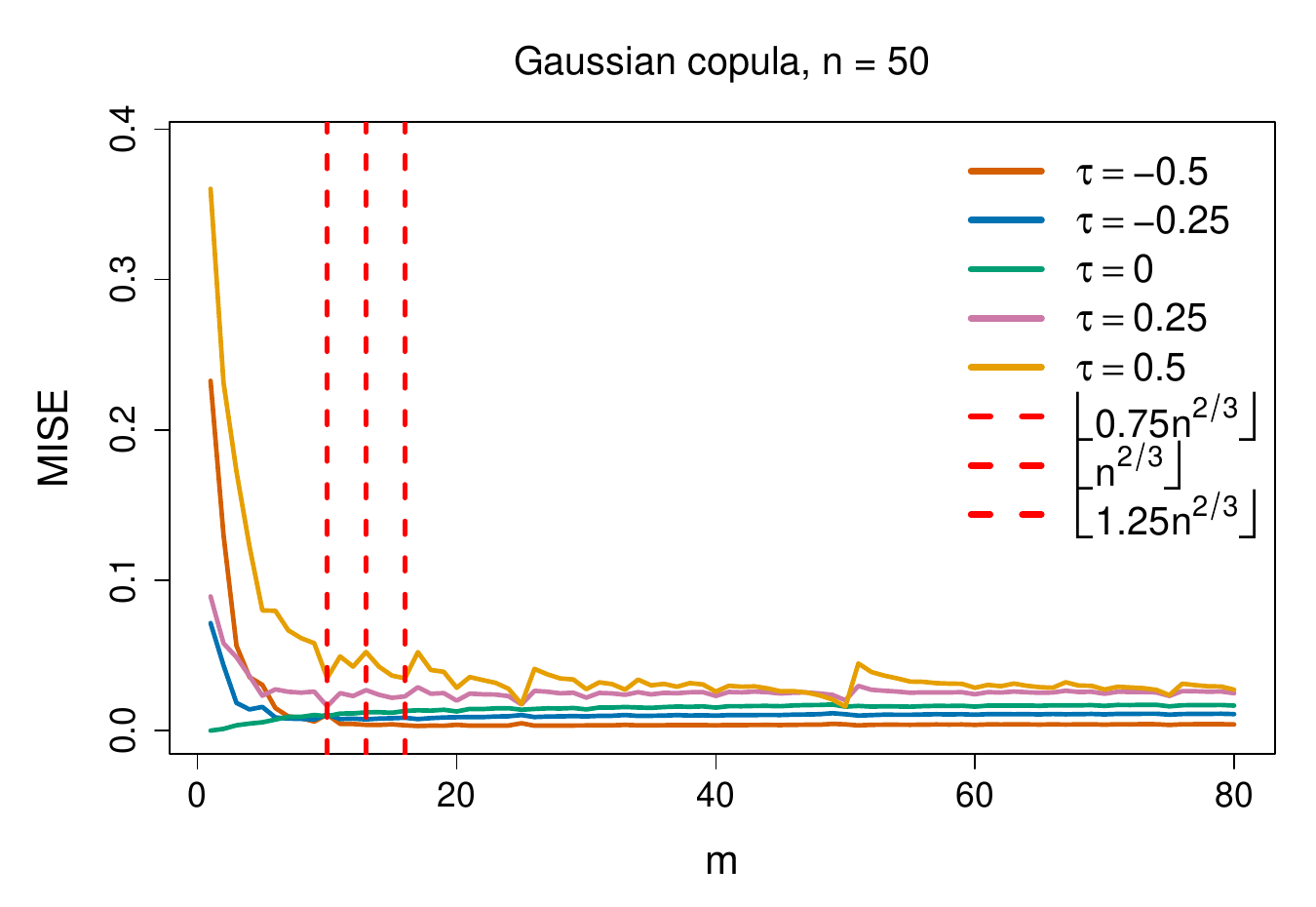}
\quad
\includegraphics[trim = 3mm 2mm 7mm 6mm, clip, width = 0.48\textwidth]{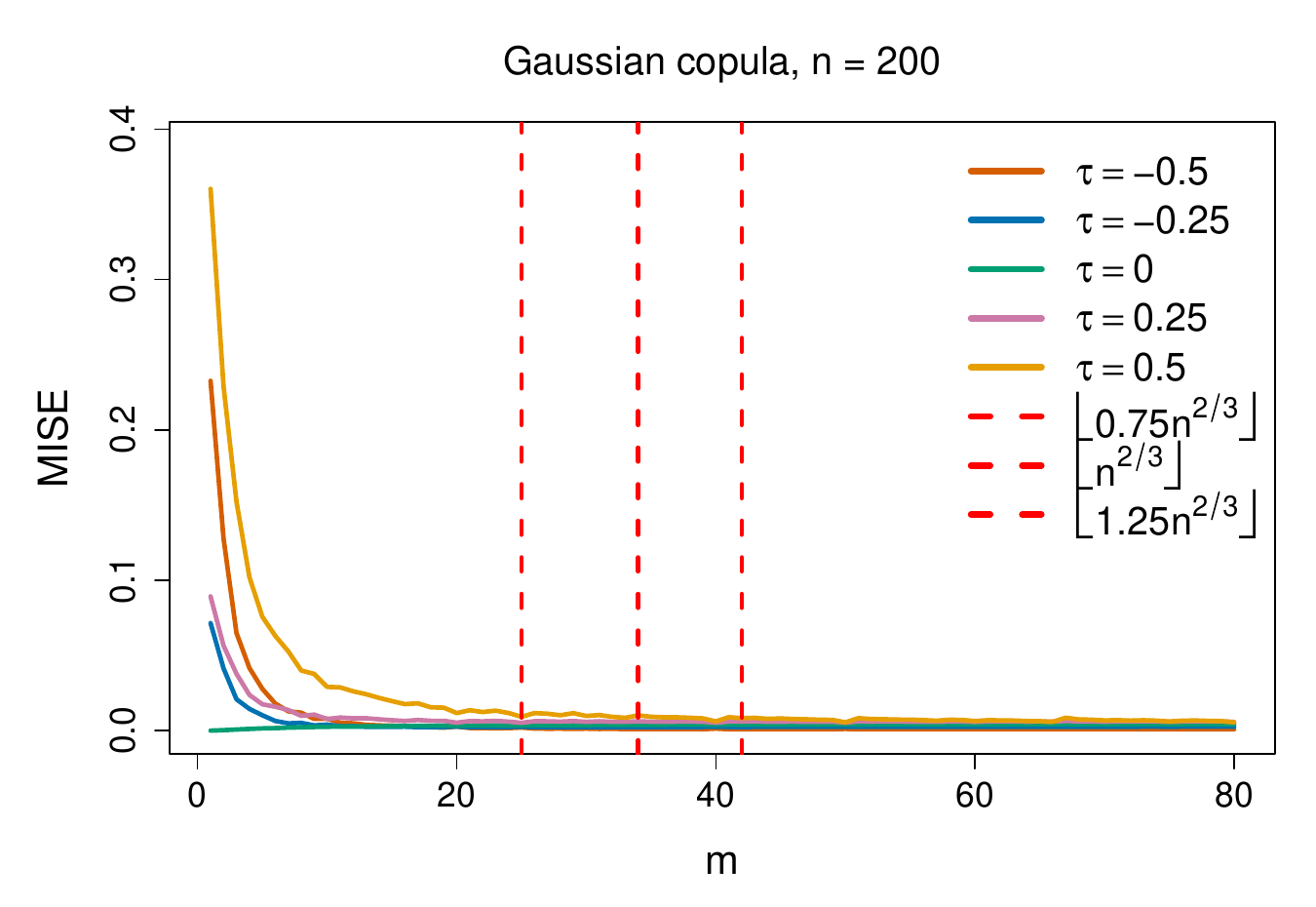} \\[-3mm]
\caption{Effect of the Bernstein polynomial degree $m$ on finite-sample accuracy. Rows correspond to the Farlie--Gumbel--Morgenstern (FGM) and Gaussian copulas, and columns correspond to $n = 50$ and $n = 200$. The vertical red dashed lines indicate $m = \lfloor c n^{2/3}\rfloor$ for $c \in \{0.75, 1, 1.25\}$. The middle line, $c = 1$, is the default used in the preceding simulations, while the other two lines are reference perturbations used only in this sensitivity analysis.}
\label{fig:sim_MISE}
\end{figure}

\begin{figure}[H]
\centering
\includegraphics[trim = 3mm 2mm 7mm 6mm, clip, width = 0.48\textwidth]{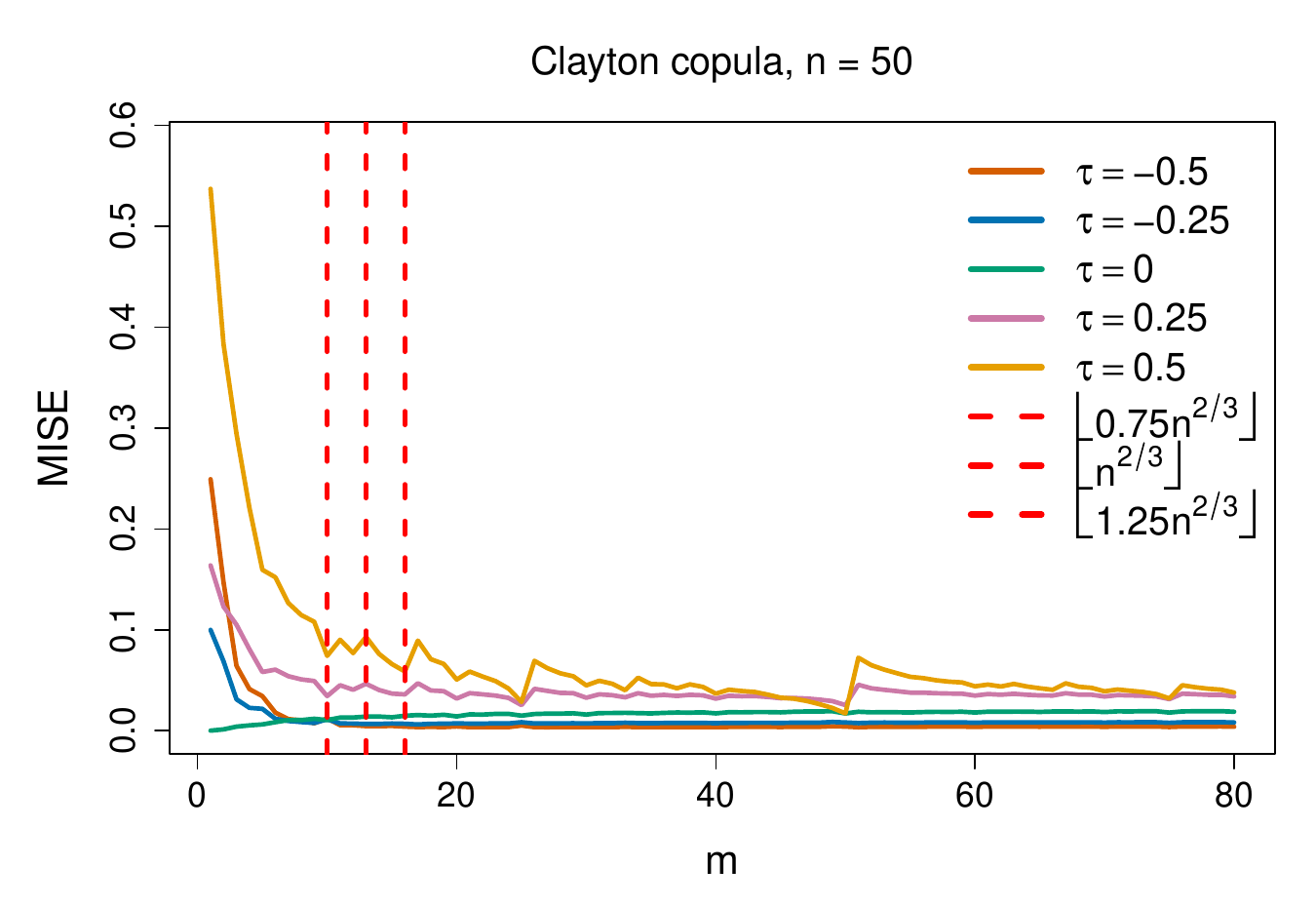}
\quad
\includegraphics[trim = 3mm 2mm 7mm 6mm, clip, width = 0.48\textwidth]{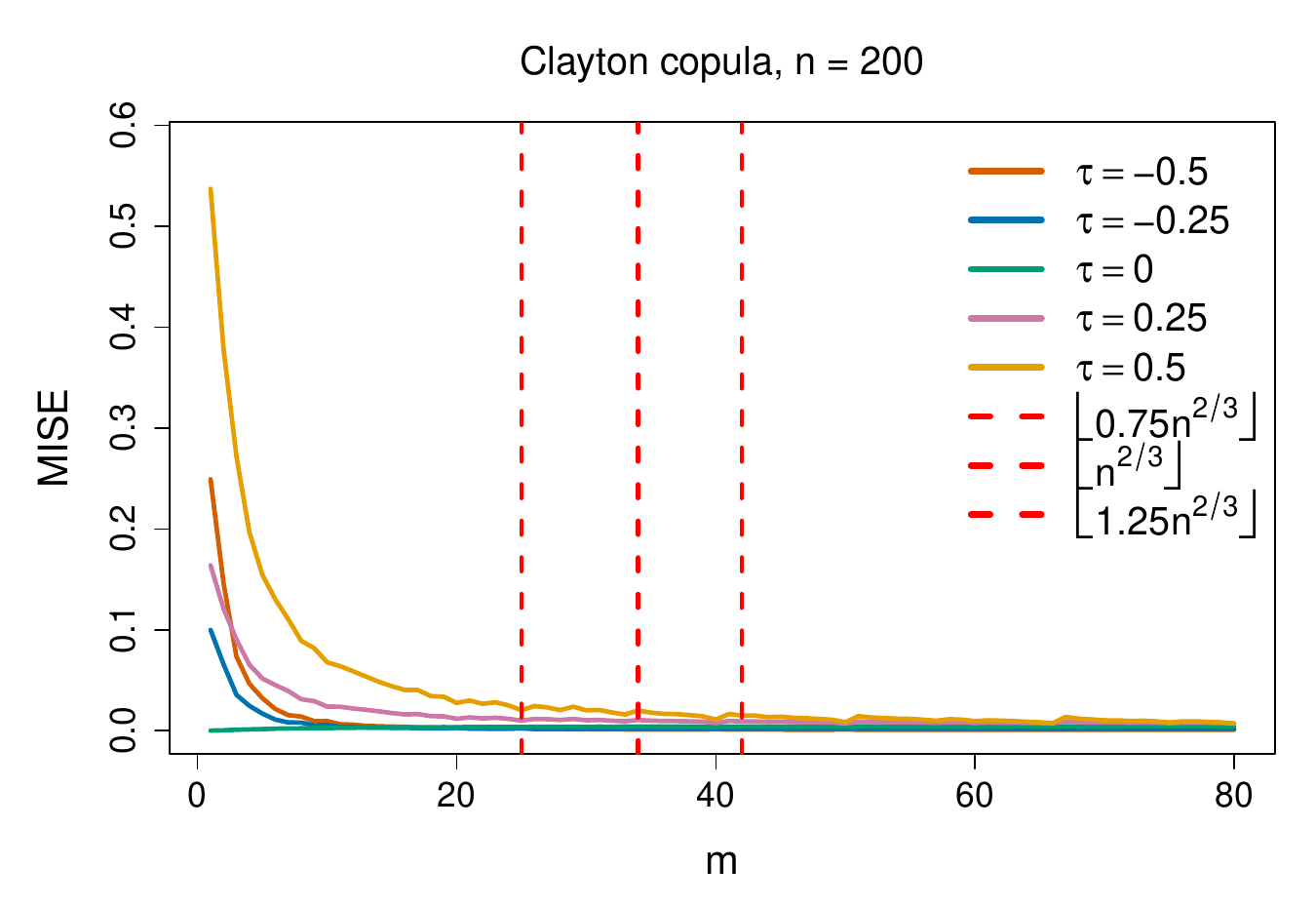} \\
\includegraphics[trim = 3mm 2mm 7mm 6mm, clip, width = 0.48\textwidth]{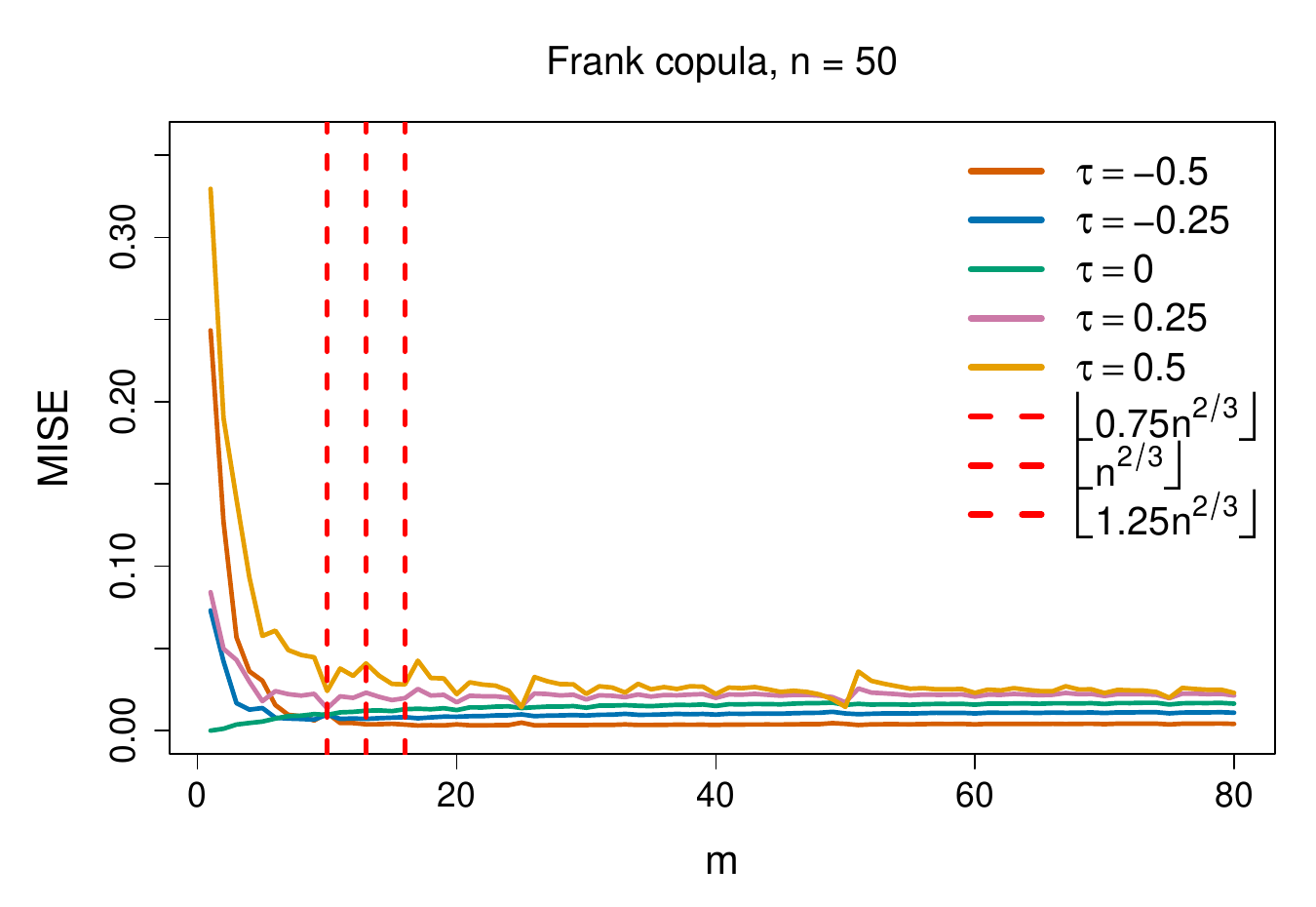}
\quad
\includegraphics[trim = 3mm 2mm 7mm 6mm, clip, width = 0.48\textwidth]{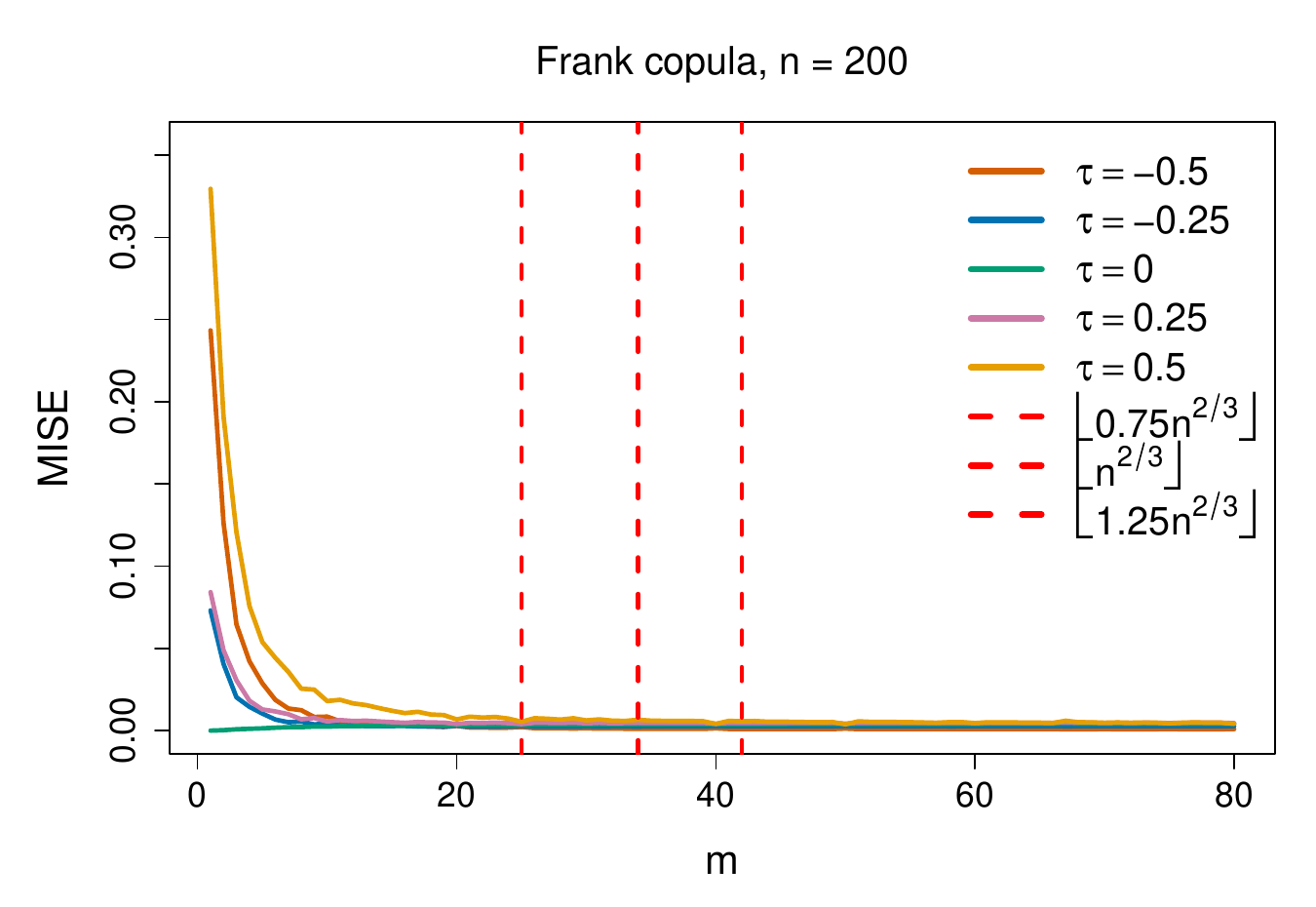} \\[-3mm]
\caption{Effect of the Bernstein polynomial degree $m$ on finite-sample accuracy. Rows correspond to the Clayton and Frank copulas, and columns correspond to $n = 50$ and $n = 200$. The vertical red dashed lines indicate $m = \lfloor c n^{2/3}\rfloor$ for $c \in \{0.75, 1, 1.25\}$. The middle line, $c = 1$, is the default used in the preceding simulations, while the other two lines are reference perturbations used only in this sensitivity analysis.}
\label{fig:sim_MISE_2}
\end{figure}

We used the same Monte Carlo design as in Section~\ref{sec:MISE_comparison}. Figures~\ref{fig:sim_MISE}~and~\ref{fig:sim_MISE_2} summarize the sensitivity of the MISE to the Bernstein degree $m$. The effect of very small degrees depends strongly on the underlying dependence. Under strong positive or negative dependence, especially for the Gaussian, Clayton, and Frank copulas, very small $m$ often produces large MISE, which decreases rapidly as $m$ increases and the smoothing bias is reduced. This pattern is not universal. Under independence, $B_1(C_n) = \Pi$ and hence $\widehat{\rho}_{1,n}(p) = 0$, so the target curve is estimated without error. In several weak-dependence settings, the MISE is also already small at the lowest degrees or initially increases.

The vertical reference lines in each panel correspond to the three values of $c$. The middle reference line, corresponding to $c = 1$, gives the default rule $m = \lfloor n^{2/3}\rfloor$, with $m = 13$ for $n = 50$ and $m = 34$ for $n = 200$. For $n = 50$, the default lies beyond the region of most severe oversmoothing in many weak and moderate dependence settings. It is not uniformly near a stable MISE plateau, however. In particular, the strongly positively dependent Gaussian, Clayton, and Frank curves continue to decrease materially for degrees well above $m = 13$. The nearby choices $c \in \{0.75, 1, 1.25\}$ therefore show only local sensitivity around the default and do not establish that the default is close to the finite-sample optimum in these cases.

The visible changes near $m = n = 50$ in the left panels are finite-sample grid-alignment effects. At $m = n$, all Bernstein grid points coincide with the empirical rank grid and $C_{n,n}$ is the empirical beta copula. For $m > n$, repeated empirical-copula levels occur on the finer Bernstein grid and $C_{m,n}$ need not be a copula, so the MISE need not vary smoothly across $m = n$.

For $n = 200$, the MISE profiles are flatter near the three reference choices in many weak and moderate dependence settings, and the default value $m = 34$ avoids the high-error region associated with excessive smoothing. The strong positive cases can nevertheless continue to benefit from larger degrees. Thus, the greater local flatness at $n = 200$ does not imply that $m = 34$ is uniformly near the finite-sample minimizer.

The sensitivity analysis therefore supports $m = \lfloor n^{2/3}\rfloor$ only as a transparent, sample-size-dependent starting value that avoids the most severe oversmoothing and satisfies the degree regime in Theorem~\ref{thm:functional_clt_rho}. It does not support a claim that this rule is uniformly stable or close to the exact finite-sample minimizer. Moderate perturbations around the default have limited impact in many weak and moderate dependence settings, whereas materially larger degrees can be preferable under strong positive dependence. This limitation is consistent with the MISE and coverage comparisons above and motivates data-adaptive degree selection or bias diagnostics in settings with pronounced copula curvature.

The deterministic rule and the empirical Bayes selector of \citet{LuGhosh2023} serve different purposes. The former is transparent, computationally inexpensive, and places the estimator in a prespecified asymptotic regime, whereas the latter adapts dimension-specific degrees to the data when estimating the full copula. Section~\ref{sec:alternative.comparison} makes this distinction concrete by fitting ECBC with data-dependent degrees. Its exploratory coverage experiment uses an adaptively expanded local-grid joint-posterior maximization in each bootstrap sample, centered at the original Markov chain Monte Carlo selection, rather than rerunning the original selector. An adaptive degree $\widehat m$ can also be incorporated directly into the proposed estimator by evaluating $\widehat{\rho}_{\widehat m,n}(p)$, but a corresponding curve-level validity theorem would require stochastic control uniformly over the candidate degrees. The deterministic degree is therefore retained for the proposed inferential procedure, while the adaptive ECBC analysis is reported as a separate finite-sample benchmark.

\section{Real-data application}\label{sec:application}

We illustrate the proposed Bernstein-smoothed lower-tail Spearman's rho estimator using the Loss--ALAE insurance claims data. This dataset is a standard benchmark in actuarial copula modeling and was used by~\citet{FreesValdez1998} to study dependence between insurance losses and claim-specific expenses. It has subsequently been analyzed in several copula and dependence-modeling studies, including bivariate copula fitting~\citep{KlugmanParsa1999}, Archimedean copula modeling with censoring~\citep{DenuitPurcaruVanKeilegom2006}, copula goodness-of-fit assessment~\citep{GenestQuessyRemillard2006}, and semiparametric copula model selection under censoring~\citep{ChenFanPouzoYing2010}. The data are publicly available in \textsf{R} as the \texttt{lossalae} dataset in the \textbf{evd} package~\citep{EvdPackage}. The dataset contains $n = 1500$ general liability claims, with two variables: the indemnity payment, denoted by $X$, and the allocated loss adjustment expense (ALAE), denoted by $Y$. Both variables are measured in US dollars. The ALAE variable records claim-specific settlement expenses, such as investigation expenses and legal fees. The dataset's \texttt{capped} attribute identifies indemnity payments capped at their policy limits. The main analysis uses the observed capped values descriptively and does not model the underlying censoring mechanism.

This dataset provides a useful benchmark for illustrating the proposed method because it contains paired indemnity payments and allocated loss adjustment expenses, two quantities whose dependence is of actuarial interest. While the co-occurrence of large indemnity payments and large ALAE is important for reserving, pricing, and reinsurance, the present analysis focuses on the normalized copula integral over the lower-left square $[0,p]^2$, which emphasizes joint behavior at small marginal ranks. The goal of this application is therefore not to model very large claims directly, but to demonstrate how the lower-tail Spearman's rho curve can describe local dependence as the threshold $p$ varies. In particular, comparing the estimated lower-tail Spearman's rho curve with the global Spearman correlation illustrates how a single global rank-correlation measure may obscure dependence patterns that are specific to a particular region of the copula.

We computed both the empirical copula-based estimator $\widehat{\rho}_n(p)$ and the Bernstein-smoothed estimator $\widehat{\rho}_{m,n}(p)$. For the Bernstein estimator, we used the rule $m = \lfloor n^{2/3}\rfloor$, which gives $m = 131$ for the Loss--ALAE sample. Pointwise 95\% confidence intervals were obtained from $B = 999$ empirical-beta bootstrap samples with $m_\# = n = 1500$, and every bootstrap sample was reranked before re-estimation. Since both variables contain ties, the two marginal rank permutations were obtained by reproducible random tie breaking. The global rank-correlation benchmark was computed as the Pearson correlation of those same two rank permutations, so the curve and the benchmark use one common tie convention. The symmetric absolute-root interval is used in the main application figure, while the equal-tailed basic interval is reported as a sensitivity check. The theory in Sections~\ref{sec:asymptotic}~and~\ref{sec:CI_bootstrap} assumes continuous margins and therefore does not cover these tied data. Accordingly, the intervals in this application are conditional on the recorded tie breaking and should be interpreted as an illustration rather than as a formally justified analysis for mixed or discrete margins.

\begin{figure}[!ht]
\centering
\includegraphics[trim = 3mm 2mm 7mm 6mm, clip, width = 0.46\textwidth]{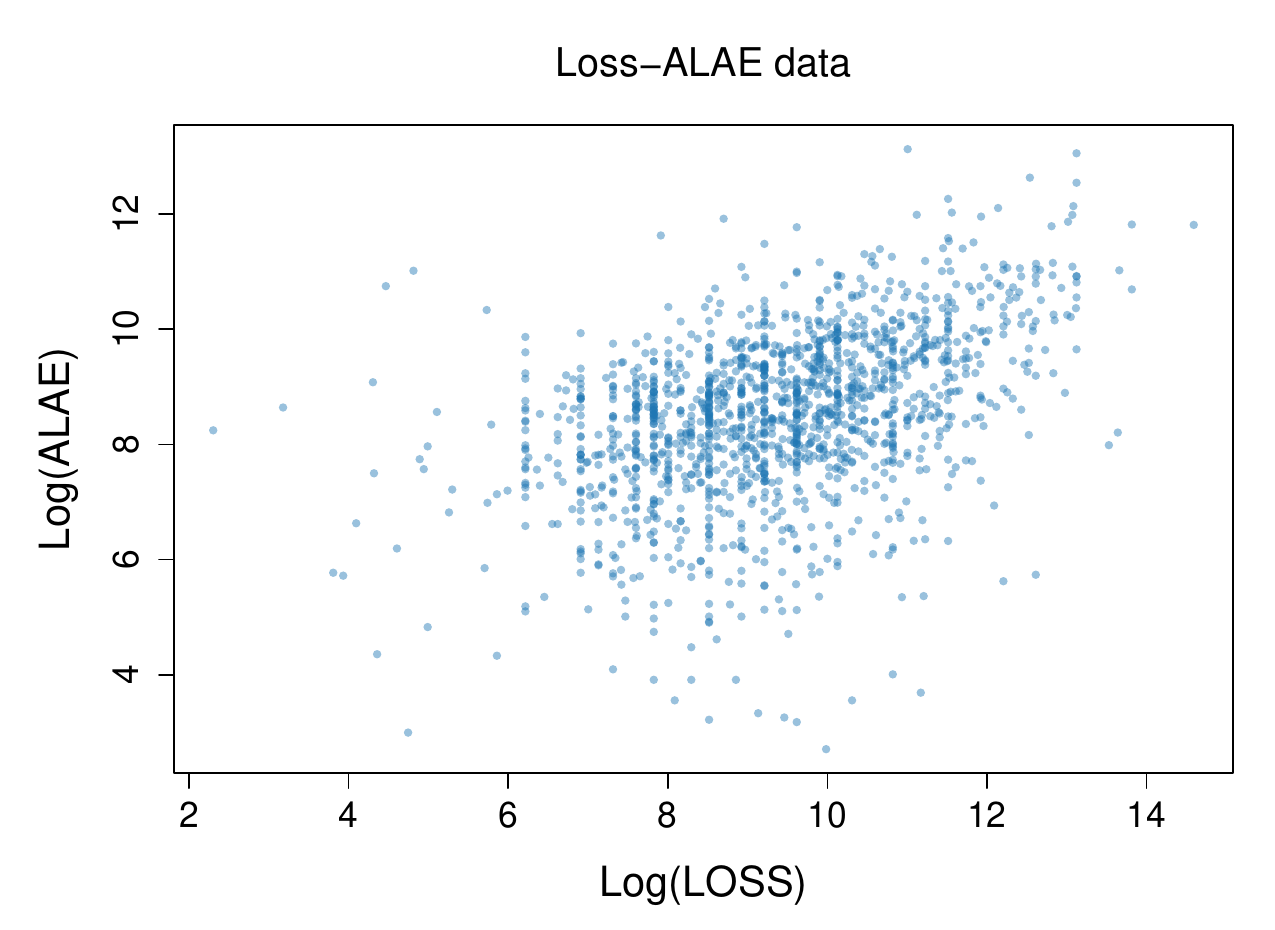}
\quad
\includegraphics[trim = 3mm 2mm 7mm 6mm, clip, width = 0.46\textwidth]{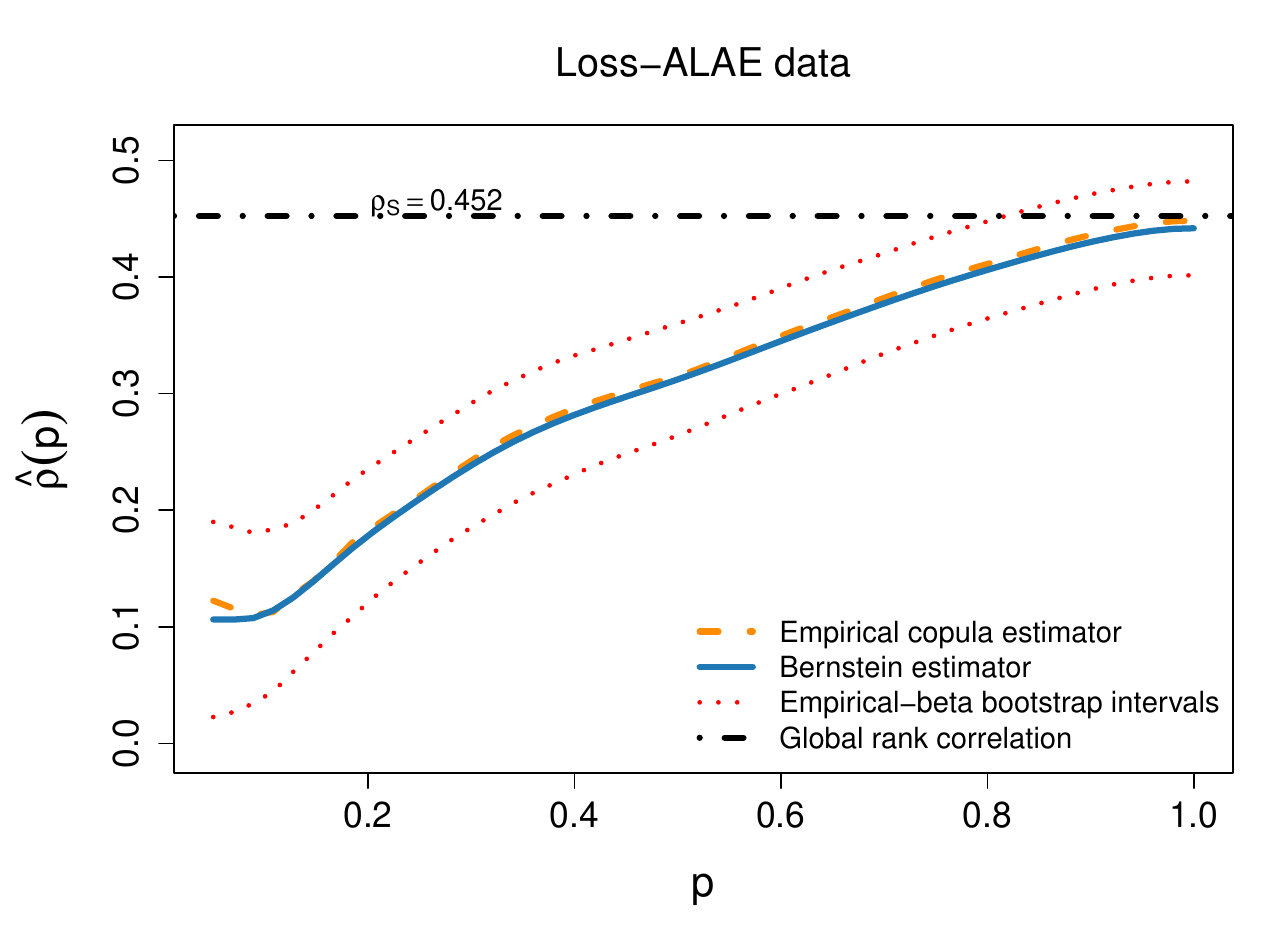}
\caption{Application to the Loss--ALAE insurance claims data. The left panel shows the scatter plot of $\log(\mathrm{Loss})$ and $\log(\mathrm{ALAE})$. The right panel displays the estimated lower-tail Spearman's rho curve. The orange dashed curve is the empirical copula-based estimator $\widehat{\rho}_{n}(p)$, and the blue solid curve is the Bernstein-smoothed estimator $\widehat{\rho}_{m,n}(p)$, with $m = \lfloor n^{2/3}\rfloor = 131$. The red dotted curves give the endpoints of the nominal 95\% symmetric absolute-root empirical-beta bootstrap intervals. The black dash-dotted horizontal line gives the global rank correlation computed from the same reproducibly tie-broken rank permutations as the two curve estimators.}
\label{fig:lossalae_application}
\end{figure}

Figure~\ref{fig:lossalae_application} displays the dependence structure in the Loss--ALAE data and the estimated lower-tail Spearman's rho curves. The scatter plot of $\log(\mathrm{Loss})$ and $\log(\mathrm{ALAE})$ shows a clear positive association, and the corresponding global rank correlation is displayed as a benchmark. The empirical copula-based and Bernstein-smoothed lower-tail curves are nearly identical across the range of $p$, which is expected for the relatively large sample size $n = 1500$ and suggests that the estimated pattern is not an artifact of smoothing. Compared with the global rank correlation, the lower-tail Spearman's rho curve is noticeably smaller for small and moderate values of $p$ in this application. This suggests that concordance among observations with relatively small ranks in both Loss and ALAE is weaker than the overall rank association. As $p$ increases, the curve approaches the global measure, since the population lower-tail functional increasingly incorporates a larger portion of the copula and satisfies $\rho(1) = \rho_S$. The empirical-copula plug-in value at $p = 1$ need not equal the usual finite-sample rank correlation exactly because the plug-in uses the pseudo-observations $R_i/n$, but the discrepancy is of finite-sample order and is small here. Because the normalization $D(p)$ changes with $p$, however, this comparison does not isolate which parts of the complement of $[0,p]^2$ account for the difference from the global correlation.

\begin{figure}[!ht]
\centering
\includegraphics[trim = 3mm 4mm 7mm 6mm, clip, width = 0.45\textwidth]{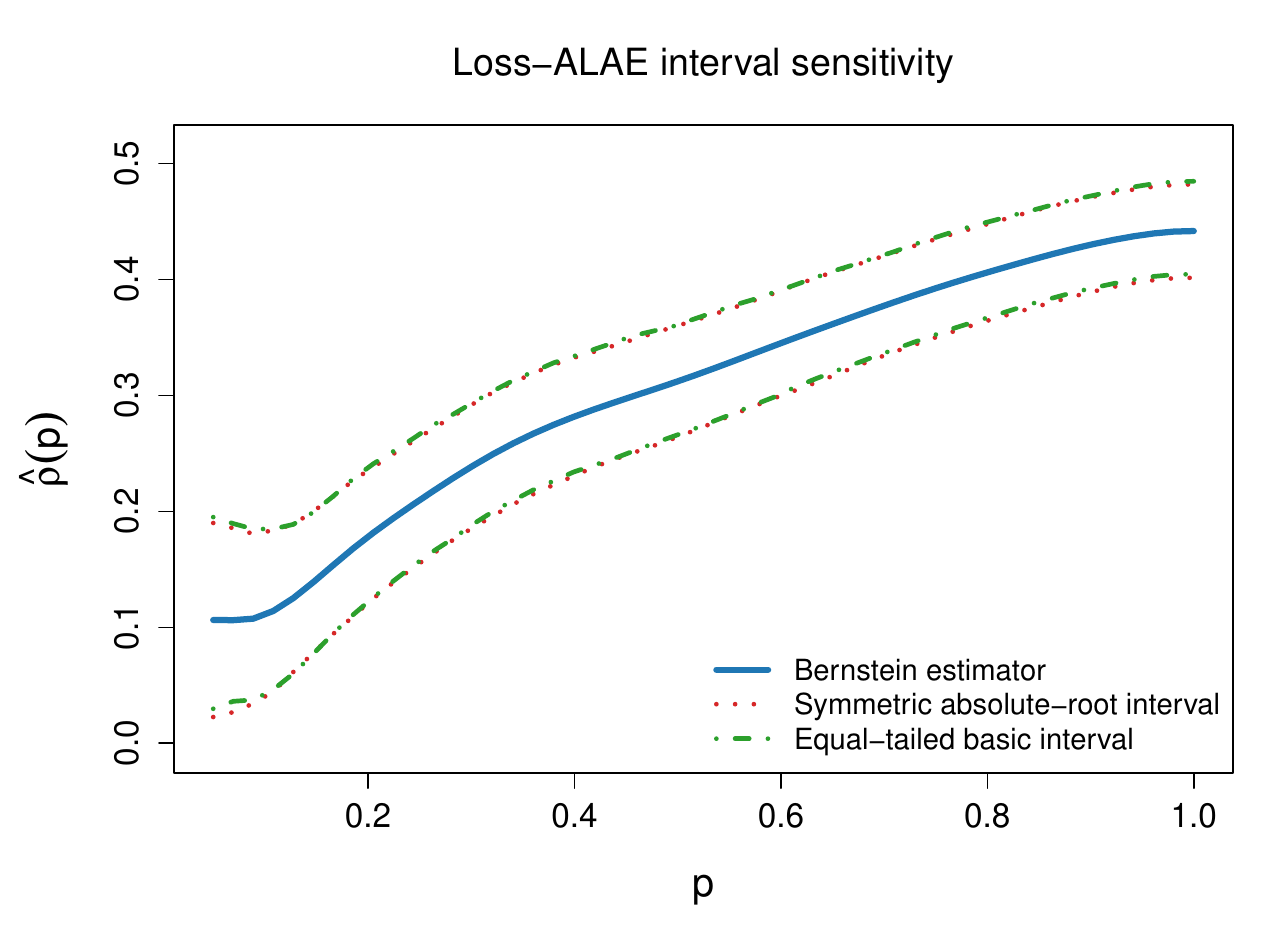}
\caption{Sensitivity of the Loss--ALAE pointwise bootstrap intervals to the bootstrap construction. The blue solid curve is the Bernstein-smoothed estimator, the red dotted curves are the symmetric absolute-root interval endpoints, and the green dash-dotted curves are the equal-tailed basic interval endpoints. Both nominal 95\% intervals use the same empirical-beta bootstrap samples with $m_\# = n = 1500$.}
\label{fig:lossalae_interval_sensitivity}
\end{figure}

\begin{figure}[H]
\centering
\includegraphics[trim = 3mm 4mm 7mm 6mm, clip, width = 0.95\textwidth]{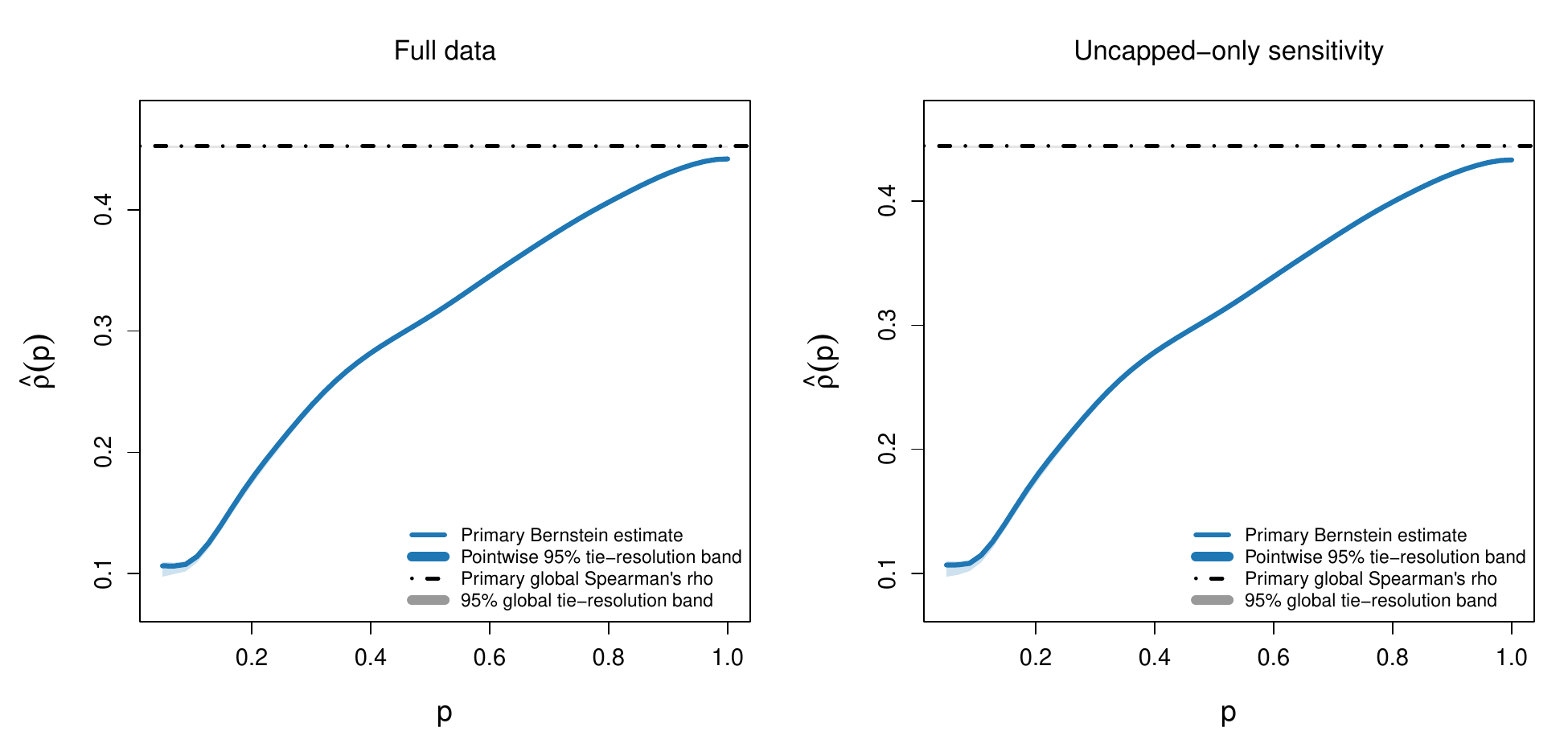}
\caption{Sensitivity of the Loss--ALAE point estimates to random tie resolution. The left panel uses the full data, and the right panel omits observations identified by the data source as policy-limit capped. The blue solid curve is the Bernstein estimate under the primary seed, and the blue shaded region gives the pointwise 2.5\% and 97.5\% quantiles over 100 reproducible tie-breaking seeds. The black dash-dotted line is the global rank correlation under the primary seed, and the gray region gives its 2.5\% and 97.5\% quantiles over the same seeds. Omitting capped observations is a descriptive sensitivity analysis rather than a censoring correction.}
\label{fig:lossalae_tie_sensitivity}
\end{figure}

Figure~\ref{fig:lossalae_interval_sensitivity} compares the symmetric absolute-root and equal-tailed basic intervals obtained from the same empirical-beta bootstrap samples. This provides a direct sensitivity check for the symmetry imposed by the main interval construction. Figure~\ref{fig:lossalae_tie_sensitivity} and Table~\ref{tab:lossalae_tie_sensitivity} report a separate point-estimate sensitivity analysis that repeats the random tie breaking over 100 deterministic seeds and records the resulting variation in the lower-tail curves and global rank correlation. It is computed for the full data and after omitting the observations identified by the data source as policy-limit capped. The resulting variation is small: the global rank-correlation intervals have width $0.002$, and the maximum pointwise curve-band widths are $0.012$ and $0.013$ for the full-data and uncapped-only analyses, respectively. The latter calculation is only a descriptive sensitivity analysis and is not a correction for censoring. Similarly, the variation across seeds assesses dependence on the chosen tie resolution but does not replace inference developed specifically for distributions with mass points.

\begin{table}[H]
\centering
\small
\setlength{\tabcolsep}{4.1pt}
\renewcommand{\arraystretch}{0.9}
\caption{Sensitivity of the Loss--ALAE analysis to random tie resolution over $100$ reproducible seeds. For every seed, the same randomized rank vectors define both the Bernstein curve and the global Spearman benchmark. The global interval contains the pointwise $2.5\%$ and $97.5\%$ quantiles over tie resolutions, and the last column is the largest width of the corresponding pointwise band for the curve. The uncapped-only calculation is a descriptive sensitivity analysis obtained by omitting observations identified by the data source as policy-limit capped; it is not a censoring correction.}
\label{tab:lossalae_tie_sensitivity}
\begin{tabular}{lrrrrr}
\toprule[1.5pt]
Analysis & $n$ & Capped omitted & Primary $\widehat\rho_S$ & Global 95\% interval & Max. curve-band width \\
\midrule[1pt]
Full data & 1500 & 0 & 0.452 & [0.451, 0.453] & 0.012 \\
Uncapped-only sensitivity & 1466 & 34 & 0.444 & [0.443, 0.445] & 0.013 \\
\bottomrule[1.5pt]
\end{tabular}
\end{table}

\section{Summary and outlook}\label{sec:conclusion}

This paper studied the target-specific Bernstein plug-in estimator $\widehat{\rho}_{m,n}(p)$ for the lower-tail Spearman's rho curve. The empirical Bernstein copula itself is an established smoothing construction; the contribution here is its use for estimation and inference on the entire lower-tail Spearman curve, together with an exact incomplete-beta representation and degree conditions tailored to this functional. The estimator is uniformly strongly consistent on compact intervals away from zero whenever $m \to \infty$. Under the integrated Bernstein-bias condition in Assumption~\ref{ass:integrated.bias} and $\sqrt n/m \to 0$, it has the same first-order Gaussian limit as the empirical copula-based estimator. This regime includes the practical degree $m = \lfloor n^{2/3}\rfloor$.

A smoothed beta bootstrap based on the empirical beta copula was justified for pointwise confidence intervals. The interval in Algorithm~\ref{alg:bootstrap.ci} is obtained by inverting the absolute centered bootstrap root and is therefore symmetric around the point estimate; it is not a percentile interval. Pointwise coverage curves and Wilson Monte Carlo intervals show how interval performance changes with $p$, with the values at $p = 0.05$ and $p = 0.10$ reported separately. The equal-tailed basic construction provides a sensitivity check for bootstrap-root asymmetry.

The numerical study shows that Bernstein smoothing reduces integrated variance in every reported setting and often lowers MISE, especially under weak to moderate dependence and for small or moderate sample sizes. However, the improvement is not uniform: under stronger positive dependence, and in some strong negative cases, the smoothing bias may become non-negligible and can offset the variance reduction. The direct common-sample comparison with the empirical beta and ECBC estimators reports target-specific accuracy together with adaptive degrees, computation time, convergence diagnostics, and numerical failures. The separate exploratory common-bootstrap experiment uses a local-grid approximation to ECBC degree reselection and compares pointwise coverage over the threshold grid. Its small number of Monte Carlo and bootstrap replications limits its numerical precision. Lower-boundary undercoverage also occurs in several FGM settings, while the pronounced failures under strong positive dependence give a clear warning that smoothing bias can compromise interval calibration. These comparisons make the practical trade-off explicit: the proposed estimator is simple and directly evaluates a particular bivariate curve, the empirical beta estimator is tuning-free, and ECBC offers adaptive full-copula estimation at a higher computational cost. No procedure is claimed to dominate uniformly.

The Loss--ALAE application uses the empirical-beta generator with reproducible tie handling and reports symmetric and equal-tailed interval constructions. The curve estimators and global benchmark use the same tie-broken ranks, and a point-estimate sensitivity analysis records the variation across tie resolutions both with and without the observations identified as policy-limit capped. The latter omission is not a correction for censoring. Since the available theory assumes continuous margins, the application intervals remain conditional on the selected tie resolution. Subject to these qualifications, the estimated lower-tail Spearman's rho curve reveals local dependence features that are not captured by the global rank correlation.

Several directions remain open for future research. A lower-tail Kendall-type curve would be a natural companion, but its nonlinear copula functional requires separate estimation and asymptotic arguments. Data-adaptive degree selection is another important extension; it would require uniform theory over candidate degrees and reselection inside the bootstrap. It would also be useful to study simultaneous confidence bands, behavior as $p \downarrow 0$, censored or dependent data, semi-continuous outcomes with ties or mass points \citep{Arends2025rank}, and higher-dimensional local concordance functionals. Finally, triangular-array methods \citep{Lu2021BCDF,Wang2022Bdensity} may provide useful tools for refined Bernstein approximations, although their adaptation to rank-based copula functionals and smooth resampling remains to be developed.

\section*{Funding}
\addcontentsline{toc}{section}{Funding}

Fr\'ed\'eric Ouimet is supported by the Natural Sciences and Engineering Research Council of Canada (NSERC) through Discovery Grant and Launch Supplement (RGPIN-2026-04471, DGECR-2026-00449).

\section*{Data availability}
\addcontentsline{toc}{section}{Data availability}

The insurance data used in Section~\ref{sec:application} are available in the \textsf{R} package \textbf{evd} as the \texttt{lossalae} dataset. The \textsf{R} code for the simulation study is available in the GitHub repository at \url{https://github.com/FredericOuimetMcGill/Spearman}.

\section*{Disclosure statement}
\addcontentsline{toc}{section}{Disclosure statement}

The authors declare no conflicts of interest.

\section*{References}
\addcontentsline{toc}{section}{References}

\setlength{\bibsep}{0pt plus 0ex} 

\bibliographystyle{plainnat}
\bibliography{bib}

\appendix

\begin{appendices}

\renewcommand{\thefigure}{\arabic{figure}}
\renewcommand{\thetable}{\arabic{table}}

\renewcommand{\thesection}{Appendix~\Alph{section}}
\section{Proof of the main results}\label{app:proof}
\renewcommand{\thesection}{\Alph{section}}

\begin{proof}[Proof of Theorem~\ref{thm:strong.consistency}]
The Bernstein operator is positive, linear, and a contraction in the supremum norm. Since $C_{m,n} = B_m(C_n)$,
\[
\|C_{m,n} - C\|_{\infty,[0,1]^2} \leq \|B_m(C_n - C)\|_{\infty,[0,1]^2} + \|B_m(C) - C\|_{\infty,[0,1]^2} \leq \|C_n - C\|_{\infty,[0,1]^2} + \|B_m(C) - C\|_{\infty,[0,1]^2}.
\]
The empirical copula is uniformly strongly consistent, so the first term converges to zero almost surely. Since every copula is continuous on $[0,1]^2$ and $m \to \infty$, the multivariate Bernstein approximation theorem implies that the second term converges to zero. Hence $\|C_{m,n} - C\|_{\infty,[0,1]^2} \to 0$ almost surely. Moreover,
\[
\begin{aligned}
\sup_{p \in \mathcal{P}}|\widehat{\rho}_{m,n}(p) - \rho(p)| &\leq \sup_{p \in \mathcal{P}}\frac{1}{D(p)}\int_{[0,p]^2}|C_{m,n}(u,v) - C(u,v)| \, \rd u \, \rd v \\
&\leq \frac{1}{D(p_0)}\|C_{m,n} - C\|_{\infty,[0,1]^2} \to 0, \qquad \text{a.s.}
\end{aligned}
\]
The fixed-threshold conclusion follows by taking any compact interval $[p_0,1]$ containing the specified value of $p$.
\end{proof}

\begin{proof}[Proof of Proposition~\ref{prop:bias.sufficient}]
Let $L < \infty$ be a Lipschitz constant for the gradient of $C$. For fixed $(u,v) \in [0,1]^2$, let $S$ and $T$ be independent binomial random variables with parameters $(r,u)$ and $(r,v)$, respectively, and set $\mathbf X_r = (S/r,T/r)$ and $\mathbf x = (u,v)$. The first-order Taylor formula with Lipschitz remainder gives
\[
|C(\mathbf X_r) - C(\mathbf x) - \nabla C(\mathbf x)^{\mathsf T}(\mathbf X_r - \mathbf x)| \leq \frac{L}{2}\|\mathbf X_r - \mathbf x\|_2^2.
\]
Taking expectations, using $B_r(C)(u,v) = \EE\{C(\mathbf X_r)\}$, $\EE(\mathbf X_r - \mathbf x) = \mathbf 0$, and $\EE\|\mathbf X_r - \mathbf x\|_2^2 \leq 1/(2r)$, we obtain
\[
|B_r(C)(u,v) - C(u,v)| \leq \frac{L}{4r}.
\]
Consequently,
\[
\sup_{p \in \mathcal{P}}\left|\int_{[0,p]^2}\{B_r(C)(u,v) - C(u,v)\} \, \rd u \, \rd v\right| \leq \frac{L}{4r},
\]
which proves Assumption~\ref{ass:integrated.bias}.
\end{proof}

\begin{proof}[Proof of Theorem~\ref{thm:functional_clt_rho}]
For $p \in \mathcal{P}$, let $g_p(u,v) = \ind\{u \leq p,v \leq p\}$, $(u,v) \in [0,1]^2$. Let $\mathcal{B}_\infty([0,1]^2)$ denote the subspace of $\ell^\infty([0,1]^2)$ consisting of bounded Borel measurable functions, and define $\Phi:\mathcal{B}_\infty([0,1]^2) \to \ell^\infty(\mathcal{P})$ by
\[
(\Phi(h))(p) = \frac{1}{D(p)}\int_{[0,1]^2}g_p(u,v)h(u,v) \, \rd u \, \rd v, \qquad p \in \mathcal{P}.
\]
The map $\Phi$ is linear and
\[
\|\Phi(h)\|_{\infty,\mathcal{P}} \leq \frac{1}{D(p_0)}\|h\|_{\infty,[0,1]^2},
\]
so it is bounded and continuous.

Write $\mathbb{C}_n = \sqrt n(C_n - C)$. Since $C_{m,n} = B_m(C_n)$,
\[
\mathbb{Z}_n = \Phi\{B_m(\mathbb{C}_n)\} + \sqrt n\Phi\{B_m(C) - C\}.
\]
Under Assumption~\ref{ass:segers}, $\mathbb{C}_n \rightsquigarrow \mathbb{G}_C$ in $\ell^\infty([0,1]^2)$, where $\mathbb{G}_C$ has continuous sample paths. For a bounded function $h$ on $[0,1]^2$, let $\omega(h,\delta)$ denote its modulus of continuity with respect to the supremum metric. If $S$ and $T$ are independent binomial random variables with parameters $(m,u)$ and $(m,v)$, respectively, then Chebyshev's inequality gives, for every $\delta > 0$,
\[
\|B_m(h) - h\|_{\infty,[0,1]^2} \leq \omega(h,\delta) + \frac{\|h\|_{\infty,[0,1]^2}}{m\delta^2}.
\]
Take $\delta = m^{-1/4}$. Since $m \to \infty$, weak convergence to a process with continuous sample paths implies asymptotic uniform equicontinuity of $\mathbb{C}_n$ and $\|\mathbb{C}_n\|_{\infty,[0,1]^2} = \OO_{\PP}(1)$. Consequently,
\[
\|B_m(\mathbb{C}_n) - \mathbb{C}_n\|_{\infty,[0,1]^2} = \oo_{\PP}(1).
\]
It follows that
\[
\Phi\{B_m(\mathbb{C}_n)\} \rightsquigarrow \Phi(\mathbb{G}_C) \qquad \text{in }\ell^\infty(\mathcal{P}).
\]
By Assumption~\ref{ass:integrated.bias},
\[
\begin{aligned}
\left\|\sqrt n\Phi\{B_m(C) - C\}\right\|_{\infty,\mathcal{P}} &\leq \frac{\sqrt n}{D(p_0)}\sup_{p \in \mathcal{P}}\left|\int_{[0,p]^2}\{B_m(C)(u,v) - C(u,v)\} \, \rd u \, \rd v\right| \\
&= \OO(\sqrt n/m) = \oo(1).
\end{aligned}
\]
Slutsky's theorem therefore gives $\mathbb{Z}_n \rightsquigarrow \Phi(\mathbb{G}_C)$. The expression for the limit follows from the definition of $\Phi$. Since $\mathbb{Z}$ is a linear transform of the centered Gaussian process $\mathbb{G}_C$, it is centered Gaussian. Finally, Fubini's theorem yields
\[
\Cov\{\mathbb{Z}(p),\mathbb{Z}(q)\} = \frac{1}{D(p)D(q)}\int_{[0,p]^2}\int_{[0,q]^2}\Gamma((u,v),(s,t)) \, \rd u \, \rd v \, \rd s \, \rd t.
\]
This completes the proof.
\end{proof}

\begin{proof}[Proof of Theorem~\ref{thm:bootstrap.validity}]
Let $C_n^\#$ be the rank-based empirical copula computed from a bootstrap sample drawn from the empirical beta copula. By \citet[Proposition~3.3]{Kiriliouk2021resampling}, under Assumption~\ref{ass:segers}, conditionally on the data and in probability,
\[
\mathbb{C}_n^\# = \sqrt n(C_n^\# - C_n) \rightsquigarrow \mathbb{G}_C \qquad \text{in }\ell^\infty([0,1]^2).
\]
Because the same Bernstein degree is used for the original and bootstrap estimators,
\[
\mathbb{Z}_n^\# = \Phi\{B_m(\mathbb{C}_n^\#)\}.
\]
The deterministic Bernstein bound used in the proof of Theorem~\ref{thm:functional_clt_rho}, together with conditional asymptotic uniform equicontinuity and tightness of $\mathbb{C}_n^\#$, implies that, for every $\varepsilon > 0$,
\[
\PP^\#\left\{\|B_m(\mathbb{C}_n^\#) - \mathbb{C}_n^\#\|_{\infty,[0,1]^2} > \varepsilon\right\} \to 0
\]
in probability whenever $m \to \infty$. Since $\Phi$ is bounded and linear, the conditional continuous mapping theorem yields
\[
\mathbb{Z}_n^\# \rightsquigarrow \Phi(\mathbb{G}_C) = \mathbb{Z}
\]
conditionally on the data, in probability, in $\ell^\infty(\mathcal{P})$. The additional degree and bias conditions in the statement ensure that the original process $\mathbb{Z}_n$ has the same limit by Theorem~\ref{thm:functional_clt_rho}.
\end{proof}

\begin{proof}[Proof of Corollary~\ref{cor:pointwise.coverage}]
The evaluation map $h \mapsto h(p_j)$ is continuous on $\ell^\infty(\mathcal{P})$. Hence Theorems~\ref{thm:functional_clt_rho} and~\ref{thm:bootstrap.validity} imply that $\mathbb{Z}_n(p_j)$ converges weakly to $\mathbb{Z}(p_j)$ and that the conditional law of $\mathbb{Z}_n^\#(p_j)$ converges in probability to the same limit. Since $\Var\{\mathbb{Z}(p_j)\} > 0$, the distribution of $|\mathbb{Z}(p_j)|$ is continuous and strictly increasing on $[0,\infty)$. Conditional quantile consistency therefore gives
\[
\PP\{|\mathbb{Z}_n(p_j)| \leq c_{j,1-\alpha,n}^\#\} \to 1 - \alpha.
\]
The event inside the probability is equivalent to coverage by the interval in the statement. Finally, if $B = B_n \to \infty$, the empirical distribution of the independent bootstrap roots converges uniformly to their conditional distribution, while the latter converges to the continuous distribution of $|\mathbb{Z}(p_j)|$. Hence the empirical bootstrap quantile may replace $c_{j,1-\alpha,n}^\#$ without changing the coverage limit.
\end{proof}

\renewcommand{\thesection}{Appendix~\Alph{section}}
\section{Additional simulation results}\label{app:simulation}
\renewcommand{\thesection}{\Alph{section}}

\begin{figure}[H]
\centering
\includegraphics[trim = 3mm 2mm 7mm 6mm, clip, width = 0.4\textwidth]{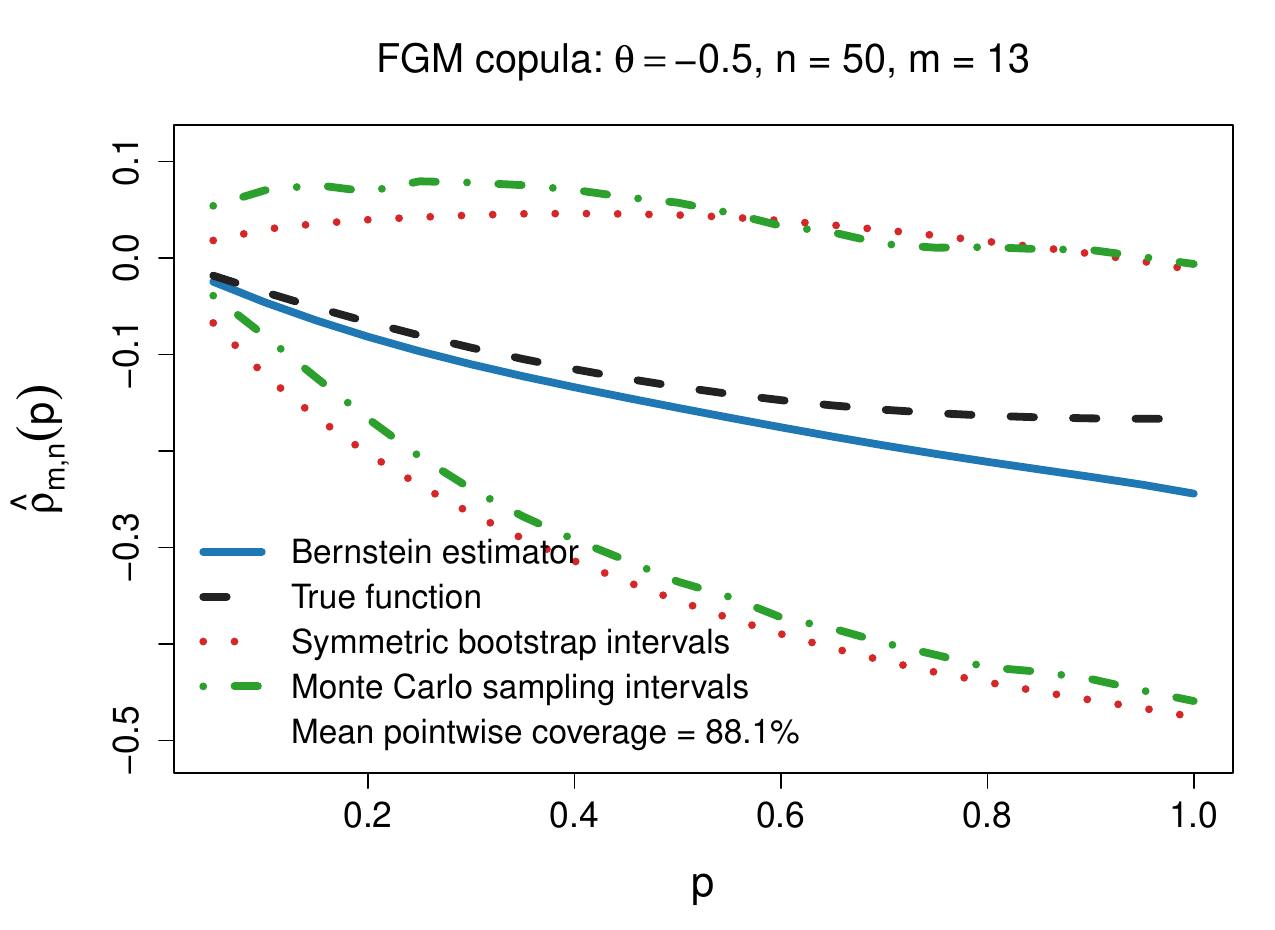}
\quad
\includegraphics[trim = 3mm 2mm 7mm 6mm, clip, width = 0.4\textwidth]{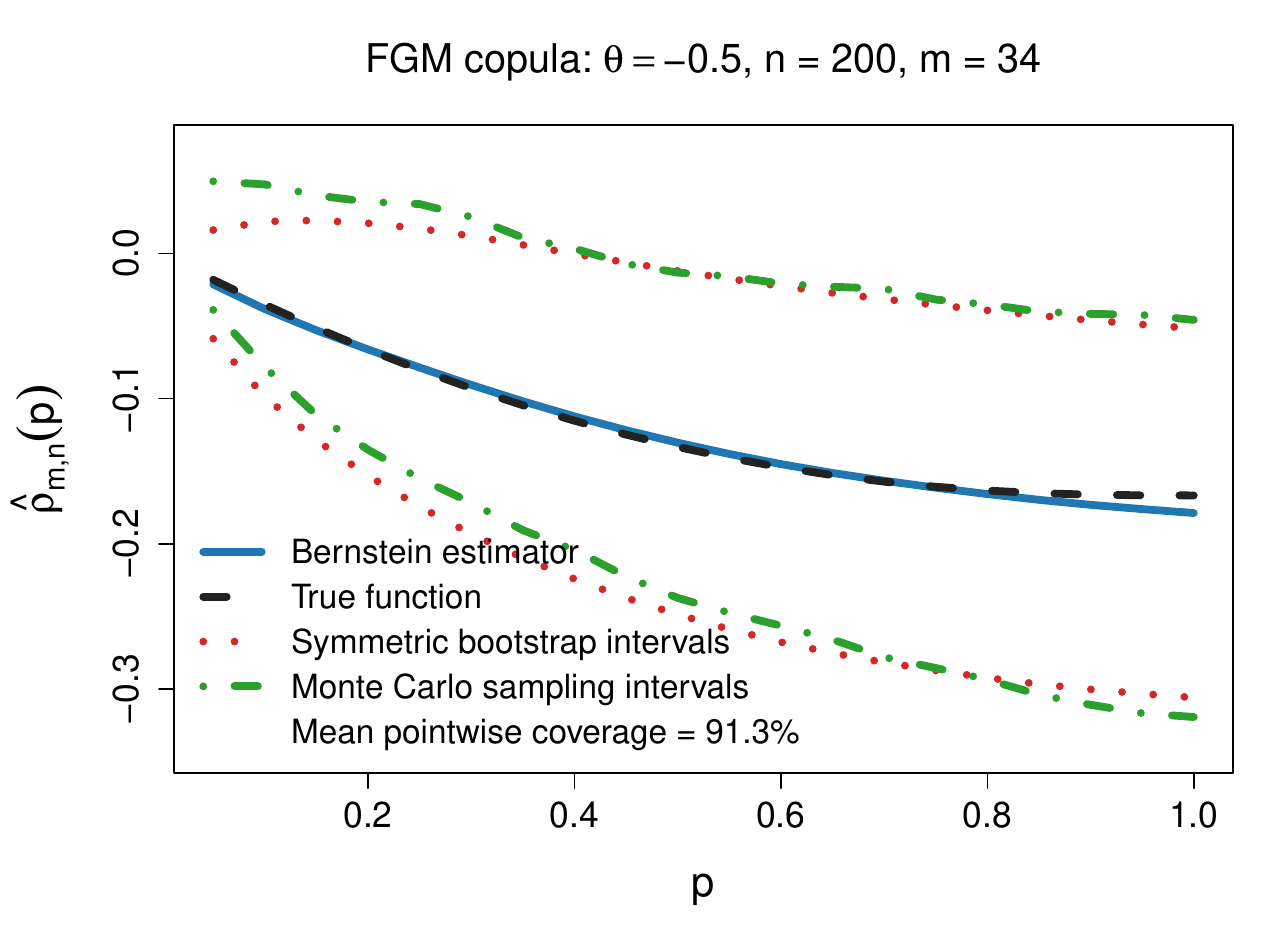} \\
\includegraphics[trim = 3mm 2mm 7mm 6mm, clip, width = 0.4\textwidth]{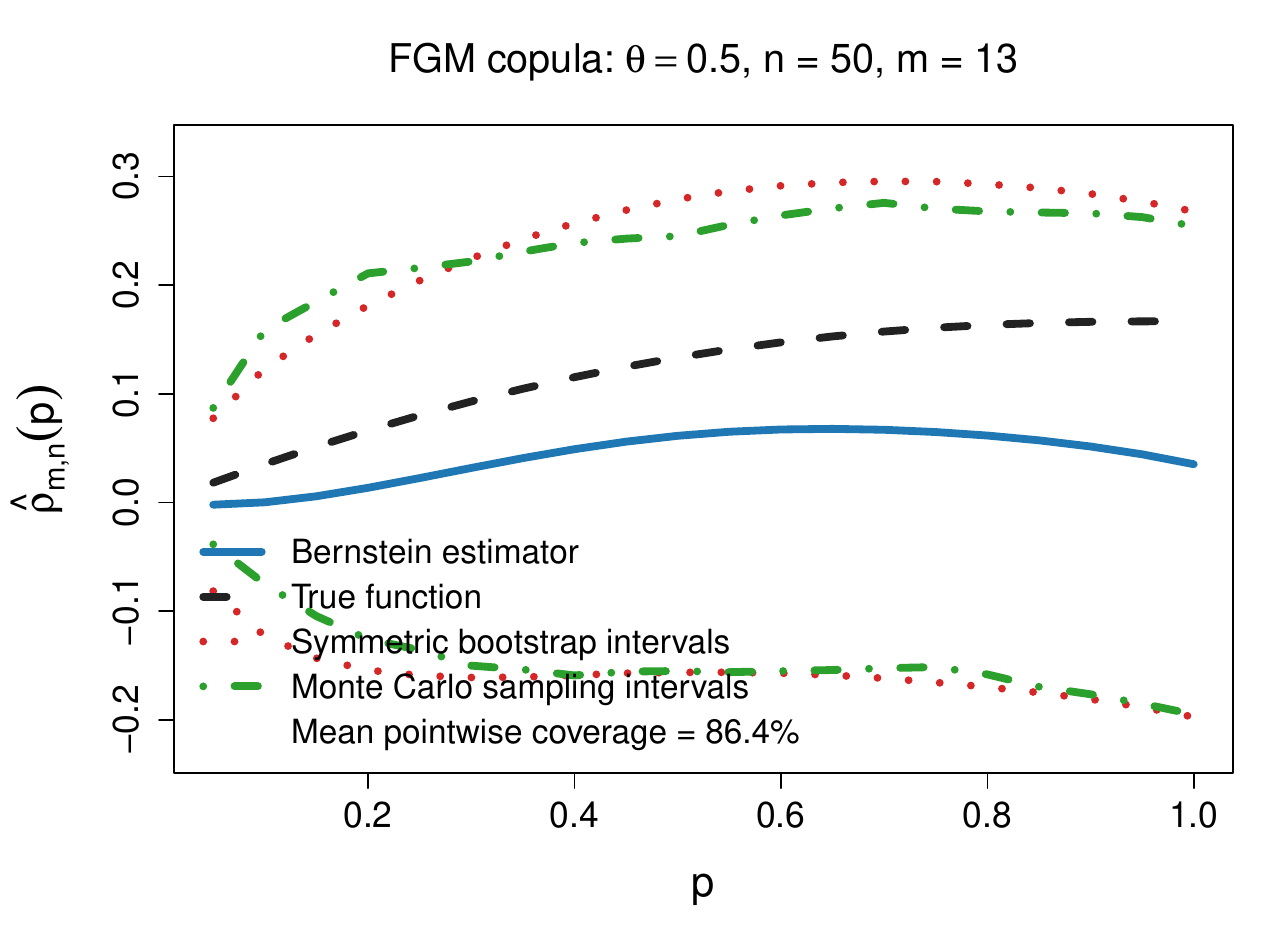}
\quad
\includegraphics[trim = 3mm 2mm 7mm 6mm, clip, width = 0.4\textwidth]{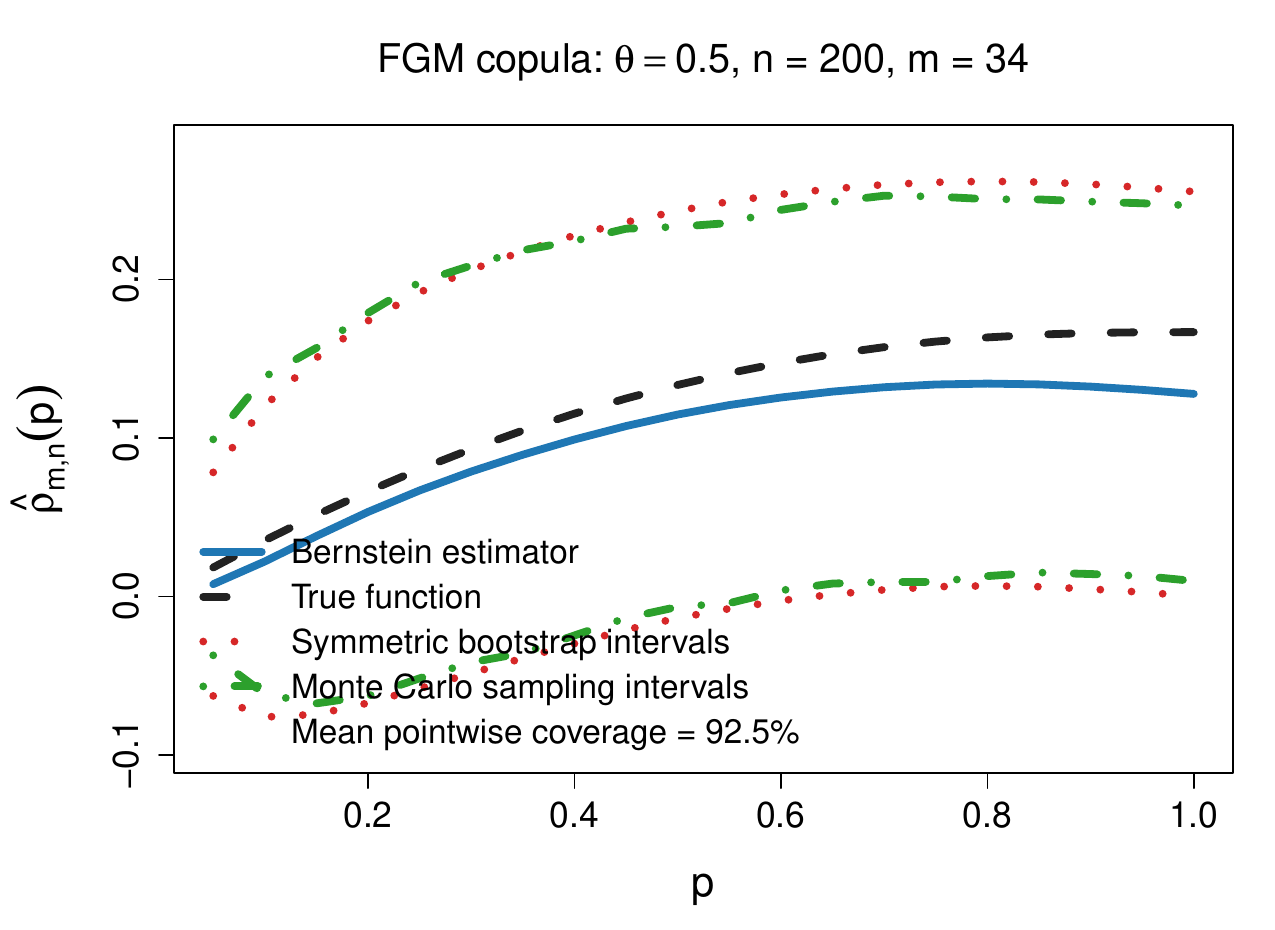} \\[-3mm]
\caption{Finite-sample assessment of the empirical-beta bootstrap pointwise confidence intervals under the FGM copula. Each panel displays the Monte Carlo average of the Bernstein-smoothed lower-tail Spearman's rho estimator $p \mapsto \widehat{\rho}_{m,n}(p)$ (solid blue), the true curve $p \mapsto \rho_\theta(p)$ (dashed black), the average lower and upper endpoints of the symmetric absolute-root intervals (red dotted), and the Monte Carlo sampling intervals obtained from the empirical pointwise quantiles of $\widehat{\rho}_{m,n}(p)$ across Monte Carlo replications (green dash-dotted). Rows correspond to $\theta = -0.5$ and $\theta = 0.5$, while columns correspond to $n = 50$ and $n = 200$, with $m = \lfloor n^{2/3}\rfloor$ and $m_\# = n$. The mean pointwise coverage reported in each panel is the empirical coverage averaged over the threshold grid and is included as a secondary summary.}
\label{fig:simulation_CI_app}
\end{figure}

\begin{figure}[H]
\centering
\includegraphics[trim = 3mm 2mm 4mm 5mm, clip, width = 0.99\textwidth]{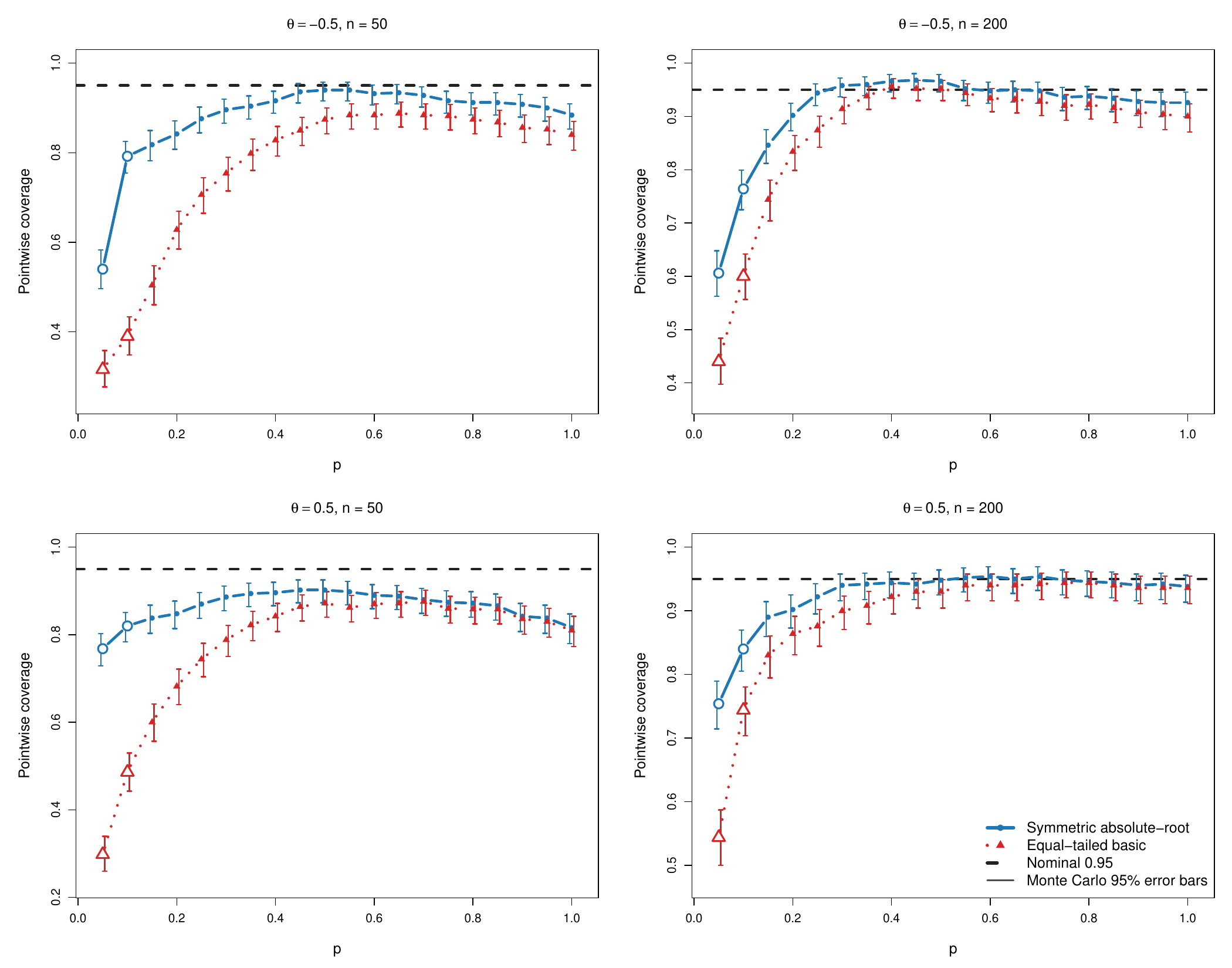}
\caption{Empirical pointwise coverage of nominal 95\% confidence intervals under the FGM copula. Rows correspond to the intermediate dependence levels $\theta = -0.5$ and $\theta = 0.5$, and columns correspond to $n = 50$ and $n = 200$. The blue solid curves represent the symmetric absolute-root intervals, the red dotted curves represent the equal-tailed basic intervals, and the horizontal dashed line gives the nominal coverage level. Vertical bars are 95\% Wilson Monte Carlo intervals, and the values at $p = 0.05$ and $p = 0.10$ are highlighted.}
\label{fig:fgm_pointwise_coverage_app}
\end{figure}

\vspace{-3mm}
\begin{table}[H]
\centering
\small
\setlength{\tabcolsep}{7pt}
\renewcommand{\arraystretch}{0.9}
\caption{Computation time per sample for the four estimators in the common-sample comparison. Entries give the median time in milliseconds, with the 95th percentile in parentheses, pooled over the copula families and dependence settings and calculated from successful evaluations. The ECBC time includes empirical-Bayes degree selection and posterior diagnostics.}
\label{tab:alternative_runtime}
\begin{tabular}{lrr}
\toprule[1.5pt]
Estimator & $n = 50$ & $n = 200$ \\
\midrule[1pt]
Empirical copula & $0.460$ ($0.600$) & $0.500$ ($0.600$) \\
Fixed-degree Bernstein & $0.700$ ($0.900$) & $3.80$ ($5.00$) \\
Empirical beta & $5.00$ ($6.00$) & $14.00$ ($16.00$) \\
ECBC & $330.0$ ($386.5$) & $1860$ ($2158$) \\
\bottomrule[1.5pt]
\end{tabular}
\end{table}

\vspace{-3mm}
\begin{table}[H]
\centering
\small
\setlength{\tabcolsep}{4pt}
\renewcommand{\arraystretch}{0.9}
\caption{Numerical diagnostics for the ECBC point-estimator comparison, pooled over the copula families and dependence settings. The selector uses the marginal posterior modes from two MCMC chains. Fits with an $\widehat R$ value above $1.20$ are classified as diagnostically nonconverged and are excluded from all pointwise and integrated accuracy summaries.}
\label{tab:alternative_ecbc_diagnostics}
\begin{tabular}{crrrrrrr}
\toprule[1.5pt]
$n$ & Attempted & Numerical failures & Nonconverged & Retained & Median max. $\widehat R$ & Max. $\widehat R$ & Mean acceptance \\
\midrule[1pt]
$50$ & 2000 & 0 & 0 & 2000 & 1.008 & 1.106 & 0.705 \\
$200$ & 2000 & 0 & 4 & 1996 & 1.019 & 1.266 & 0.738 \\
\bottomrule[1.5pt]
\end{tabular}
\end{table}

\begin{table}[H]
\centering
\small
\setlength{\tabcolsep}{5pt}
\renewcommand{\arraystretch}{0.85}
\caption{MISE decomposition for the empirical copula-based estimator and the Bernstein estimator when $n = 50$. Entries for integrated squared bias, integrated variance, and MISE are reported as $10^3$ times the unscaled values. The final column reports $100\times \MISE_{B}/\MISE_{E}$, where values below $100\%$ indicate lower MISE for the Bernstein estimator. For the FGM copula, the second column gives the copula parameter $\theta$; for the Gaussian, Clayton, and Frank copulas, it gives Kendall's $\tau$. The Bernstein degree is $m = 13$.}
\label{tab:mise_decomp_n50}
\begin{tabular}{lcrrrrrrr}
\toprule[1.5pt]
& & \multicolumn{3}{c}{Empirical estimator} & \multicolumn{3}{c}{Bernstein estimator} & \\
\cmidrule(lr){3-5}\cmidrule(lr){6-8}
Copula & $\theta$ or $\tau$
& $\IBias^2$ & $\IVar$ & $\MISE$
& $\IBias^2$ & $\IVar$ & $\MISE$
& Ratio \\
\midrule[1pt]
\multirow{5}{*}{FGM}
& $-1$ & 2.089 & 9.267 & 11.356 & 0.244 & 6.339 & 6.583 & 58.0\% \\
& $-0.5$ & 3.624 & 13.653 & 17.278 & 1.410 & 8.978 & 10.388 & 60.1\% \\
& $0$ & 3.003 & 16.950 & 19.953 & 2.596 & 10.517 & 13.113 & 65.7\% \\
& $0.5$ & 4.040 & 19.445 & 23.485 & 5.907 & 11.551 & 17.458 & 74.3\% \\
& $1$ & 6.174 & 17.276 & 23.451 & 11.784 & 9.699 & 21.483 & 91.6\% \\
\midrule
\multirow{5}{*}{Gaussian}
& $-0.5$ & 1.436 & 3.635 & 5.071 & 1.124 & 2.626 & 3.750 & 73.9\% \\
& $-0.25$ & 1.709 & 10.828 & 12.538 & 0.166 & 7.392 & 7.558 & 60.3\% \\
& $0$ & 2.388 & 16.544 & 18.932 & 2.132 & 10.149 & 12.281 & 64.9\% \\
& $0.25$ & 7.076 & 20.276 & 27.352 & 15.846 & 10.932 & 26.778 & 97.9\% \\
& $0.5$ & 13.343 & 13.878 & 27.221 & 45.712 & 5.784 & 51.495 & 189.2\% \\
\midrule
\multirow{5}{*}{Clayton}
& $-0.5$ & 1.561 & 3.233 & 4.794 & 1.845 & 2.468 & 4.313 & 90.0\% \\
& $-0.25$ & 1.865 & 7.407 & 9.272 & 0.931 & 5.495 & 6.426 & 69.3\% \\
& $0$ & 3.302 & 17.392 & 20.694 & 2.934 & 10.787 & 13.721 & 66.3\% \\
& $0.25$ & 12.054 & 22.755 & 34.809 & 34.574 & 11.371 & 45.945 & 132.0\% \\
& $0.5$ & 23.127 & 9.615 & 32.742 & 88.672 & 4.553 & 93.226 & 284.7\% \\
\midrule
\multirow{5}{*}{Frank}
& $-0.5$ & 1.652 & 3.550 & 5.202 & 1.121 & 2.508 & 3.629 & 69.8\% \\
& $-0.25$ & 1.763 & 10.448 & 12.211 & 0.118 & 7.046 & 7.164 & 58.7\% \\
& $0$ & 2.919 & 15.571 & 18.490 & 2.523 & 9.673 & 12.196 & 66.0\% \\
& $0.25$ & 6.924 & 17.165 & 24.089 & 14.350 & 9.817 & 24.167 & 100.3\% \\
& $0.5$ & 10.981 & 13.612 & 24.593 & 35.095 & 5.926 & 41.021 & 166.8\% \\
\bottomrule[1.5pt]
\end{tabular}
\end{table}

\begin{table}[H]
\centering
\small
\setlength{\tabcolsep}{5pt}
\renewcommand{\arraystretch}{0.85}
\caption{MISE decomposition for the empirical copula-based estimator and the Bernstein estimator when $n = 200$. Entries for integrated squared bias, integrated variance, and MISE are reported as $10^3$ times the unscaled values. The final column reports $100\times \MISE_{B}/\MISE_{E}$, where values below $100\%$ indicate lower MISE for the Bernstein estimator. For the FGM copula, the second column gives the copula parameter $\theta$; for the Gaussian, Clayton, and Frank copulas, it gives Kendall's $\tau$. The Bernstein degree is $m = 34$.}
\label{tab:mise_decomp_n200}
\begin{tabular}{lcrrrrrrr}
\toprule[1.5pt]
& & \multicolumn{3}{c}{Empirical estimator} & \multicolumn{3}{c}{Bernstein estimator} & \\
\cmidrule(lr){3-5}\cmidrule(lr){6-8}
Copula & $\theta$ or $\tau$
& $\IBias^2$ & $\IVar$ & $\MISE$
& $\IBias^2$ & $\IVar$ & $\MISE$
& Ratio \\
\midrule[1pt]
\multirow{5}{*}{FGM}
& $-1$ & 0.093 & 2.306 & 2.399 & 0.063 & 1.960 & 2.023 & 84.3\% \\
& $-0.5$ & 0.091 & 4.022 & 4.112 & 0.012 & 3.301 & 3.312 & 80.5\% \\
& $0$ & 0.278 & 4.560 & 4.838 & 0.250 & 3.636 & 3.885 & 80.3\% \\
& $0.5$ & 0.610 & 4.897 & 5.507 & 0.956 & 3.835 & 4.791 & 87.0\% \\
& $1$ & 0.337 & 4.721 & 5.058 & 0.992 & 3.555 & 4.547 & 89.9\% \\
\midrule
\multirow{5}{*}{Gaussian}
& $-0.5$ & 0.083 & 0.954 & 1.037 & 0.303 & 0.834 & 1.137 & 109.7\% \\
& $-0.25$ & 0.246 & 2.729 & 2.975 & 0.013 & 2.352 & 2.364 & 79.5\% \\
& $0$ & 0.161 & 4.304 & 4.465 & 0.150 & 3.401 & 3.551 & 79.5\% \\
& $0.25$ & 0.513 & 5.868 & 6.381 & 1.930 & 4.085 & 6.014 & 94.3\% \\
& $0.5$ & 1.044 & 4.523 & 5.566 & 7.310 & 2.720 & 10.031 & 180.2\% \\
\midrule
\multirow{5}{*}{Clayton}
& $-0.5$ & 0.092 & 0.773 & 0.865 & 0.458 & 0.695 & 1.153 & 133.2\% \\
& $-0.25$ & 0.089 & 1.755 & 1.844 & 0.211 & 1.555 & 1.766 & 95.8\% \\
& $0$ & 0.264 & 4.656 & 4.920 & 0.241 & 3.738 & 3.979 & 80.9\% \\
& $0.25$ & 1.055 & 6.628 & 7.683 & 6.714 & 4.543 & 11.257 & 146.5\% \\
& $0.5$ & 1.691 & 2.486 & 4.177 & 18.242 & 1.600 & 19.842 & 475.1\% \\
\midrule
\multirow{5}{*}{Frank}
& $-0.5$ & 0.096 & 0.925 & 1.021 & 0.328 & 0.796 & 1.123 & 110.0\% \\
& $-0.25$ & 0.106 & 2.825 & 2.931 & 0.044 & 2.358 & 2.402 & 81.9\% \\
& $0$ & 0.220 & 4.673 & 4.893 & 0.206 & 3.754 & 3.960 & 80.9\% \\
& $0.25$ & 0.313 & 5.022 & 5.336 & 1.212 & 3.711 & 4.923 & 92.3\% \\
& $0.5$ & 0.618 & 3.909 & 4.527 & 4.151 & 2.440 & 6.592 & 145.6\% \\
\bottomrule[1.5pt]
\end{tabular}
\end{table}

\begin{table}[H]
\centering
\small
\setlength{\tabcolsep}{4pt}
\renewcommand{\arraystretch}{0.95}
\caption{Integrated mean squared error for the four estimators when $n = 50$. MISE entries are reported as $10^3$ times the unscaled values. The final columns summarize the ECBC empirical-Bayes degrees.}
\label{tab:alternative_comparison_n50}
\begin{tabular}{lcrrrrrr}
\toprule[1.5pt]
Copula & $\theta$ or $\tau$ & Empirical & Fixed degree & Beta & ECBC & $\widetilde m_1$ & $\widetilde m_2$ \\
\midrule[1pt]
FGM & $-1$ & 10.463 & 6.009 & 8.807 & 8.267 & 15.0 & 15.0 \\
FGM & $-0.5$ & 17.997 & 10.824 & 12.905 & 10.211 & 16.0 & 16.0 \\
FGM & $0$ & 20.294 & 13.530 & 16.267 & 12.491 & 16.0 & 16.0 \\
FGM & $0.5$ & 22.445 & 15.170 & 21.012 & 15.138 & 16.0 & 16.0 \\
FGM & $1$ & 21.431 & 21.862 & 14.326 & 11.633 & 16.0 & 15.0 \\
\midrule
Gaussian & $-0.5$ & 4.784 & 4.015 & 4.208 & 6.914 & 15.0 & 15.0 \\
Gaussian & $-0.25$ & 15.512 & 8.680 & 11.692 & 10.017 & 16.0 & 16.0 \\
Gaussian & $0$ & 18.407 & 11.395 & 16.346 & 12.234 & 16.0 & 16.0 \\
Gaussian & $0.25$ & 27.833 & 26.020 & 21.677 & 16.586 & 15.0 & 16.0 \\
Gaussian & $0.5$ & 24.704 & 48.252 & 14.577 & 20.675 & 15.0 & 15.0 \\
\midrule
Clayton & $-0.5$ & 4.403 & 4.498 & 4.000 & 6.801 & 17.0 & 16.0 \\
Clayton & $-0.25$ & 9.035 & 5.829 & 6.686 & 7.199 & 16.0 & 16.0 \\
Clayton & $0$ & 18.216 & 12.024 & 16.205 & 12.622 & 16.0 & 16.0 \\
Clayton & $0.25$ & 31.388 & 44.124 & 22.103 & 24.637 & 16.0 & 16.0 \\
Clayton & $0.5$ & 32.950 & 93.291 & 16.466 & 43.931 & 15.0 & 15.5 \\
\midrule
Frank & $-0.5$ & 4.927 & 3.295 & 3.630 & 5.984 & 15.0 & 15.0 \\
Frank & $-0.25$ & 12.858 & 7.779 & 11.566 & 10.093 & 16.0 & 16.0 \\
Frank & $0$ & 17.895 & 11.630 & 15.397 & 11.467 & 16.0 & 16.0 \\
Frank & $0.25$ & 23.784 & 24.425 & 16.809 & 14.128 & 15.0 & 16.0 \\
Frank & $0.5$ & 26.630 & 41.399 & 16.482 & 17.737 & 15.0 & 15.0 \\
\bottomrule[1.5pt]
\end{tabular}
\end{table}

\bigskip
\begin{table}[H]
\centering
\small
\setlength{\tabcolsep}{4pt}
\renewcommand{\arraystretch}{0.95}
\caption{Integrated mean squared error for the four estimators when $n = 200$. MISE entries are reported as $10^3$ times the unscaled values. The final columns summarize the ECBC empirical-Bayes degrees.}
\label{tab:alternative_comparison_n200}
\begin{tabular}{lcrrrrrr}
\toprule[1.5pt]
Copula & $\theta$ or $\tau$ & Empirical & Fixed degree & Beta & ECBC & $\widetilde m_1$ & $\widetilde m_2$ \\
\midrule[1pt]
FGM & $-1$ & 2.437 & 2.145 & 2.416 & 2.292 & 48.0 & 49.0 \\
FGM & $-0.5$ & 4.525 & 3.709 & 4.492 & 4.167 & 48.0 & 49.0 \\
FGM & $0$ & 4.612 & 3.698 & 4.268 & 3.776 & 49.0 & 49.0 \\
FGM & $0.5$ & 5.269 & 4.702 & 4.654 & 4.172 & 49.0 & 49.0 \\
FGM & $1$ & 4.735 & 4.311 & 4.421 & 3.845 & 48.0 & 47.5 \\
\midrule
Gaussian & $-0.5$ & 1.289 & 1.212 & 1.149 & 1.367 & 44.0 & 44.0 \\
Gaussian & $-0.25$ & 3.152 & 2.272 & 2.638 & 2.339 & 48.0 & 48.0 \\
Gaussian & $0$ & 4.014 & 3.172 & 3.707 & 3.339 & 48.0 & 49.0 \\
Gaussian & $0.25$ & 5.558 & 5.463 & 4.972 & 4.437 & 47.0 & 47.0 \\
Gaussian & $0.5$ & 5.794 & 10.385 & 4.781 & 5.904 & 44.0 & 44.0 \\
\midrule
Clayton & $-0.5$ & 1.110 & 1.108 & 0.900 & 1.101 & 50.0 & 50.0 \\
Clayton & $-0.25$ & 2.224 & 2.399 & 2.424 & 2.664 & 47.0 & 48.0 \\
Clayton & $0$ & 5.071 & 4.172 & 4.949 & 4.457 & 49.5 & 49.0 \\
Clayton & $0.25$ & 8.947 & 13.474 & 7.616 & 9.001 & 48.0 & 48.0 \\
Clayton & $0.5$ & 3.962 & 19.863 & 2.898 & 8.978 & 45.0 & 45.0 \\
\midrule
Frank & $-0.5$ & 1.028 & 0.979 & 0.879 & 1.130 & 45.0 & 44.0 \\
Frank & $-0.25$ & 2.410 & 2.197 & 2.537 & 2.560 & 48.0 & 48.0 \\
Frank & $0$ & 5.330 & 4.347 & 5.135 & 4.570 & 49.0 & 48.5 \\
Frank & $0.25$ & 6.198 & 5.861 & 5.624 & 5.041 & 48.0 & 47.0 \\
Frank & $0.5$ & 5.368 & 7.508 & 4.630 & 4.676 & 45.0 & 44.0 \\
\bottomrule[1.5pt]
\end{tabular}
\end{table}

\begin{figure}[p]
\centering
\includegraphics[height = 0.92\textheight, keepaspectratio]{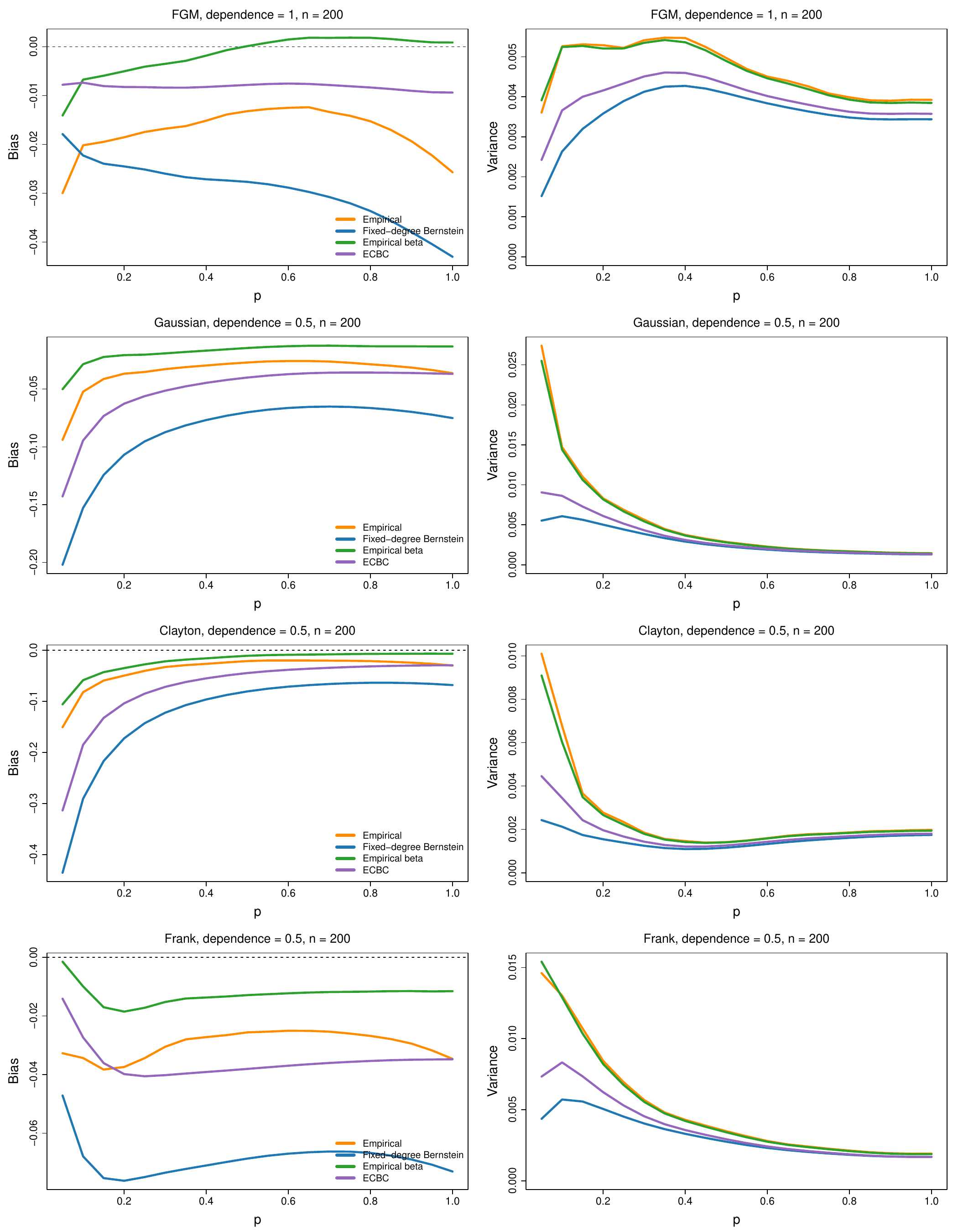}
\caption{Pointwise bias and variance for the empirical copula-based estimator, the fixed-degree Bernstein estimator, the empirical beta estimator, and the empirical checkerboard Bernstein copula estimator in the strongest positive-dependence setting considered for each copula family. The FGM row corresponds to $\theta = 1$, while the Gaussian, Clayton, and Frank rows correspond to $\tau = 0.5$. The left and right columns display bias and variance, respectively, for $n = 200$.}
\label{fig:alternative_pointwise}
\end{figure}

\begin{figure}[p]
\centering
\includegraphics[height = 0.92\textheight, keepaspectratio]{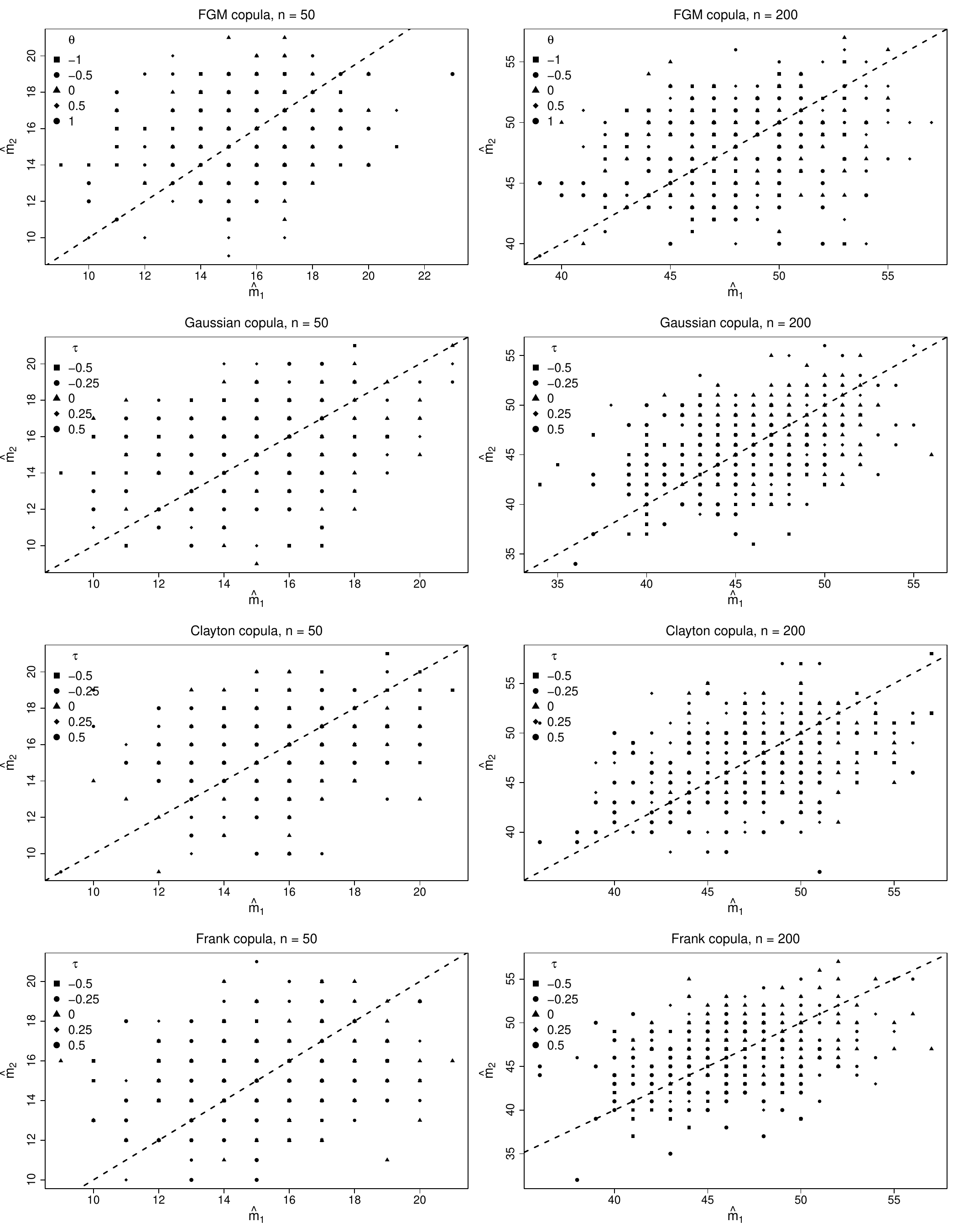}
\caption{Selected empirical checkerboard Bernstein copula degrees in the common-sample comparison. Each panel displays $(\widehat m_1,\widehat m_2)$, with the diagonal indicating equal marginal degrees. Rows correspond to copula families and columns correspond to sample sizes.}
\label{fig:alternative_degrees}
\end{figure}

\begin{figure}[p]
\centering
\includegraphics[height = 0.87\textheight, keepaspectratio]{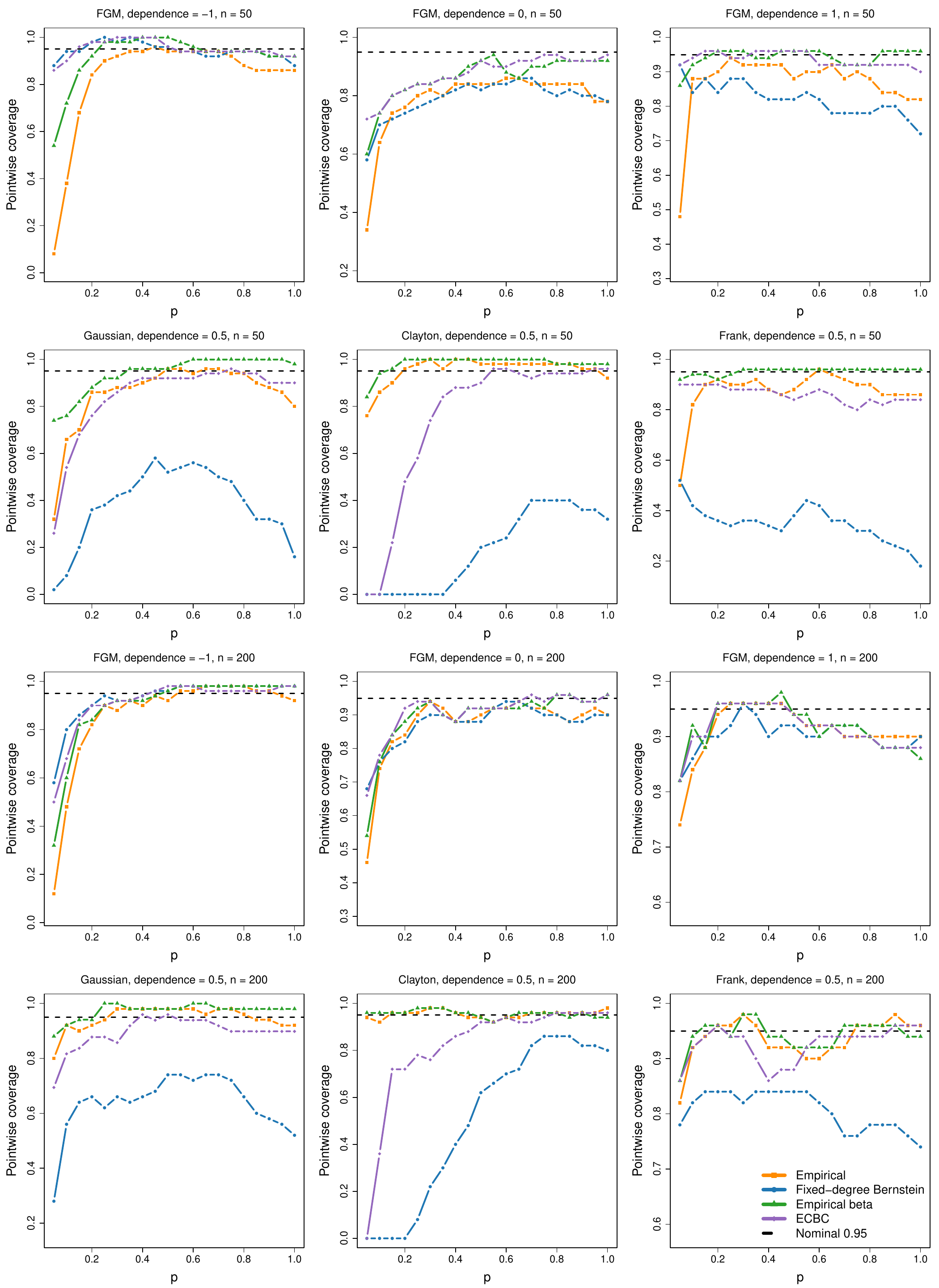}
\caption{Pointwise coverage of the symmetric absolute-root intervals in the common-bootstrap comparison. The first row displays the FGM cases $\theta = -1,0,1$ for $n = 50$, and the second row displays the Gaussian, Clayton, and Frank cases with $\tau = 0.5$ for $n = 50$. The third and fourth rows display the same respective cases for $n = 200$. Within each Monte Carlo replication, all available method-specific intervals are based on the same empirical-beta bootstrap samples. Whenever the observed ECBC fit is diagnostically converged, its degrees are the marginal posterior modes from the Markov chain Monte Carlo fit, and ECBC degree selection in each bootstrap sample is approximated by taking the joint posterior maximizer over an adaptively expanded local integer grid centered at the original degree pair, with initial radius four and at most two radius doublings after non-global boundary hits.}
\label{fig:alternative_coverage}
\end{figure}

\end{appendices}

\end{document}